\documentclass[12pt]{article}

\usepackage{amsmath,amsthm,amsfonts,amssymb,epsf,epsfig,psfrag,enumerate,hyperref}

\title{A Classification of Automorphisms\break of Compact 3-Manifolds}

\author{Leonardo Navarro Carvalho\footnote{Supported by FAPESP fellowship 03/06914-2, Brazil.} \hbox{\ \ }\&\hbox{\ } Ulrich Oertel}

\date{October, 2005}

\newtheorem{thm}{Theorem}[section] \newtheorem{lemma}[thm]{Lemma}
 \newtheorem{corollary}[thm]{Corollary}

\newtheorem{proposition}[thm]{Proposition}
\newtheorem{question}[thm]{Question} \newtheorem*{claim*}{Claim}
\newtheorem*{restate}{}

\newtheorem{conj}[thm]{Conjecture} 
 \theoremstyle{definition}
\newtheorem{defn}[thm]{Definition}
\newtheorem{defns}[thm]{Definitions} 
\newtheorem{com}[thm]{Remark} \newtheorem{ex}[thm]{Example}
 \newtheorem{remark}[thm]{Remark}
\theoremstyle{remark}

\def\figurefontsize{12}
\def\figuresmallfontsize{10}
\def\figuresmallerfontsize{6}

\begin{document}

\maketitle


\def\HDS{half-disc sum}

\def\irred{irreducible}
\def\half{spinal pair }
\def\spinal{\half}
\def\spinals{\halfs}
\def\halfs{spinal pairs }
\def\reals{\mathbb R}
\def\rationals{\mathbb Q}
\def\complex{\mathbb C}
\def\naturals{\mathbb N}
\def\integers{\mathbb Z}

\def\proj{P}
\def\hyp {\hbox {\rm {H \kern -2.8ex I}\kern 1.15ex}}

\def\weight#1#2#3{{#1}\raise2.5pt\hbox{$\centerdot$}\left({#2},{#3}\right)}
\def\intr{{\rm int}}
\def\inter{\ \raise4pt\hbox{$^\circ$}\kern -1.6ex}
\def\Cal{\cal}
\def\from{:}
\def\inverse{^{-1}}
\def\id{{\rm Id}}
\def\Max{{\rm Max}}
\def\Min{{\rm Min}}
\def\sing{{\rm Sing}}
\def\fr{{\rm fr}}
\def\embed{\hookrightarrow}
\def\Genus{{\rm Genus}}
\def\Z{Z}
\def\X{X}

\def\roster{\begin{enumerate}}
\def\endroster{\end{enumerate}}
\def\intersect{\cap}
\def\definition{\begin{defn}}
\def\enddefinition{\end{defn}}
\def\subhead{\subsection\{}
\def\theorem{thm}
\def\endsubhead{\}}
\def\head{\section\{}
\def\endhead{\}}
\def\example{\begin{ex}}
\def\endexample{\end{ex}}
\def\ves{\vs}
\def\mZ{{\mathbb Z}}
\def\M{M(\Phi)}
\def\bdry{\partial}
\def\hop{\vskip 0.15in}
\def\mathring{\inter}
\def\trip{\vskip 0.09in}

\begin{abstract}
We classify isotopy classes of automorphisms (self-homeomorphisms) of
3-manifolds satisfying the Thurston Geometrization Conjecture.  The
classification is similar to the classification of automorphisms of
surfaces developed by Nielsen and Thurston, except an automorphism of
a reducible manifold must first be written as a suitable composition
of two automorphisms, each of which fits into our classification.
Given an automorphism, the goal is to show, loosely speaking, either
that it is periodic, or that it can be decomposed on a surface
invariant up to isotopy, or that it  has a ``dynamically nice"
representative, with invariant laminations that ``fill" the manifold.

We consider automorphisms of irreducible and boundary-irreducible
3-mani\-folds as being already classified, though there are some
exceptional manifolds for which the automorphisms are not understood.
Thus the paper is particularly aimed at understanding automorphisms of
reducible and/or boundary reducible 3-manifolds.

Previously unknown phenomena are found even in the case of connected
sums of $S^2\times S^1$'s.  To deal with this case, we prove that a
minimal genus Heegaard decomposition is unique up to isotopy, a result
which apparently was previously unknown.

Much remains to be understood about some of the automorphisms of the
classification.

\end{abstract}

\section{Introduction}\label{Intro}

There is a large body of work studying the mapping class groups and
homeomorphism spaces of reducible and $\bdry$-reducible manifolds.
Without giving references, some of the researchers who have
contributed are:  E. C{\'e}sar de S{\'a}, D. Gabai, M.  Hamstrom,
A. Hatcher, H. Hendriks, K. Johannson, J. Kalliongis, F.  Laudenbach,
D. McCullough, A. Miller, C. Rourke. Much less attention has been
given to the problem of classifying or describing individual
automorphisms up to isotopy. The paper \cite{JWW:Homeomorphisms} shows
that if $M$ is a Haken manifold or a manifold admitting one of
Thurston's eight geometries, then any homeomorphism of $M$ is isotopic
to one realizing the Nielsen number.

In this paper we use the word ``automorphism" interchangeably to mean
either a self-homeomorphism or its isotopy class.

The paper \cite{UO:Autos} gives a classification up to isotopy of
automorphisms (self-homeo\-morph\-isms) of 3-di\-men\-sion\-al
handlebodies and compression bodies, analogous to the Nielsen-Thurston
classification of automorphisms of surfaces.  Indecomposable
automorphisms analogous to pseudo-Anosov automorphisms were identified
and called {\it generic}. Generic automorphisms were studied further
in \cite{LNC:Generic} and \cite{LNC:Tightness}.  Our understanding of
generic automorphisms remains quite limited.

The goal of this paper is to represent an isotopy class of an
automorphism of any compact, orientable 3-manifold $M$ by a
homeomorphism which is either periodic, has a suitable invariant
surface on which we can decompose the manifold and the automorphism,
or is dynamically nice in the sense that it has invariant laminations
which ``fill" the manifold.  We do not quite achieve this goal. In
some cases we must represent an automorphism as a composition of two
automorphisms each of which fit into the classification.

Throughout this paper $M$ will denote a compact, orientable
3-manifold, possibly with boundary.

If an automorphism $f\from M \to M$ is not periodic (up to
isotopy), then we search for a {\it reducing surface}, see Section
\ref{Reducing}.  This is a suitable surface $F\embed M$ invariant
under $f$ up to isotopy.  In many cases, we will have to decompose
$M$ and $f$ by cutting $M$ on the surface $F$. Cutting $M$ on $F$
we obtain a manifold $M|F$ ($M$-cut-on-$F$) and if we isotope $f$
such that it preserves $F$, we obtain an induced automorphism
$f|F$ ($f$-cut-on-$F$) of $M|F$.

If there are no reducing surfaces for $f$ (and $f$ is not periodic),
then the goal is to find a ``dynamically nice" representative in the
isotopy class of $f$.  This representative will have associated
invariant laminations.  Also, just as in the case of pseudo-Anosov
automorphisms, we will associate to the dynamically nice
representative a growth rate.  A ``best" representative should, of
course, have minimal growth rate.  In the case of automorphisms of
handlebodies it is still not known whether, in some sense, there is a
canonical representative.

The reducing surfaces used in the analysis of automorphisms of
handlebodies and compression bodies, see \cite{UO:Autos}, have a nice
property we call {\it rigidity}:  A reducing surface $F$ for $f\from M
\to M$  is {\it rigid} if $f$ uniquely determines the isotopy class of
the induced automorphism $f|F$ on $M|F$.  As we shall see, the
analysis of automorphisms of many 3-manifolds requires the use of
non-rigid reducing surfaces $F$.  This means that there are many
choices for the induced isotopy class of automorphisms $f|F$, and it
becomes more difficult to identify a ``canonical" or ``best"
representative of $f|F$.  For example, if $F$ cuts $M$ into two
components, then on each component of $M|F$, the induced automorphism
$f|F$ is not uniquely determined.  If there are no reducing surfaces
for $f|F$ on each of these components, and they are not periodic, we
again seek a dynamically nice representative, and in particular, the
best representative should at least achieve the minimal growth rate of
all possible ``nice" representatives.  But now we search for the best
representative not only from among representatives of the isotopy
class of $f|F$, but from among all representatives of all possible
isotopy classes of $f|F$.  This means one may have to isotope $f$
using an isotopy not preserving $F$ to find a better representative of
$f|F$.

There is a subtlety in the definition of reducing surface.  One
must actually deal with automorphisms of {\it pairs} or {\it Haken
pairs}, of the form $(M,V)$ where $M$ is irreducible and $V$ is an
incompressible surface in $\bdry M$.  An automorphism $f:(M,V)\to
(M,V)$ of a Haken pair is an automorphism which restricts on $V$
to an automorphism of $V$.  Even in the case of automorphisms of
handlebodies, Haken pairs play a role: The kind of Haken pair
which arises is called a compression body.

A {\it compression body} is a Haken pair $(Q,V)$ constructed from
a product $V\times [0,1]$ and a collection of balls by attaching
1-handles to $V\times 1$ on the product and to the boundaries of
the balls to obtain a 3-manifold $Q$.  Letting $V=V\times 0$, we
obtain the pair $(Q,V)$. Here $V$ can be a surface with boundary.
The surface $V$ is called the {\it interior boundary} $\bdry_iQ$
of the compression body, while $W=\bdry Q-\intr(V)$ is called the
{\it exterior boundary}, $W=\bdry_eQ$. We regard a handlebody $H$
as a compression body whose interior boundary $V$ is empty. An
automorphism of a compression body is an automorphism of the pair
$(Q,V)$.

A reducing surface $F$ for $f:(M,V)\to (M,V)$ is always $f$-invariant
up to isotopy, but may lie entirely in $\bdry M$. If $F\subseteq\bdry
M$, $F$ is incompressible in $M$ and $f$-invariant up to isotopy, and
if $F\supset V$ up to isotopy, then $F$ is a reducing surface.  In
other words, if it is possible to rechoose the $f$-invariant Haken
pair structure for $(M,V)$ such that $V$ is replaced by a larger $F$
to obtain $(M,F)$, then $f$ is reducible, and $F$ is regarded as a
reducing surface, and is called a {\it peripheral reducing surface}.

We shall give a more complete account of reducing surfaces in Section
\ref{Reducing}.

An automorphism is {\it (rigidly) reducible} if it has a (rigid)
reducing surface.

Returning to the special case of automorphisms of handlebodies, when
the reducing surface is not peripheral, then the handlebody can
actually be decomposed to yield an automorphism $f|F$ of $H|F$.  When
the reducing surface is peripheral, and incompressible in the
handlebody, the reducing surface makes it possible to regard the
automorphism as an automorphism of a compression body whose underlying
topological space is the handlebody, and the automorphism is then
deemed ``reduced."

We continue this introduction with an overview of our current
understanding of automorphisms of handlebodies.  Without
immediately giving precise definitions (see Sections
\ref{Reducing} and \ref{Laminations} for details), we state one of
the main theorems of \cite{UO:Autos}, which applies to
automorphisms of handlebodies.

\begin{thm}\label{HandlebodyClassificationThm}
Suppose $f\from H\to H$ is an automorphism of a connected
handlebody. Then the automorphism is:

1) rigidly reducible,

2) periodic, or

3) generic on the handlebody.

\end{thm}

A {\it generic automorphism} $f$ of a handlebody is defined as an
automorphism which is not periodic and does not have a reducing
surface (see Section \ref{Reducing} for the precise definition).
One can show that this amounts to requiring that $f|_{\bdry H}$ is
pseudo-Anosov and that there are no closed reducing surfaces.

There is a theorem similar to Theorem
\ref{HandlebodyClassificationThm}, classifying automorphisms of
compression bodies.

In the paper \cite{UO:Autos}, the first steps were taken towards
understanding generic automorphisms using invariant laminations of
dimensions 1 and 2.  The theory has been extended and refined in
\cite{LNC:Generic} and \cite{LNC:Tightness}.  We state here only the
theorem which applies to automorphisms of handlebodies; there is an
analogous theorem for arbitrary compression bodies, though in that
case no 1-dimensional invariant lamination was described.  We shall
remedy that omission in this paper.

If $H$ is a handlebody, in the following theorem $H_0$ is a {\it
concentric} in $H$, meaning that $H_0$ is embedded in the interior of
$H$ such that $H-\intr(H_0)$ has the structure of a product $\bdry
H\times I$.

\begin{thm}\label{HandlebodyLaminationThm}
Suppose $f\from H\to H$ is a generic automorphism of a $3$-dimensional
handlebody. Then there is a 2-dimensional measured lamination
$\Lambda\embed \intr (H)$ with transverse measure $\mu$ such that, up
to isotopy, $f((\Lambda,\mu))=(\Lambda,\lambda \mu)$ for some
$\lambda> 1$. The lamination has the following properties:

1) Each leaf $\ell$ of $\Lambda$ is an open $2$-dimensional disk.

2) The lamination $\Lambda$ fills $H_0$, in the sense that the
components of $H_0-\Lambda$ are contractible.

3) For each leaf $\ell$ of $\Lambda$,  $\ell-\intr (H_0)$ is
incompressible in $H-\intr(H_0)$.

4) $\Lambda\cup \bdry H$ is a closed subset of $H$.

There is also a 1-dimensional lamination $\Omega$, transverse to
$\Lambda$, with transverse measure $\nu$ and a map $\omega\from
\Omega\to \intr (H_0)$ such that
$f(\omega(\Omega,\nu))=\omega(\Omega,\nu/\lambda)$.  The map $f$ is an
embedding on $f\inverse(N(\Lambda))$ for some neighborhood
$N(\Lambda)$. The statement that
$f(\omega(\Omega,\nu))=\omega(\Omega,\nu/\lambda)$ should be
interpreted to mean that there is an isomorphism $h\from
(\Omega,\nu)\to (\Omega,\nu/\lambda)$ such that
$f\circ\omega=\omega\circ h$.
\end{thm}

We note that the laminations in the above statement are ``essential"
only in a rather weak sense.  For example, lifts of leaves of
$\Lambda$ to the universal cover of $H$ are not necessarily properly
embedded.  The map $\omega$ need not be proper either:  If we regarded
$w\from \Omega\to H$ as a homotopy class, in some examples there would
be much unravelling.  The lamination $\Omega$ can also be regarded as
a ``singular" embedded lamination, with finitely many singularities.

\begin{remark}
The laminations $\Lambda$ and $\Omega$ should be regarded as a
pair of dual laminations.  Both are constructed from a chosen
``complete system of essential discs," $\cal E$, see Section
\ref{Laminations}. One attempts to find the most desirable pairs
$\Lambda$ and $\Omega$. For example, it is reasonable to require
laminations with minimal $\lambda$, but even with this requirement
it is not known whether the invariant laminations are unique in
any sense.
\end{remark}

One of the authors has shown in his dissertation,
\cite{LNC:Generic} and in the preprint \cite{LNC:Tightness}, that
it is often possible to replace  the invariant laminations in the
above theorem by laminations satisfying a more refined condition.
The condition on the 2-dimensional lamination is called {\it
tightness} and is defined in terms of both invariant laminations.
He shows that invariant lamination having this property
necessarily achieve the minimum eigenvalue $\lambda$ (growth rate)
among all (suitably constructed) measured invariant laminations.
One remaining problem is to show that every generic automorphism
admits a tight invariant lamination; this is known in the case of
an automorphism of a genus 2 handlebody and under certain other
conditions.  The growth rate $\lambda$ for a tight invariant
lamination associated to a generic automorphism $f$ is no larger
than the growth rate for the induced pseudo-Anosov automorphism on
$\bdry H$, see \cite{LNC:Generic}.  In particular, this means that
for a genus 2 handlebody a growth rate can be achieved which is no
larger than the growth rate of the pseudo-Anosov on the boundary.
We hope that the tightness condition on invariant laminations will
be sufficient to yield some kind of uniqueness.  The most
optimistic goal (perhaps too optimistic) is to show that a tight
invariant 2-dimensional lamination for a given generic
automorphism $f\from H\to H$ is unique up to isotopy.

\hop

In our broader scheme for classifying automorphisms of arbitrary
compact 3-manifolds, an important test case is to classify up to
isotopy automorphisms $f\from M\to M$, where $M$ is a connected sum of
$S^2\times S^1$'s.  In calculations of the mapping class group, this
was also an important case, see \cite{FL:Spheres} and
\cite{FL:HomotopyIsotopy}. The new idea for dealing with this case is
that there is {\em always} a surface invariant up to isotopy, a
reducing surface which, as we shall see, is non-rigid. This invariant
surface is the Heegaard splitting surface coming from the Heegaard
splitting in which $M$ is expressed as the double of a handlebody.  We
call this a {\it symmetric} Heegaard splitting.  A remark by
F. Waldhausen at the end of his paper \cite{FW:Heegaard} outlines a
proof that every minimal genus Heegaard splitting of a connected sum
of $S^2\times S^1$'s is symmetric. We prove that there is only one
symmetric splitting up to isotopy, and we call this the {\it canonical
splitting} of $M$.  We were surprised not to find this result in the
literature, and it is still possible that we have overlooked a
reference.

\begin{thm} \label{CanonicalHeegaardTheorem}
Any two symmetric Heegaard splittings of a connected sum of $S^2\times
S^1$'s are isotopic.
\end{thm}

Thus when an automorphism $f$ of $M$ is isotoped so it preserves $F$,
it induces automorphisms of $H_1$ and $H_2$, which are mirror images
of each other.   However, there is not a unique way of isotoping an
automorphism $f$ of $M$ to preserve $F$, see
\cite{FL:ConnectedS2TimesS1} and Example \ref{NonUniqueExample}.  Thus
the canonical splitting of a connected sum of $S^2\times S^1$'s gives
a class of examples of non-rigid reducing surfaces.

In the important case where $M$ is a connected sum of $k$ copies
of $S^2\times S^1$, if $f:M\to M$ is an automorphism, we next
consider the best choice of representative of $f|F$, which is the
induced automorphism on $M|F$, a pair of genus $k$ handlebodies.
From the symmetry of the canonical splitting, the induced
automorphisms of the two handlebodies are conjugate, and we denote
each of these by $g$.

 To state our theorem about automorphisms of connected sums of
$S^2\times S^1$'s, we first mention that given a generic automorphism
$g$ of a handlebody $H$ and an invariant 2-dimensional lamination as
in Theorem \ref{HandlebodyLaminationThm}, the lamination gives a
quotient map $H\to G$ where $G$ is a graph, together with a homotopy
equivalence of $G$ induced by $g$.  Typically, the 2-dimensional
lamination {\em cannot} be chosen such that the induced homotopy
equivalence on the quotient graph is a train track map in the sense of
\cite{BH:Tracks}.  In the exceptional cases when there exists a
2-dimensional invariant lamination so the quotient homotopy
equivalence is a train track map, we say $g$ is a {\it train track
generic automorphism}.

\begin{thm}[Train Track Theorem]\label{TrainTrackTheorem}
If $M$ is a connected sum of $S^2\times S^1$'s, $M=H_1\cup H_2$ is the
canonical splitting, and $f\from M\to M$ is an automorphism, then $f$
is non-rigidly reducible on the canonical splitting surface, and can
be isotoped so that $f$ preserves the splitting and the induced
automorphism $g$ on $H_1$ is

1) periodic,

2) rigidly reducible, or

3) train track generic.
\end{thm}

In the above, not so much information is lost in passing to the
induced automorphism on $\pi_1(H_1)$, which is not surprising in view
of Laudenbach's calculation of the mapping class group of $M$.
F. Laudenbach showed in \cite{FL:Spheres} that the map from the
mapping class group of $M$ to $\text{Out}(\pi_1(M))$ has kernel
$\integers_2^n$ generated by rotations in a maximal collection of $n$
disjoint non-separating 2-spheres in $M$.

Note that when the automorphism $g$ is generic, we have not shown that
there is a unique pair of invariant laminations such that the quotient
is a train track map.

For a typical generic automorphism, any invariant 1-dimensional
lamination has much self-linking and is somewhat pathological.  Even
for a train-track generic map, some of the self-linking and pathology
may remain.  We say a 1-dimensional lamination $\Omega\embed M$ is
{\it tame} if any point of $\Omega$ has a flat chart neighborhood in
$M$.

\begin{conj}[Tameness Conjecture]\label{TamenessConjecture}
If $M$ is a connected sum of $S^2\times S^1$'s, $M=H_1\cup H_2$ is the
canonical splitting, and $f\from M\to M$ is a train-track generic
automorphism, then $f$ can be isotoped so that $f$ preserves the
splitting and the induced automorphism $g$ on $H_1$ is train-track
generic and has a tame 1-dimensional invariant lamination.
\end{conj}

To proceed, we must say something about automorphisms of
compression bodies.  Generic automorphisms of compression bodies
are defined like generic automorphisms of handlebodies; in
particular, if $f\from Q\to Q$ is an automorphism, the induced
automorphism $\bdry_ef$ on $\bdry_eQ$ is pseudo-Anosov.  Also,
there should be no (suitably defined) reducing surfaces (see
Section \ref{Reducing}).

We will adopt some further terminology for convenience.  A connected
sum of $S^2\times S^1$'s is a {\it sphere body}.  Removing open
ball(s) from a manifold is called {\it holing} and the result is a
{\it holed manifold}, or a manifold with sphere boundary components.
Thus holing a sphere body yields a {\it holed sphere body}.  A {\it
mixed body} is a connected sum of $S^2\times S^1$'s and handlebodies.
If we accept a ball as a special case of a handlebody, then a holed
mixed body is actually also a mixed body, since it can be obtained by
forming a connect sums with balls.  Nevertheless, we distinguish a
{\it holed mixed body}, which has sphere boundary components, from a
mixed body, which does not.

It is often convenient to cap sphere boundary components of a
holed manifold with balls, distinguishing the caps from the
remainder of the manifold.  Thus given a manifold $M$ with sphere
boundary components, we cap with balls, whose union we denote
$\mathcal{D}$, to obtain a {\it spotted manifold} (with
3-dimensional balls).  For classifying automorphisms, spotted
manifolds are just as good as holed manifolds.

Suppose now that $M$ is a mixed body.  Again, there is a canonical
Heegaard splitting $M=H\cup Q$, where $H$ is a handlebody and $Q$ is a
compression body, see Theorem \ref{CanonicalHeegaardTheorem2}.
Whereas in the decomposition of an irreducible manifold with
compressible boundary the Bonahon characteristic compression body has
exterior boundary in the boundary of the manifold, our compression
body $Q$ here has interior boundary $\bdry_iQ=\bdry M$.  Again,  this
splitting lacks rigidity for the same reasons.

In Section \ref{MixedBodies}, we undertake the study of automorphisms
of mixed bodies. We prove a theorem analogous to Theorem
\ref{TrainTrackTheorem} (the Train Track Theorem), see Theorem
\ref{CompressionTrainTrackTheorem}.

In a classification of automorphisms of compact 3-manifolds, the
most difficult case is the case of an arbitrary reducible
3-manifold.  We will assume the manifold has no sphere boundary
components; if the manifold has boundary spheres, we cap them with
balls. Let $B$ be the 3-ball. Suppose that a (possibly holed)
mixed body $R\neq B$ is embedded in $M$. We say that $R$ is {\em
essential} if any essential sphere in $R$ or sphere component of
$\bdry R$ is also essential in $M$. An automorphism $f$ of $M$
which preserves an essential mixed body $R$ up to isotopy and
which is periodic on $M-\intr(R)$ is called an {\it adjusting
automorphism}.  The surface $\bdry R$ is $f$-invariant up to
isotopy, but we  do not regard it as a reducing surface, see
Section \ref{Reducing}.

Note that an essential mixed body $R$ in an irreducible manifold
is a handlebody.  Since an adjusting automorphism preserving $R$
is isotopic to a periodic map on $\bdry R$, it is then also
isotopic to a periodic map on $R$, so such an adjusting
automorphism is not especially useful.

For the classification of automorphisms of reducible compact
3-manifolds, we will use the following result, which is based on a
result due to E. C{\'e}sar de S{\'a} \cite{EC:Automorphisms}, see also
M. Scharlemann in Appendix A of \cite{FB:CompressionBody}, and
D. McCullough in \cite{DM:MappingSurvey} (p. 69).

\begin{thm}\label{AdjustingTheorem}
Suppose $M$ has irreducible summands $M_i$, $i=1,\ldots k$, and
suppose $\hat M_i$ is $M_i$ with one hole, and one boundary sphere
$S_i$.  Suppose the $\hat M_i$ are disjoint submanifolds of $M$,
$i=1,\ldots k$.   Then the closure of $M-\cup_i \hat M_i$ is an
essential holed sphere body $\hat M_0$.  ($\hat M_0$ is a holed
connected sum of $S^2\times S^1$'s or possibly a holed sphere.) If
$f\from M\to M$ is an automorphism, then there is an adjusting
automorphism $\hslash$ preserving an essential (possibly holed)
mixed body $R\supset \hat M_0$ in $M$ with the property that
$g=\hslash\circ f$ is rigidly reducible on the spheres of $\bdry
\hat M_0$.

The automorphism $g$ is determined by $f$ only up to application of an
adjusting automorphism preserving $\bdry \hat M_0$.

Note that we can write $f$ as a composition $f=\hslash\inverse\circ
g=h\circ g$, where $h=\hslash\inverse$ is an adjusting automorphism
and $g$ preserves each $\hat M_i$.

The automorphism $\hat g_i=g|_{\hat M_i}$ clearly uniquely determines
an automorphism $g_i$ of $M_i$ by capping (for $i=0$, $M_0$ is $\hat
M_0$ capped). On the other hand $g_i$ determines $\hat g_i$ only up to
application of an adjusting automorphism in  $\hat M_i$.
\end{thm}

The main ingredients used in this result are {\em slide automorphisms
of $M$}. There are two types, one of which is as defined in
\cite{DM:MappingSurvey}. We shall describe them briefly in Section
\ref{Adjusting}.

With a little work, see Section \ref{Adjusting}, we can reformulate
Theorem \ref{AdjustingTheorem} in terms of 4-dimensional compression
bodies as follows.  Let $Q$ be a 4-dimensional compression body
constructed from a disjoint union of products $M_i\times I$ and a
4-dimensional ball $K_0$ by attaching one 1-handle joining $\bdry
K_0=M_0$ to $M_i\times 1$ and one 1-handle from $\bdry K_0$ to itself
for every $S^2\times S^1$ summand of $M$.  Then $\bdry_eQ$ is
homeomorphic with and can be identified with $M$, while $\bdry_iQ$ is
homeomorphic to and can be identified with the disjoint union of
$M_i$'s, $i\ge 1$.

\begin{corollary} \label{AdjustingCorollary}  With the
hypotheses of Theorem \ref{AdjustingTheorem} and with the
4-dimensional compression body $Q$ constructed as above, there is
an automorphism $\bar f:Q\to Q$ of the compression body such that
$\bar f|_{\bdry_eQ}=f$. Also, there is an automorphism $\bar
\hslash:Q\to Q$ of the compression body such that $\bar
\hslash|_{\bdry_eQ}=\hslash$, $\bar
\hslash|_{\bdry_iQ}=\text{id}$. Also $\bar\hslash\circ \bar f$
gives the adjusted automorphisms on $\bdry_iQ$, which is the union
of $M_i$'s, $i\ge 1$.

In other words, the automorphism $f:M\to M$ is cobordant over $Q$ to
an automorphism of the disjoint union of the $M_i$'s via an automorphism
$\bar \hslash$ of $Q$ rel $\bdry_iQ$.
\end{corollary}

We will refer to the automorphism $\bar \hslash$ as an {\it adjusting automorphism of the compression body $Q$}.

Before turning to the classification of automorphisms, we mention
another necessary and important ingredient.

\begin{thm} \label{BonahonTheorem} (F. Bonahon's characteristic compression body)  Suppose $M$
is a compact manifold with boundary.  Then there is a compression body
$Q\embed M$ which is unique up to isotopy with the property that
$\bdry_eQ=\bdry M$ and $\bdry_iQ$ is incompressible in
$\overline{M-Q}$.
\end{thm}

Bonahon's theorem gives a canonical decomposition of any compact
irreducible 3-manifold with boundary, dividing it into a compression
body and a compact irreducible manifold with incompressible boundary.
The interior boundary of the characteristic compression body is a
rigid reducing surface.

Finally, in the classification, we often encounter automorphisms of
{\it holed manifolds}. That is to say, we obtain an automorphism
$f:\hat M\to \hat M$ where $\hat M$ is obtained from $M$ by removing a
finite number of open balls.  We have already observed that an
automorphism of a holed manifold $\hat M$ is essentially the same as
an automorphism of manifold pair $(M,\mathcal{D})$ where $\mathcal{D}$
is the union of closures of the removed open balls,
i.e. $\mathcal{D}=\overline{M-\hat M}$.  An automorphism of a spotted
manifold $(M,\mathcal{D})$ or holed manifold $\hat M$ determines an
automorphism of the manifold $M$ obtained by forgetting the spots or
capping the boundary spheres of $\hat M$.  Clearly, there is a loss of
information in passing from the automorphism $\hat f$ of the spotted
or holed manifold to the corresponding automorphism $f$ of the
unspotted or unholed manifold. The relationship is that $\hat f$
defined on the spotted manifold is isotopic to $f$ if we simply regard
$\hat f$ as an automorphism of $M$.  To obtain $\hat f$ from $f$, we
compose with an isotopy $h$ designed to ensure that $h\circ f$
preserves the spots.  There are many choices for $h$, and two such
choices $h_1$ and $h_2$ differ by an automorphism of $(M,B)$ which can
be realized as an isotopy of $M$.  This means $h_1\circ h_2\inverse$
is an automorphism of $(M,B)$ realized by an isotopy of $M$. In
describing the classification, we will not comment further on the
difference between automorphisms of manifolds and automorphisms of
their holed counterparts.  Of course any manifold with sphere boundary
components can be regarded as a holed manifold.  Much more can be said
about automorphisms of spotted manifolds, but we will leave this topic
to another paper.

\subsection*{Automorphisms of Irreducible and $\partial$-irreducible Manifolds}

The following outlines the classification of automorphisms of
irreducible and $\partial$-irre\-ducible manifolds. If the
characteristic manifold (Jaco-Shalen, Johannson) is non-empty or
the manifold is Haken, we use \cite{WJPS:Characteristic},
\cite{KJ:Characteristic}.  For example, if the characteristic
manifold is non-empty and not all of $M$, then the boundary of the
characteristic manifold is a rigid incompressible reducing
surface. There are finitely many isotopy classes of automorphisms
of the complementary ``simple'' pieces \cite{KJ:Characteristic}.
For automorphisms of Seifert-fibered pieces we refer to
\cite{PS:Geometries,JWW:Homeomorphisms}. The idea is that in most
cases an automorphism of such a manifold preserves a Seifert
fibering,  therefore the standard projection yields an
automorphism of the base orbifold. The article
\cite{JWW:Homeomorphisms} gives a classification of automorphisms
of orbifolds similar to Nielsen-Thurston's classification of
automorphisms of surfaces. It is then used for a purpose quite
distinct from ours: their goal is to realize the Nilsen number in
the isotopy class of a $3$-manifold automorphism. Still, their
classification of orbifold automorphisms is helpful in our
setting. The idea is that invariant objects identified by the
classification in the base orbifold lift to invariant objects in
the original Seifert fibered manifold.

In the following we let $T^n$ be the $n$-torus, $D^2$ the
$2$-disc, $I$ the interval ($1$-disc), $\mathbb{P}(2,2)$ the
orbifold with the projective plane as underlying space and two
cone points of order $2$, and $M_{\mathbb{P}(2,2)}$ the Seifert
fibered space over $\mathbb{P}(2,2)$.

\begin{thm}[\cite{JWW:Homeomorphisms}]\label{T:PreservedFibering}
Suppose that $M$  is a compact  orientable Seifert fibered space
which is not $T^3$, $M_{\mathbb{P}(2,2)}$, $D^2\times S^1$ or
$T^2\times I$. Given $f\colon M\to M$ there exists a Seifert
fibration which is preserved by $f$ up to isotopy.
\end{thm}

We discuss the classification of automorphisms of these four
exceptional cases later. In the general case we consider the
projection $M\to X$ over the fiber space $X$, which is an
orbifold, and let the surface $Y$ be the underlying space of $X$.
If $f\colon M\to M$ preserves the corresponding fiber structure we
let $\hat f\colon X\to X$ be the projected ``automorphism'' of
$X$. An {\em automorphism of $X$} is an automorphism of $Y$
preserving the orbifold structure. Let $\sing(X)$ be the set of
singular points of $X$. Note that we are assuming that $M$ is
orientable therefore there are no reflector lines in $\sing(X)$,
which then consists of isolated cone points. Moreover $\hat f$ has
to map a cone point to another of the same order. An isotopy of
$\hat f$ is an isotopy restricted to $X-\sing(X)$.

The classification of automorphisms of $X$ is divided in three
cases, depending on the sign of $\chi(X)$ (as an orbifold Euler
characteristic). Below we essentially follow
\cite{JWW:Homeomorphisms}.

\begin{enumerate}
    \item $\chi(X)>0$.
    In this case, unless $\sing(X)=\emptyset$,
    the underlying space $Y$ is either $S^2$ or $\mathbb{P}^2$ and $X$ has a few singular points.
    Any map is isotopic to a periodic one.

    \item $\chi(X)=0$.
    For most orbifolds $X$ any automorphism is isotopic to a periodic one.
    The only two exceptions are $T^2$ (without singularities) and $S^2(2,2,2,2)$, the sphere with
    four cone points of order 2.

    If $X=T^2$ then an automorphism $\hat f\colon X\to X$
    is isotopic to either i) a periodic one, ii) reducible, preserving a curve or
    iii) an Anosov map. In ii) the lift of an invariant curve in $X$
    yields a torus reducing surface for $f$ in $M$ (if the lift of the curve is a 1-sided Klein
    bottle we consider the torus boundary of its regular neighborhood). Such a torus may separate
    a $K\widetilde{\times} I$ piece, which is the orientable $I$-bundle over the Klein bottle $K$.
    An automorphism of such a piece is isotopic to a periodic map.
    In iii) the invariant foliations lift to foliations on $M$ invariant
    under $f$. The leaves are either open annuli or open M\"obius bands.

    If $X=S^2(2,2,2,2)$ then $X'=X-\sing(X)$ is a four times
    punctured sphere. Since $\chi(X')<0$ the automorphism $\hat f|_{X'}$ is subject to
    Nielsen-Thurston's classification, which isotopes it to
    be either i) periodic, ii) reducible or iii) pseudo-Anosov. As
    above, in ii) a reduction of $\hat f$ yields a reduction of $f$ along
    tori, maybe separating $K\widetilde{\times} I$ pieces. In iii)
    we consider the foliations on $X$ invariant under $\hat f$.
    They may have 1-prong singularities in $\sing(X)$. We can assume that there are
    no other singularities since $\chi(X)=0$, with $X$ covered by
    $T^2$. But all cone points have order 2, therefore the lifts of the foliations
    to $M$ are non-singular foliations invariant under $f$.

    \item $\chi(X)<0$. Also in this case $X'=X-\sing(X)$ is a surface
    with $\chi(X')<0$, therefore $\hat f|_{X'}$ is subject to
    Nielsen-Thurston's classification. A reducing curve system for
    $\hat f|_{X'}$ is a reducing curve system for $\hat
    f$. As before, this yields a reducing surface for $f$ consisting of
    tori. Some may be thrown out, in case there are parallel
    copies.

    In the pseudo-Anosov case we consider the corresponding
    invariant laminations on $X'$, which lift to invariant laminations on
    $M$ by vertical open annuli or M\"obius bands.
\end{enumerate}

In the above we did not say anything about $f\colon M\to M$ when
$\hat f\colon X\to X$ is periodic. There is not much hope in
trying to isotope $f$ to be periodic, since there are
homeomorphisms of $M$ preserving each fiber but not periodic (e.g,
``vertical'' Dehn twists along vertical tori). What can be done is
to write $f$ as a composition $f=g\circ h$ where $h$ is periodic
and $g$ preserves each fiber.

The outline above deals with Seifert fibered spaces $M$ for which
there is a fibering structure preserved by $f\colon M\to M$.
Recall that this excludes $T^3$, $M_{\mathbb{P}(2,2)}$, $D^2\times
S^1$ or $T^2\times I$ \ref{T:PreservedFibering}, which we consider
now. If $M=D^2\times S^1$ then $f$ is isotopic to a power of a
Dehn twist along the meridional disc. If $M=T^2\times I$ then $f$
is isotopic to a product $g\times(\pm\id_I)$, where $\id_I$ is the
identity and $-\id_I$ is the orientation reversing isometry of
$I$. Assume that $g\colon T^2\to T^2$ is either periodic,
reducible (along a curve) or Anosov.

In the case $M=T^3$ any automorphism is isotopic to a {\em linear
automorphism} (i.e., an automorphism with a lift to
$\widetilde{M}=\mathbb{R}^3$ which is linear). The study of
eigenspaces of such an automorphism yield invariant objects.
Roughly, such an automorphism may be periodic, reducible
(preserving an embedded torus) or may have invariant 2-dimensional
foliations, whose leaves may consist of dense open annuli or
planes. In some cases it may also be relevant to further consider
invariant 1-dimensional foliations.

In the case $M=M_{\mathbb{P}(2,2)}$ there is a cell decomposition
which is preserved by any automorphism up to isotopy. It thus
imply that any automorphism is isotopic to a periodic one.

This last case completes the outline of the classification for
Seifert fibered spaces.

If the irreducible and $\bdry$-irreducible manifold $M$ is closed
and hyperbolic, it is known that any automorphism is periodic, see
\cite{WPT:Notes}.  If the manifold is closed and spherical, recent
results of D. McCullough imply the automorphisms are periodic, see
\cite{DM:Elliptic}.

\subsection*{The Classification} The following outline
classifies automorphisms of 3-manifolds in terms of automorphisms of
compression bodies, automorphisms of mixed bodies, and automorphisms
of irreducible, $\bdry$-irreducible 3-manifolds.  (Here we regard
handlebodies as compression bodies.) In every case, we consider an
automorphism $f\from M\to M$.

i) If $M$ is reducible use Theorem \ref{AdjustingTheorem} and write
$f=h\circ g$, where $g$ is rigidly reducible on essential spheres, and
$h$ is non-rigidly reducible. The reductions yield automorphisms of
holed irreducible 3-manifolds (the irreducible summands of $M$ with
one hole each), and of (holed) mixed bodies, which we consider in the
following cases.  The automorphism $h$ also yields, after cutting on
$\bdry R$ of Theorem \ref{AdjustingTheorem}), periodic automorphisms
of (holed) irreducible manifolds with incompressible boundary.

Alternatively, using Corollary \ref{AdjustingCorollary}, the
automorphism $f:M\to M$ is cobordant over a 4-dimensional compression
body to an ``adjusted" automorphism $f_a$ of the disjoint union of the
irreducible summands $M_i$.

ii) If $M$ is irreducible, with $\bdry M\ne \emptyset$ let $Q$ be the
characteristic compression body of Theorem \ref{BonahonTheorem}, then
$f$ induces unique automorphisms on $Q$ and $\overline{M-Q}$, , see
\cite{FB:CompressionBody} and \cite{UO:Autos}. $Q$ is a compression
body and $\overline{M-Q}$ is irreducible and $\bdry$-irreducible,
whose automorphisms we consider in following cases.

iii)  If $M$ is a holed sphere body (a connected sum of $S^2\times
S^1$'s with open balls removed), as is the case when the automorphism
is obtained by decomposition from an adjusting automorphism for a
different automorphism of a different 3-manifold, then we first cap
the boundary spheres to obtain an induced automorphism of a sphere
body.  Then we use the Train Track Theorem (Theorem
\ref{TrainTrackTheorem} and Theorem
\ref{CompressionTrainTrackTheorem}), to obtain a particularly nice
representative of $f$. If the Tameness Conjecture,  Conjecture
\ref{TamenessConjecture}, is true, then one can find an even nicer
representative.

iv)  If $M$ is a (holed) mixed body, as is the case when the
automorphism is an adjusting automorphism for a different
automorphism of a different 3-manifold, then we use the canonical
splitting surface (Theorem \ref{CanonicalHeegaardTheorem2}) as a
non-rigid reducing surface for $f$ which cuts $M$ into a
handlebody $H$ and a compression body $Q$.  In the handlebody,
there is another non-rigid reducing surface (see Section
\ref{MixedBodies}), which makes it possible to cut the handlebody
$H$ into a compression body $Q'$ isomorphic to $Q$, and a possibly
disconnected handlebody $H'$, with a component of genus $g$
corresponding to each genus $g$ component of $\bdry M$. Theorem
\ref{CompressionTrainTrackTheorem} gives a best representative,
yielding the simplest  invariant laminations for $f$ restricted to
$Q$, $Q'$, and $H'$.  The restrictions of $f$ to $Q$ and $Q'$ are
then automorphisms of compression bodies periodic on the interior
boundary.

v)  If $M$ is a compression body, the automorphisms are classified by
Theorem \ref{HandlebodyClassificationThm} and its analogue, see
\cite{UO:Autos}, which says the automorphism is periodic, rigidly
reducible or generic.  See Section \ref{Laminations} for a description
of generic automorphisms.

vi)  If $M$ is irreducible and $\bdry$-irreducible, then the
automorphism $f$ is classified following the previous subsection.

\centerline {---}

\hop

If it were true that every automorphism of a reducible manifold had a
closed reducing surface, possibly a collection of essential spheres,
then decomposition on the reducing surface would make it possible to
analyze the automorphism of the resulting manifold-with-boundary in a
straightforward way, without the necessity of first writing the
automorphism as a composition of an adjusting automorphism and an
automorphism preserving a system of reducing spheres.  This would
yield a dramatic improvement of our classification.  For this reason,
we pose the following:

\begin{question}\label{ReducingQuestion} Is it true that every automorphism of a reducible manifold admits a closed reducing
surface, which cannot be isotoped to be disjoint from a system of essential spheres cutting $M$ into holed irreducible
manifolds? In particular, does every automorphism of a closed reducible 3-manifold admit a closed reducing surface of this
kind?
\end{question}

To understand the significance of this question, one needs a more
detailed understanding of reducing surfaces.  Assuming that an
automorphism $f:M\to M$ has no reducing surfaces consisting of an
$f$-invariant collection of essential spheres, a closed reducing
surface $F$ is required to have the property that $M|F$ is
irreducible, see Section \ref{Reducing}.

We briefly describe the much improved classification of automorphisms
which would be possible if the answer to Question
\ref{ReducingQuestion} were positive.  If $f:M\to M$ is an
automorphism and $M$ is reducible, we first seek reducing surfaces
consisting of essential spheres.  If these exist, we decompose and cap
the spheres.  If there are no sphere reducing surfaces, we seek closed
reducing surfaces.  Decomposing on such a reducing surface gives an
automorphism of an irreducible manifold with non-empty boundary, where
the irreducibilty of the decomposed manifold comes from the
requirement that cutting on a reducing surface should yield an
irreducible manifold.  For the resulting irreducible manifold with
boundary, we decompose on the interior boundary of the characteristic
compression body as before, obtaining an automorphism of a manifold
with incompressible boundary and an automorphism of a compression
body.  (Either one could be empty, and the compression body could have
handlebody component(s).)  Automorphisms of compression bodies are
analyzed using the existing (incomplete) theory, the automorphism of
the irreducible, $\bdry$-irreducible 3-manifold is analyzed as in the
classification given above, using mostly well-established 3-manifold
theory.

\section{Canonical Heegaard splittings}\label{Splitting}

Given a handlebody $H$ of genus $g$ we can double the handlebody to
obtain a connected sum $M$ of $g$ copies of $S^2\times S^1$, a genus
$g$ sphere body. We choose $g$ non-separating discs in $H$ whose
doubles give $g$ non-separating spheres in $M$.  We then say that the
Heegaard splitting $M=H\cup (M- \intr(H))$ is a {\it symmetric
Heegaard splitting}.  The symmetric splitting is characterized by the
property that any curve in $\bdry H$ which bounds a disc in $H$ also
bounds a disc in $M- \intr(H)$.

The following is our interpretation of a comment at the end of a paper
of F. Waldhausen, see \cite{FW:Heegaard}.

\begin{lemma} (F. Waldhausen)
Any minimal genus Heegaard splitting of a sphere body is symmetric.
\end{lemma}

We restate a theorem which we will prove in this section:

\begin{restate} [Theorem \ref{CanonicalHeegaardTheorem}]
Any two symmetric Heegaard splittings of a sphere body are isotopic.
\end{restate}

One would expect, from the fundamental nature of the result, that this
theorem would have been proved much earlier.  A proof may exist in the
literature, but we have not found it.

In the first part of this section, $M$ will always be a genus $g$
sphere body.  An {\it essential system of spheres} is the isotopy
class of a system of disjointly embedded spheres with the property
that none of the spheres bounds a ball.   The system is {\it
primitive} if no two spheres of the system are isotopic.  A system
$\cal S$ is said to be {\it complete} if cutting $M$ on $\cal S$
yields a collection of holed balls.  There is a cell complex of
weighted primitive systems of essential spheres in $M$ constructed as
follows.  For any primitive essential system $\cal S$ in $M$, the
space contains a simplex corresponding to projective classes of
systems of non-negative weights on the spheres.  Two simplices
corresponding to primitive essential systems ${\cal S}_1$ and ${\cal
S}_2$ are identified on the face corresponding to the maximal
primitive essential system ${\cal S}_3$ of spheres common to ${\cal
S}_1$ and ${\cal S}_2$.  We call the resulting space the {\it sphere
space} of $M$.  In the sphere space, the union of open cells
corresponding to complete primitive systems is called the {\it
complete sphere space} of $M$

\begin{lemma}\label{SphereSpaceLemma}
The complete sphere space of a sphere body $M$ is connected.  In other
words, it is possible to pass from any primitive complete essential
system to any other by moves which replace a primitive complete
essential system by a primitive complete subsystem or a primitive
complete essential supersystem.
\end{lemma}

Note: A. Hatcher proved in \cite{AH:HomologicalStability} that the
sphere space is contractible.   It may be possible to adapt the
following proof to show that the complete sphere space is contractible.

\begin{proof}
Suppose ${\cal S}={\cal S}_0$ and ${\cal S}'$ are two primitive
complete systems of essential spheres.  We will abuse notation by
referring to the union of spheres in a system by the same name as the
system. Without loss of generality, we may assume ${\cal S}$ and
${\cal S}'$ are transverse.

Consider ${\cal S}\cap{\cal S}'$ and choose a curve of intersection
innermost in ${\cal S}'$.  Remove the curve of intersection by isotopy
of ${\cal S}$ if possible.  Otherwise, surger $\cal S$ on a disc
bounded in ${\cal S}'$ by the curve to obtain ${\cal S}_2$.   Now
${\cal S}_2$ is obtained from ${\cal S}={\cal S}_0$ by removing one
sphere, the surgered sphere $S_0$ of ${\cal S}_0$, say, and adding two
spheres, $S_{21}$ and $S_{22}$, which result from the surgery.  We can
do this in two steps, first adding the two new spheres to obtain
${\cal S}_1$, then removing the surgered sphere to obtain ${\cal
S}_2$.  If one of the spheres $S_{21}$ and $S_{22}$ is inessential,
this shows that the curve of intersection between ${\cal S}={\cal
S}_0$ and ${\cal S}'$ could have been removed by isotopy.  We note
that $S_{21}$ or $S_{22}$, or both, might be isotopic to spheres of
${\cal S}_0$. If the new spheres are essential, we claim that ${\cal
S}_1$ and ${\cal S}_2$ are complete systems of essential spheres.  It
is clear that ${\cal S}_1$ is complete, since it is obtained by adding
two essential spheres $S_{21}$ or $S_{22}$ to ${\cal S}_0$.  We then
remove $S_0$ from ${\cal S}_1$ to obtain ${\cal S}_2$.  There is, of
course a danger that removing an essential sphere from an essential
system ${\cal S}_1$ yields an incomplete system.  However, if we
examine the effect of surgery on the complementary holed balls, we see
that first we cut one of them on a properly embedded disc (the surgery
disc), which clearly yields holed balls.  Then we attach a 2-handle to
the boundary of another (or the same) holed ball, which again yields a
holed ball.  Thus, ${\cal S}_2$ is complete.  If ${\cal S}_2$ is not
primitive, i.e. contains two isotopic spheres, we can remove one,
until we are left with a primitive complete essential system.

Inductively, removing curves of intersection of ${\cal S}_{k}$ with
${\cal S}'$ either by isotopy or by surgery on an innermost disc in
${\cal S}'$ as above, we must obtain a sequence ending with ${\cal
S}_{n-2}$ disjoint from ${\cal S}'$. Next, we add the spheres of
${\cal S}'$ not in ${\cal S}_{n-2}$ to ${\cal S}_{n-2}$ to obtain
${\cal S}_{n-1}$, which is still complete.  Finally we remove all
spheres from ${\cal S}_{n-1}$ except those of ${\cal S}'$ to get
${\cal S}_{n}={\cal S}'$.
\end{proof}

A {\it ball with $r$ spots}, or an {\it $r$-spotted ball} is a pair
$(K,E)$ where $K$ is a 3-ball and $E$ is consists of $r$ discs in
$\bdry K$. Let $P$ be a holed ball, with $r\ge 1$ boundary components
say.  Then we say $(K,E)\embed (P,\bdry P)$ is {\it standard} if the
closure of $P-K$ in $P$ is another $r$-spotted ball.  In that case,
the complementary spotted ball is also standard.  In practical terms,
a standard $r$-spotted ball is ``unknotted."   Another point of view
is the following:  If $(K,E)\embed (P,\bdry P)$ is standard, then
$(P,\bdry P)=(K,E)\cup (K',E')$, where $(K',E')$ complementary to
$(K,E)$, is a ``symmetric" splitting of $(P,\bdry P)$, similar to a
symmetric splitting of a sphere body.

\begin{lemma} \label{SpottedBallStandardLemma}
(a) A standard embedded $r$-spotted ball $(K,E)$ in an $r$-holed
3-sphere $(P,\bdry P)$ ($r\ge 1$) is unique up to isotopy of pairs.

(b) Suppose $(K,E)$ and $(K',E')$ are two different standard spotted
balls in $(P,\bdry P)$ and suppose $\hat E\subseteq E$ is a union of
spots in $E$, which coincides with a union of spots in $E'$.  Then
there is an isotopy rel $\hat E$ of  $(K,E)$ to $(K',E')$.
\end{lemma}
\begin{proof}
(a) We will use the following result.  (It is stated, for example, in
\cite{DM:MappingSurvey}.)  The (mapping class) group of isotopy
classes of automorphisms of $(P,\bdry P)$ which map each component of
$\bdry P$ to itself is trivial.  It is then also easy to verify that
the mapping class group of $(P,\bdry P)$ rel $\bdry P$ is trivial.

Suppose now that $(K,E)$ and $(K',E')$ are two different standard
spotted balls in $(P,\bdry P)$.  We will construct an automorphism $f$
of $(P,\bdry P)$ by defining it first on $K$, then extending to the
complement.  We choose a homeomorphism $f:(K,E)\to (K',E')$. Next, we
extend to $\bdry P- E$ such that $f(\bdry P- E)=\bdry P-E'$.  Since
$(K,E)$ is standard, the sphere $\fr(K)\cup (\bdry P-E)$ bounds a ball
$K^c$ (with interior disjoint from $\intr(K)$), and similarly
$\fr(K')\cup (\bdry P-E')$ bounds a ball $(K')^c$.  Thus we can extend
$f$ so that the ball $K^c$ is mapped homeomorphically to $(K')^c$.
Then $f$ is isotopic to the identity, whence $(K,E)$ is isotopic to
$(K',E')$.

(b) To prove this, we adjust the isotopy of (a) in a collar of $\bdry
P$.  To see that such an adjustment is possible, note that we may
assume the collar is chosen sufficiently small so that by
transversality the isotopy of (a) is a ``product isotopy" in the
collar, i.e., has the form $h_t\times \text{id}$ for an isotopy $h_t$
of $\bdry P$.  The adjustment of the isotopy in the collar is achieved
using isotopy extension.
\end{proof}

\begin{proof} [Proof of Theorem \ref{CanonicalHeegaardTheorem}]

Any complete system of essential spheres $\cal S$ gives a symmetric
Heegaard splitting as follows.  Each holed 3-sphere $(P_i,\bdry P_i)$
obtained by cutting $M$ on $\cal S$ (where $P_i$ has $r_i$ holes, say)
can be expressed as a union of two standard $r_i$-spotted balls
$(K_i,E_i)$ and $(K_i^c,E_i^c)$. We assemble the spotted balls
$(K_i,E_i)$ by isotoping and identifying the two spots on opposite
sides of each essential sphere of the system $\cal S$, extending the
isotopy to $K_i$ and to $K_i^c$.  The result is a handlebody
$H=\cup_iK_i$ with a complementary handlebody $H^c=\cup_iK_i^c$.

We claim that $H$ (depending on $\cal S$) is unique up to isotopy in
$M$.  First replace each sphere $S$ of $\cal S$ by $\bdry N(S)$,
i.e. by two parallel copies of itself.  The result is that no holed
ball obtained by cutting $M$ on $\cal S$ has two different boundary
spheres identified in $M$.  Suppose $H'$ is constructed in the same
way as $H$ above.  Inductively, we isotope $H'\cap P_i=K_i'$ to $K_i$
using Lemma \ref{SpottedBallStandardLemma}, extending the isotopy to
$H'$.  Having isotoped $K_i'$ to $K_i$ for $i<k$, we isotope $K_k'$ to
$K_k$ in $P_k$ rel $\cup_{i<k} K_i'\cap \bdry P_k$.  This means we
perform the isotopy rel spots on the boundary of $P_k$ where previous
$K_i'$'s have already been isotoped to coincide with $K_i$'s.  Now it
is readily verified that the $H$ associated to the non-primitive $\cal
S$ obtained by replacing each sphere with two copies of itself is the
same as the $H$ associated to the original $\cal S$.

Next we must show that $H$ does not depend on the choice of a complete
system $\cal S$ of essential spheres.  By Lemma \ref{SphereSpaceLemma}
it is possible to find a sequence $\{{\cal S}_i\}$ of complete systems
of essential spheres such that any two successive systems differ by
the addition or removal of spheres.  It is then easy to verify that
after each move, the handlebody $H_i$, constructed from ${\cal S}_i$
as $H$ was constructed from $\cal S$ above, will be the same up to
isotopy.  For example, when removing an essential sphere from the
system ${\cal S}_i$ to obtain ${\cal S}_{i+1}$, two spotted balls are
identified on a disc on the removed sphere, yielding a new spotted
ball in a new holed sphere, and $H_{i+1}=H_i$.  The reverse process of
adding an essential sphere must then also leave $H_{i+1}=H_i$
unchanged.
\end{proof}

We offer an alternative simpler proof of the theorem.  Although it is
simpler, this approach has the disadvantage that it does not appear to
work for our more general result about canonical splittings of mixed
bodies.

\begin{proof}[Alternate proof of Theorem \ref{CanonicalHeegaardTheorem}]
Let $M=H_1\cup_{\partial H_1=\partial H_2} H_2$ and
$M=H_1'\cup_{\partial H_1'=\partial H_2'} H_2'$ determine two
symmetric Heegaard splittings of $M$. Our goal is to build an
automorphism $h\colon M\to M$ taking $H_1$ to $H_1'$ (and hence $H_2$
to $H_2'$) which is isotopic to the identity.

Choose any automorphism $f\colon H_1\to H_1'$. It doubles to an
automorphism $Df\colon M\to M$ which takes the first Heegaard
splitting to the second. Now consider $Df_*\colon \pi_1(M)\to
\pi_1(M)$. It is a known fact that any automorphism of a free group is
realizable by an automorphism of a handlebody (e.g., \cite{HG:64}).
Let $g\colon H_1\to H_1$ be such that $g_*=(Df_*)^{-1}$ and double
$g$, obtaining $Dg\colon M\to M$. Then it is clear that $h= (Df\circ
Dg)\colon M\to M$ takes the first Heegaard splitting to the
second. Also, $h_*=\text{Id}_{\pi_1(M)}$.

Now let $\mathcal{S}$ be a complete system of spheres symmetric with
respect to the splitting $M=H_1\cup H_2$. From
$h_*=\text{Id}_{\pi_1(M)}$ follows that $h$ is isotopic to a
composition of twists on spheres of $\mathcal{S}$
\cite{FL:Spheres}. But such a composition preserves $H_1$, $H_2$ (a
twist on such a sphere restricts to each handlebody as a twist on a
disc), therefore $h$ is isotopic to an automorphism $h'$ preserving
$H_1$, $H_2$. The isotopy from $h$ to $h'$ then takes the second
splitting to the first one.
\end{proof}

Our next task is to describe canonical Heegaard splittings of mixed
bodies.  Until further notice, $M$ will be a mixed body without holes.
We obtain a canonical splitting of a holed mixed body simply by
capping the holes.

\begin{thm}\label{CanonicalHeegaardTheorem2}
Let $M$ be a mixed body. There is a surface $F$ bounding a handlebody
$H$ in $M$ with the property that every curve on $F$ which bounds a
disc in $H$ also bounds a disc in $M-\intr (H)$, or it is boundary
parallel in $M$. $F$ and $H$ are unique up to isotopy.  $M-\intr(H)$
is a compression body.  Also $\pi_1(H)\to \pi_1(M)$ is an isomorphism.
\end{thm}

The proof is entirely analogous to the proof of Theorem
\ref{CanonicalHeegaardTheorem}, but we must first generalize the
definitions given above, for $M$ as sphere body, to definitions for
$M$ a mixed body.  Let $M$ be a connected sum of handlebodies $H_i$,
$i=1,\ldots, k$, $H_i$ having genus $g_i$, and $N$, which is a
connected sum of $n$ copies of $S^2\times S^1$.

An {\it essential reducing system} or {\it essential system} for $M$
is the isotopy class of a system of disjointly embedded spheres and
discs with the property that none of the spheres bounds a ball and
none of the discs is $\bdry$-parallel.   The system is {\it primitive}
if no two surfaces of the system are isotopic.  A system $\cal S$ is
said to be {\it complete} if cutting $M$ on $\cal S$ yields a
collection of holed balls $P_i$ satisfying an additional condition
which we describe below.  Cutting on $\cal S$ yields components of the
form $(P_i,S_i)$, where $S_i$ is the union of discs and spheres in
$\bdry P_i$ corresponding to the cutting locus:  if $P_i$ is viewed as
an immersed submanifold of $M$, then $\intr(S_i)$ is mapped to the
interior of $M$.  For the system to be complete we also require that
$S_i$ contain all but at most one boundary sphere of $P_i$.  Thus at
most one boundary component contains discs of $S_i$. There is a cell
complex of weighted primitive systems of essential spheres and discs
in $M$ constructed as follows.  For any primitive essential system
$\cal S$ in $M$, the space contains a simplex corresponding to
projective classes of systems of non-negative weights on the spheres.
Two simplices corresponding to primitive essential systems ${\cal
S}_1$ and ${\cal S}_2$ are identified on the face corresponding to the
maximal primitive essential system ${\cal S}_3$ of spheres common to
${\cal S}_1$ and ${\cal S}_2$.  We call the resulting space the {\it
reducing space} of $M$.  In the reducing space, the union of open
cells corresponding to complete primitive systems is called the {\it
complete reducing space} of $M$

\begin{lemma}\label{ReducingSpaceLemma}
The complete reducing space of a mixed body $M$ is connected.  In
other words, it is possible to pass from any primitive complete
essential system to any other by moves which replace a primitive
complete essential system by a primitive complete subsystem or a
primitive complete essential supersystem.
\end{lemma}

\begin{proof}
The proof is, of course, similar to that of Lemma
\ref{SphereSpaceLemma}. Suppose ${\cal S}={\cal S}_0$ and ${\cal S}'$
are two primitive complete systems.  We will use the term ``reducing
surface" or just ``reducer" to refer to a properly embedded surface
which is either a sphere or a disc.

Consider ${\cal S}\cap{\cal S}'$ and choose a curve of intersection
innermost in ${\cal S}'$.  Remove the curve of intersection by isotopy
of ${\cal S}$ if possible.  Otherwise, surger $\cal S$ on a disc
(half-disc) bounded in ${\cal S}'$ by the curve to obtain ${\cal
S}_2$.   Now ${\cal S}_2$ is obtained from ${\cal S}={\cal S}_0$ by
removing one reducer, the surgered reducer $S_0$ of ${\cal S}_0$, say,
and adding two reducers, $S_{21}$ and $S_{22}$, which result from the
surgery.  We can do this in two steps, first adding the two new
reducers to obtain ${\cal S}_1$, then removing the surgered reducer to
obtain ${\cal S}_2$.  If one of the reducers $S_{21}$ and $S_{22}$ is
inessential, this shows that the curve of intersection between ${\cal
S}={\cal S}_0$ and ${\cal S}'$ could have been removed by isotopy.
Note that surgery on arcs is always between two different discs of the
system, but surgery on closed curves may involve any pair of reducers
(disc and sphere, two discs, or two spheres).  We note that $S_{21}$
or $S_{22}$, or both, might be isotopic to reducers of ${\cal
S}_0$. If the new reducers are essential, we claim that ${\cal S}_1$
and ${\cal S}_2$ are complete systems.  It is clear that ${\cal S}_1$
cuts $M$ into holed balls, since it is obtained by adding two reducers
$S_{21}$ or $S_{22}$ to ${\cal S}_0$.  The property that cutting on
reducers yields pairs $(P_i,S_i)$ with at most one boundary component
of $P_i$ not in $S_i$ is preserved.  This is clearly true when adding
a sphere.  When adding a disc, a sphere of $\bdry P_i$ not contained
in $S_i$ is split, and so is $P_i$, yielding two holed balls with
exactly one boundary sphere not contained in $S_i$.  We then remove
$S_0$ from ${\cal S}_1$ to obtain ${\cal S}_2$.  There is, of course a
danger that removing an essential reducer from an essential system
${\cal S}_1$ yields a system which does not cut $M$ into holed balls,
or does not preserve the other completeness property.  As before, if
the surgery is on a disc, then the effect on the complementary
$(P_i,S_i)$'s is first to cut a $P_i$ on a disc with boundary in
$S_i$.  Clearly this has the effect of cutting $P_i$ into two holed
balls.  If the surgery is on a disc of $S_i$, we obtain one new holed
ball $(P_j,S_j)$ with all of its boundary in  $S_j$, and one holed
ball with one boundary component not contained in $S_j$. Then we
attach a 2-handle to an $S_i$ of another $(P_i,S_i)$, which, as
before, yields the required types of holed balls. If the surgery is on
a half-disc, we first cut a holed ball $(P_i,S_i)$ on a half disc,
i.e. we cut on a disc intersecting the distinguished component of
$\bdry P_i$ which is not in $S_i$.  This yields the required types of
holed balls.  Then we attach a half 2-handle whose boundary intersects
$S_i$ in a single arc to a $(P_i,S_i)$, which has the effect of
splitting a disc of $S_i$ into two discs. Thus, ${\cal S}_2$ cuts $M$
into holed balls of the desired type.  If ${\cal S}_2$ is not
primitive, i.e. contains two isotopic reducers, we can remove one,
until we are left with a primitive complete essential system.

Inductively, removing curves of intersection of ${\cal S}_{k}$ with
${\cal S}'$ either by isotopy or by surgery on an innermost disc
(half-disc) in ${\cal S}'$ as above, we must obtain a sequence ending
with ${\cal S}_{n-2}$ disjoint from ${\cal S}'$.  Next, we add the
reducers of ${\cal S}'$ not in ${\cal S}_{n-2}$ to ${\cal S}_{n-2}$ to
obtain ${\cal S}_{n-1}$, which is still complete.  Finally we remove
all reducers from ${\cal S}_{n-1}$ except those of ${\cal S}'$ to get
${\cal S}_{n}={\cal S}'$.
\end{proof}

Let $P$ be a 3-sphere with $r\ge 1$ holes and let $S\subseteq \bdry P$
be a union of discs and spheres such that each of $r-1\ge 1$
components of $\bdry P$ is a sphere of $S$, and the remaining sphere
of $\bdry P$ either is also in $S$ or contains a non-empty set of
discs of $S$.  ($(P,S)$ is the type of pair obtained by cutting on a
complete system.)  If $S$ contains discs, a {\it standard} $r$-spotted
ball $(K,E)$ in $(P,S)$ is one with the property that each component
of $S$ contains a component of $E$, and the complement
$\overline{P-K}$ is a spotted product, $(V\times I, F)$, where $F$ is
a union of discs in $V\times 1$.  This means that $\overline{P-K}$ has
the form $V\times I$, where $V$ is a planar surface, $V\times 0$ is
mapped to $\bdry P-\intr(S)$, $\bdry V\times I$ is mapped to the
annular components of $S-\intr(E)$, and exactly one disc component of
$F$ is mapped to each sphere of $S$.  We have already defined what it
means for $(K,E)$ to be standard in $(P,S)$ if $\bdry P=S$.  Note that
the spots of $F$ are in one-one correspondence with the sphere
components of $S$, and each is mapped to the complement of $E$ in its
sphere component.

\begin{lemma} \label{SpottedBallStandardLemma2}
Let $P$ be a 3-sphere with $r\ge 1$ holes and let $S\subseteq \bdry P$
be a union of discs and spheres such that each of $r-1\ge 1$
components of $\bdry P$ is a sphere of $S$, and the remaining sphere
of $\bdry P$ either is also in $S$ or contains a non-empty set of
discs of $S$

(a) A standard embedded $r$-spotted ball $(K,E)$ in $(P,S)$ is unique
up to isotopy of pairs.

(b) Suppose $(K,E)$ and $(K',E')$ are two different standard spotted
balls in $(P,S)$ and suppose $\hat E\subseteq E$ is a union of spots
in $E$, which coincides with a union of spots in $E'$.  Then there is
an isotopy rel $\hat E$ of  $(K,E)$ to $(K',E')$.
\end{lemma}
\begin{proof}
(a) The (mapping class) group of isotopy classes of automorphisms of
$(P,S)$ which map each component of $S$ to itself is trivial.  It is
then also easy to verify that the mapping class group of $(P,S)$ rel
$S$ is trivial.

We may assume that $S$ contains discs, otherwise we have already
proved the result.  Suppose now that $(K,E)$ and $(K',E')$ are two
different standard spotted balls in $(P,S)$.  We will construct an
automorphism $f$ of $(P,S)$ by defining it first on $K$, then
extending to the complement.  Isotope $E'$ to $E$ extending the
isotopy to $K'$.  We choose a homeomorphism $f:(K,E)\to (K',E')$
sending a disc of $E$ to the corresponding disc of $E'$. Next, we
extend to $S- E$ such that $f(S- E)=S-E'$.  Since $(K,E)$ is standard,
its complement has a natural structure as a spotted product, and the
complement of $K'$ has an isomorphic structure.  Thus we can extend
$f$ so that the complementary spotted product $K^c$ is mapped
homeomorphically to $(K')^c$.  Then $f$ is isotopic to the identity,
whence $(K,E)$ is isotopic to $(K',E')$.

(b) To prove (b) we adjust the isotopy of (a) in a collar of $\hat E$
in $K$.
\end{proof}

\begin{proof} {(\it Theorem \ref{CanonicalHeegaardTheorem2})}
Any complete system $\cal S$ for $M$ gives a Heegaard splitting as
follows.  In each holed 3-sphere $(P_i,S_i)$ obtained by cutting $M$
on $\cal S$ we insert an $s_i$-spotted ball $(K_i,E_i)$, where $s_i$
is the number of components of $S_i$. Each component of $S_i$ contains
exactly one disc $E_i$.  We assemble the spotted balls $(K_i,E_i)$ by
isotoping and identifying the two spots on opposite sides of each
essential sphere or disc of the system $\cal S$, extending the isotopy
to $K_i$.  The result is a handlebody $H=\cup_iK_i$ with a
complementary handlebody $H^c=\cup_iK_i^c$.

We claim that $H$ (depending on $\cal S$) is unique up to isotopy in
$M$.  First replace each sphere or disc $S$ of $\cal S$ by $\bdry
N(S)$, i.e. by two parallel copies of itself.  The result is that no
$(P_i,S_i)$ obtained by cutting $M$ on $\cal S$ has two different disc
or sphere components of $S_i$ identified in $M$.  Suppose $H'$ is
constructed in the same way as $H$ above.  Inductively, we isotope
$H'\cap P_i=K_i'$ to $K_i$ using Lemma
\ref{SpottedBallStandardLemma2}, extending the isotopy to $H'$.
Having isotoped $K_i'$ to $K_i$ for $i<k$, we isotope $K_k'$ to $K_k$
in $P_k$ rel $\cup_{i<k} K_i'\cap \bdry P_k$.  This means we perform
the isotopy rel spots on the boundary of $P_k$ where previous $K_i'$'s
have already been isotoped to coincide with $K_i$'s.  Now it is
readily verified that the $H$ associated to the non-primitive $\cal S$
obtained by replacing each sphere or disc with two copies of itself is
the same as the $H$ associated to the original $\cal S$.

Next we must show that $H$ does not depend on the choice of a complete
system $\cal S$.  By Lemma \ref{ReducingSpaceLemma} it is possible to
find a sequence $\{{\cal S}_i\}$ of complete systems such that any two
successive systems differ by the addition or removal of reducers.  It
is then easy to verify that after each move, the handlebody $H_i$,
constructed from ${\cal S}_i$ as $H$ was constructed from $\cal S$
above, will be the same up to isotopy.  For example, when removing an
essential reducer from the system ${\cal S}_i$ to obtain ${\cal
S}_{i+1}$, two spotted balls are identified on a disc on the removed
reducing surface, yielding a new spotted ball in a new holed sphere,
and $H_{i+1}=H_i$.  The reverse process of adding an essential reducer
must then also leave $H_{i+1}=H_i$ unchanged.
\end{proof}

The canonical splitting surface for $M$ a sphere body or a mixed body
yields an automatic reducing surface for any automorphism of $f:M\to
M$.  However, this reducing surface is non-rigid.  There are
automorphisms of $M$ isotopic to the identity which induce
automorphisms not isotopic to the identity on the splitting
surface. The non-rigidity (without our terminology) is well-known, but
rather than giving the example described in
\cite{FL:ConnectedS2TimesS1}, we describe a simpler example, which in
some sense is the source of the phenomenon.

\begin{example}\label{NonUniqueExample}
Let $M$ be a sphere body and let $S$ be a non-separating sphere in $M$
intersecting the canonical splitting surface in one closed curve.  Let
$H$ be one of the handlebodies bounded by the splitting surface, so
the splitting surface is $\bdry H$.  We describe an isotopy of $H$ in
$M$ which moves $H$ back to coincide with itself and induces an
automorphism not isotopic to the identity on $\bdry H$.  In other
words, we describe $f:M\to M$ isotopic to the identity, such that
$f|_{\bdry H}$ is not isotopic to the identity.  Figure
\ref{Aut3UhrMove} shows the intersection of $H$ with $S$.  We have
introduced a kink in a handle of $H$ such that the handle runs
parallel to $S$ for some distance.  The ends of the handle project to
(neighborhoods of) distinct points while the handle projects to a
neighborhood of an arc $\alpha$ in $S$.  There is an isotopy of
$\alpha$ in $S$ fixing endpoints, so the arc rotates about each point
by an angle of $2\pi$ and the arc sweeps over the entire 2-sphere,
returning to itself.  We perform the corresponding isotopy of the
1-handle in $H$.  Calling the result of this isotopy $f$, we note that
$f|_{\bdry H}$ is a double Dehn twist as shown.  As one can easily
imagine, by composing automorphisms of this type on different spheres
intersecting $H$ in a single meridian, one can obtain $f$ such that
$f$ restricted to $\bdry H$ is a more interesting automorphism, e.g. a
pseudo-Anosov, see \cite{FL:ConnectedS2TimesS1}.

\begin{figure}[ht]
\centering
\psfrag{S}{\fontsize{\figurefontsize}{12}$S$}\psfrag{H}{\fontsize{\figurefontsize}{12}$H$}
\psfrag{a}{\fontsize{\figurefontsize}{12}$\alpha$}

\scalebox{1.0}{\includegraphics{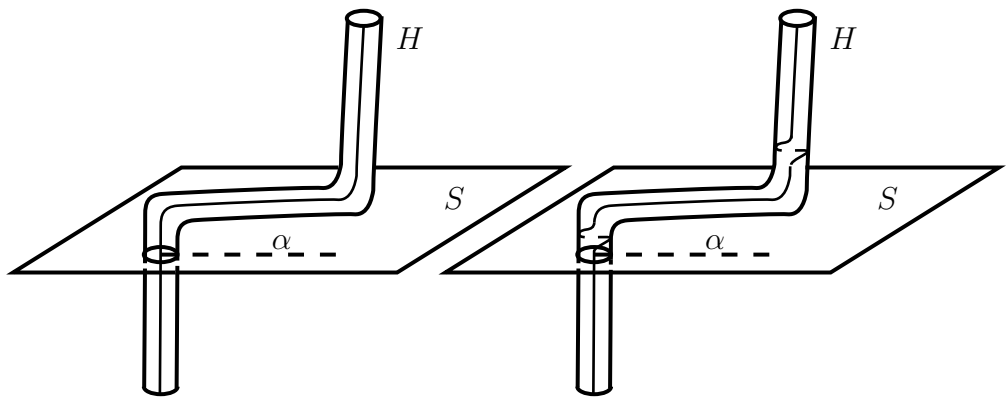}} \caption{\small
Non-rigidity of canonical splitting.} \label{Aut3UhrMove}
\end{figure}
\end{example}

We shall refer to an isotopy of the type described above as a {\it
double-Dehn isotopy}.  If, as above,  $M$ is a sphere body or
mixed body, and $f\from M\to M$ is an automorphism, we shall use
double-Dehn isotopies to choose a good representative of the
isotopy class of $f: M\to M$ preserving the canonical splitting.

\section{Reducing surfaces} \label{Reducing}

The purpose of this section is to define reducing surfaces in the
general setting of automorphisms of compact 3-manifolds. In addition,
we review reducing surfaces which arise in the setting of
automorphisms of handlebodies and compression bodies, \cite{UO:Autos},
using the somewhat more inclusive definitions given here.

We will also consider manifolds with some additional structure; we
will consider holed or spotted manifolds.

Suppose $M$ is a compact 3-manifold without sphere boundary
components, and $f:M\to M$ is an automorphism.  If $\mathcal{S}$ is an
$f$-invariant (up to isotopy) union of essential spheres, then
$\mathcal{S}$ is a special kind of reducing surface.  We can decompose
$M$ on $\mathcal{S}$ to obtain $M|\mathcal{S}$ and cap boundary
spheres with balls to obtain a  {\it manifold $(\hat M, {\cal D})$
with spots in the interior}, where $\mathcal{D}$ denotes the union of
the capping balls, and an induced automorphism $(\hat M, {\cal D})\to
(\hat M, {\cal D})$.   We could equally well not cap boundary spheres
after decomposing on $\mathcal{S}$ and consider holed manifolds rather
than spotted manifolds.

\begin{definition}   Let
$(M,\mathcal{D})$ be a spotted 3-manifold  with spots in the interior
($\mathcal{D}$ a union of 3-balls disjointly embedded in the interior
of $M$).  We say a sphere is {\it essential} in $(M,\mathcal{D})$ if
it is embedded in $M-\mathcal{D}$, it does not bound bound a ball in
$M-\mathcal{D}$, and it is not isotopic in $M-\mathcal{D}$ to a
component of $\bdry \mathcal{D}$.   The spotted manifold is {\it
irreducible} if it has no essential spheres.
\end{definition}

Clearly a connected spotted manifold $(M,\mathcal{D})$  is irreducible
if and only if the underlying manifold $M$ is irreducible and
$\mathcal{D}$ is empty or consists of a single ball.

\begin{definition}  Let $(M,\mathcal{D})$ be a spotted manifold with spots in the interior.  Let $f:(M,\mathcal{D})\to (M,\mathcal{D})$
be an automorphism of the spotted manifold.  Then $f$ is {\it
reducible on spheres} if there is an $f$-invariant union $\mathcal{S}$
of disjoint essential spheres in $(M,\mathcal{D})$.  $\cal S$ is
called a {\it sphere reducing surface}.  {\it Decomposing} on
$\mathcal{S}$ and capping with balls gives a new spotted manifold with
an induced automorphism.  The caps of duplicates of spheres of the
reducing surface are additional spots on the decomposed manifold.
\end{definition}

Clearly by repeatedly decomposing on sphere reducing surfaces we
obtain automorphisms of manifolds with interior spots which do not
have sphere reducing surfaces.

Automorphisms of a spotted 3-manifold $(M,\mathcal{D})$ (with
3-dimensional spots) induce automorphisms of the underlying spotless
3-manifold $M$, and in the cases of interest to us it is relatively
easy to understand the relationship between the two.  For this reason,
we can henceforth consider only 3-manifolds without 3-dimensional
spots and without sphere reducing surfaces.

\begin{definition} Suppose $f:M\to M$ is an automorphism,
without sphere reducing surfaces, of a 3-manifold.  Then a closed surface
$F\embed M$ which is $f$-invariant up to isotopy is called a {\it
closed reducing surface} provided $M|F$ is irreducible.  The reducing
surface $F$ is {\it rigid} if $f$ induces uniquely an automorphism of
the decomposed manifold $M|F$.
\end{definition}

In particular, the canonical splitting surface for $M$ a sphere body
or mixed body is a (non-rigid) reducing surface for any automorphism
of $M$.

The existence of any closed reducing surface for $f:M\to M$ (without
sphere components) gives a decomposition into automorphisms of
compression bodies and manifolds with incompressible boundary, as the
following theorems show.  We give a simple version for closed
manifolds separately to illustrate the ideas.

\begin{thm} \label{ReducingDecompositionTheorem}
Let $M$ be a closed reducible manifold, $f:M\to M$ an automorphism,
and suppose $F$ is a closed reducing surface, with no sphere
components.  Then either

i) $F$ is incompressible,

ii) there is a reducing surface $F'$ such that the reducing surface
$F\cup F'$ cuts $M$ into a non-product compression body $Q$ (possibly
disconnected) and a manifold $M'$ with incompressible boundary, or

iii) $F$ is a Heegaard surface.
\end{thm}

\begin{proof}  We suppose $M$ is a closed manifold, $f:M\to M$ is an automorphism
and suppose $F$ is a closed reducing surface.  Cutting $M$ on $F$ we
obtain $M|F$ which is a manifold whose boundary is not in general
incompressible.  Let $Q$ be the characteristic compression body in
$M|F$.  This exists because $M|F$ is irreducible by the definition of
a reducing surface.  Then $f|F$ preserves $Q$ and $F'=\bdry_i Q$ is
another reducing surface, unless it is empty.  In case $\bdry_iQ$ is
empty, we conclude that $Q$ must consist of handlebodies, and we have
a Heegaard decomposition of $M$.  Otherwise, letting
$M'=\overline{M-Q}$, $M'$ has incompressible boundary.
\end{proof}

Here is a more general version:

\begin{thm} \label{ReducingDecompositionTheorem2}
Let $M$ be a reducible manifold,  $f:M\to M$ an automorphism, and
suppose $F$ is a closed reducing surface, with no sphere components.
Then either

i) $F$ is incompressible,

ii) there is a reducing surface $F'$ such that the reducing surface
$F\cup F'$ cuts $M$ into a non-product compression body $Q$ (possibly
disconnected) and a manifold $M'$ with incompressible boundary, or

iii) $F$ is a Heegaard surface dividing $M$ into two compression
bodies.
\end{thm}

A priori, there is no particular reason why every automorphism $f:M\to
M$ without sphere reducing surfaces should have a closed reducing
surface.  Recall we posed Question \ref{ReducingQuestion}.  In any
case, if $M$ is reducible, we deal only with automorphisms for which
this is true. Therefore we henceforth consider only automorphisms of
irreducible 3-manifolds of the types resulting from the previous two
theorems:  Namely, we must consider automorphisms of irreducible,
$\bdry$-irreducible manifolds, and we must consider automorphisms of
compression bodies (and handlebodies).

It might be interesting to learn more about all possible reducing
surfaces associated to a given automorphism of a given manifold,
especially for answering Question \ref{ReducingQuestion}.  In our
classification, however, we often use reducing surfaces guaranteed by
the topology of the 3-manifold: reducing surfaces coming from
F. Bonahon's characteristic compression body and the
Jaco-Shalen-Johannson characteristic manifold, and reducing surfaces
coming from the canonical Heegard splittings of sphere bodies or mixed
bodies.  To analyze automorphisms of compression bodies, one uses
reducing surfaces not guaranteed by the topology, see \cite {UO:Autos}.

The definition of reducing surfaces for irreducible,
$\bdry$-irreducible manifolds should be designed such that the
frontier of the characteristic manifold is a reducing surface.  The
following is a natural definition.

\begin{definition} Suppose $f:M\to M$ is an automorphism of an irreducible, $\bdry$-irreducible compact manifold. Then an invariant
surface $F$ is a {\it reducing surface} if it is incompressible and
$\bdry$-incompressible and not boundary-parallel.
\end{definition}

It only remains to describe reducing surfaces for automorphisms of
handlebodies and compression bodies, and this has already been done in
\cite{UO:Autos}.  We give a review here, taken from that paper.

\begin{defns}  We have already defined a {\it spotless compression body} $(Q,V)$:
The space $Q$ is obtained from a disjoint union of balls and a product
$V\times I$ by attaching 1-handles to the boundaries of the balls and
to $V\times 1$; the surface $V\embed \bdry Q$ is the same as $V\times
0$.  We require that $V$ contain no disc or sphere components.  When
$Q$ is connected and $V=\emptyset$, $(Q,V)$ is a handlebody or ball.
If $(Q,V)$ has the form $Q=V\times I$ with $V\times 0=V$, then $(Q,V)$
is called a {\it spotless product compression body} or a {\it trivial
compression body}.  As usual, $\bdry_eQ=\bdry Q-\intr (V)$ is the {\it
exterior boundary}, even if $V=\emptyset$.

A {\it spotted compression body} is a triple $(Q,V,\mathcal{D})$ where
$(Q,V)$ is a spotless compression body and $\mathcal{D}\ne \emptyset$
denotes a union of discs or ``spots" embedded in $\bdry_e Q=\bdry
Q-\intr(V)$.

A {\it spotted product} is a spotted compression body
$(Q,V,\mathcal{D})$ of the form $Q=V\times I$ with $V=V\times 0$ and
with $\mathcal{D}\ne \emptyset$ a union of discs embedded in $V\times
1$.

A {\it spotted ball} is a spotted compression body whose underlying
space is a ball.  It has the form $(B,\mathcal{D})$ where $B$ is a
ball and $\mathcal{D}\ne \emptyset$ is a disjoint union of discs in
$\bdry B=\bdry_eB$.

An {\it $I$-bundle pair} is a pair $(H,V)$ where $H$ is the total
space of an $I$-bundle $p\from H \to S$ over a surface $S$ and $V$ is
$\bdry_iH$, the total space of the associated $\bdry I$-bundle.

\end{defns}

\begin{com}
A spotless compression body or $I$-bundle pair is an example of a {\it
Haken pair}, i.e. a pair $(M,F)$ where $M$ is an irreducible
3-manifold and $F\subset \bdry M$ is incompressible in $M$.
\end{com}

\begin{definition} If  $(Q, V)$ is a compression body, and $f\from(Q, V)\to (Q, V)$ is an automorphism, then an
$f$-invariant surface $(F,\bdry F)\embed (Q,\bdry Q)$, having no
sphere components, is a {\it reducing surface for $f$}, if:

i) $F$ is a union of essential discs, $\bdry F\subset \bdry_eQ$,

or no component of $F$ is a disc and:

ii) $F$ is the interior boundary of an invariant compression body
$(X,F)\embed (Q,V)$, with $\bdry_eX\embed \bdry_eQ$ but $F$ is not
isotopic to $V$, or

iii) $F\subset \bdry Q$, $F$ is incompressible, $F$ strictly contains
$V$, and $(Q,F)$ has the structure of an $I$-bundle pair, or

iv) $F$ is a union of annuli, each having one boundary component in
$\bdry_eQ$ and one boundary component in $\bdry_iQ$

In case i) we say $F$ is a {\it disc reducing surface}; in case ii) we
say $F$ is a {\it compressional reducing surface}; in case iii) we say
$F$ is a {\it peripheral reducing surface}; in case iv) we say $F$ is
an {\it annular reducing surface}.

In case i) we {\it decompose} on $F$ by replacing $(Q,V)$ by the {\it
spotted compression body $(M|F,V,\mathcal D)$}, and replacing $f$ by
the induced automorphism on this spotted compression body. ($\cal D$
is the union of duplicates of the discs in $F$.)  An automorphism $g$
of a spotted compression body $(Q,V,\mathcal D)$ can be {\it
decomposed} on a {spot-removing reducing surface}, which is $\bdry_e
Q$ isotoped to become a properly embedded surface.  This yields an
automorphism of a {\it spotted product} and an automorphism of a
spotless compression body.

In case ii) we {\it decompose} on $F$ by cutting on $F$ to obtain
$Q|F$, which has the structure $(X,F)\cup (Q',V')$ of a compression
body, and by replacing $f$ by the induced automorphism.  The duplicate
of $F$ not attached to $X$ belongs to $\bdry_eQ'$.  We throw away any
trivial product compression bodies in $(Q',V')$.

In case iii) we {\it decompose on $F$} by replacing  $(Q,V)$ by
$(Q,F)$.

In case iv), cutting on $F$ yields another compression body with
induced automorphism.

For completeness, if $(H,V)$ is a spotless $I$-bundle pair, we define
reducing surfaces for automorphisms $f\from (H,V)\to (H,V)$.  Any
$f$-invariant non-$\bdry$-parallel union of annuli with boundary in
$V$ is such an {\it annular reducing surface}, and these are the only
reducing surfaces.
\end{definition}

Non-rigid reducing surfaces are quite common, as the following example
shows.

\begin{example}
Let $M$ be a mapping torus of a pseudo-Anosov automorphism $\phi$,
$M=S\times I/ \sim$ where $(x,0)\sim (\phi(x),1)$ describes the
identification.  Let $\tau$ be an invariant train track for the
unstable lamination of $\phi$ in $S=(S\times 0)/\sim$.  Let $f:M\to M$
be the automorphism induced on $M$ by $\phi\times \text{Identity}$.
Then $N(\tau)$ is invariant under $f$ up to isotopy and $F=\bdry
N(\tau)$ is a reducing surface for $F$.  Isotoping $\tau$ and $F=\bdry
N(\tau)$ through levels $S_t=S\times t$ from $t=0$ to $t=1$ and
returning it to its original position induces a non-trivial
automorphism on $N(\tau)$.  Hence $F$ is non-rigid.

It is easy to modify this $M$ by Dehn surgery operations such that
$F=\bdry N(\tau)$ becomes incompressible.
\end{example}

The results of Section \ref{Splitting} guarantee reducing surfaces in
the special case where $M$ is a mixed body without sphere boundary
components, i.e. a connected sum of $S^2\times S^1$'s and
handlebodies, and this is what we use in the classification.  The
canonical splitting surfaces of Section \ref{Splitting} are examples
of reducing surfaces such that applying Theorem
\ref{ReducingDecompositionTheorem} and Theorem
\ref{ReducingDecompositionTheorem2} yields a Heegaard splitting.

We state the obvious:

\begin{proposition} If $M$ is a mixed body,
then there is a canonical Heegaard decomposition $M=H\cup Q$ on a
reducing surface $F=\bdry H$, which is unique up to isotopy.  $F$ is a
non-rigid closed reducing surface for any automorphism $f:M\to M$.
\end{proposition}

\section{Adjusting automorphisms} \label{Adjusting}

In this section we will prove Theorem \ref{AdjustingTheorem} and
Corollary \ref{AdjustingCorollary}.

The adjusting automorphisms are constructed from {\em slide
automorphisms of $M$}. There are two types, one of which is as
defined in \cite{DM:MappingSurvey}, see the reference for details.
For a fixed $i$, construct $\check M$ by removing
$\intr(\hat{M}_i)$ from $M$ and capping the resulting boundary
sphere with a ball $B$. Let $\alpha\subseteq \check M$ be an
embedded arc meeting $B$ only at its endpoints. One can ``slide''
$B$ along $\alpha$ and back to its original location.  By then
replacing $B$ with $\hat M_i$ we obtain an automorphism of $M$
which is the identity except in a holed solid torus $T$ (a closed
neighborhood of $\alpha\cup S_i$ in $M-\intr(\hat M_i)$). Call
such an automorphism a {\em simple slide}. It is elementary to
verify that $\hat M_0 \cup T$ is a mixed body.

We also need another type of slide which is not described explicitly
in \cite{DM:MappingSurvey}. First, consider the classes of $\hat
M_i$'s which are homeomorphic. For each class, fix an abstract
manifold $N$ homeomorphic to the elements of the class.  For each
$\hat M_i$ in the class fix a homeomorphism $f_i\colon \hat M_i\to N$.

If $i\ne j$, and $\hat M_i$ is homeomorphic to $\hat M_j$,
construct $\check M$ by removing $\intr(\hat{M}_i\cup\hat{M}_j)$
and capping the boundary spheres $S_i$, $S_j$ with balls $B_i$,
$B_j$. Let $\alpha$, $\alpha'\subseteq\check M$ be two disjoint
embedded arcs in $\check M$ linking $S_i$ and $S_j$. One can slide
$B_i$ along $\alpha$ and $B_j$ along $\alpha'$ simultaneously,
interchanging them. Now remove the balls and re-introduce $\hat
M_i$, $\hat M_j$, now interchanged. Denote by $s$ the resulting
automorphism of $M$.  Isotope $s$ further to ensure that it
restricts to $\hat M_i$ as $f_j^{-1}\circ f_i$ (and to $\hat M_j$
as $f_i^{-1}\circ f_j$). Call such an automorphism an {\em
interchanging slide}. The support of an interchanging slide is a
closed neighborhood of $\hat M_i\cup\hat
M_j\cup\alpha\cup\alpha'$. Note that a sequence of interchanging
slides that fixes an $\hat M_i$ restricts to it as the identity.

\begin{remark}
An argument in \cite{DM:MappingSurvey} makes use of automorphisms
called {\em interchanges of irreducible summands}. Our interchanging
slides are a particular type of such automorphisms.
\end{remark}

\begin{proof} (Theorem \ref{AdjustingTheorem}.)
We start by noting that $\hat M_0$ is an essential holed sphere
body.

Let $\mathcal{S}=\bigcup_i S_i$ and consider $f(\mathcal{S})$. The
strategy is to use isotopies and both types of slides to take
$f(S_i)$ to $S_i$. If only isotopes are necessary we make
$R=\emptyset$, $\hslash=\id_M$ and $f=g$ in the statement.
Therefore suppose that some $f(S_i)$ is not isotopic to $S_i$. The
first step is to simplify $f(\mathcal{S})\cap\mathcal{S}$. The
argument of \cite{DM:MappingSurvey} yields a sequence of simple
slides (and isotopies) whose composition is $\sigma$ satisfying
$\sigma\circ f(\mathcal{S})\cap\mathcal{S}=\emptyset$. We skip the
details. We can then assume that $\sigma\circ
f(\mathcal{S})\subseteq\hat M_0$, otherwise some $S_i$ would bound
a ball. Note that $\sigma$ preserves $K\subseteq M$, where $K$
consists of $\hat M_0$ to which is attached a collection of
disjoint 1-handles $J\subseteq\cup\hat M_i$ along $\mathcal{S}$.
Each 1-handle consists of a neighborhood (in $\hat M_i$) of a
component of $\alpha\cap\hat M_i$, where $\alpha$ is a path along
which a simple slide is performed. Perturb these paths so that
they are pairwise disjoint.

Now consider $\sigma\circ f(S_i)\subseteq\hat M_0$. It must bound an
(once holed) irreducible summand on one side. It is then easy to check
that $\sigma\circ f(S_i)$ is isotopic to an $S_j$.  Indeed,
$\sigma\circ f(S_i)$ separates $\hat M_0$, bounding (in $\hat M_0$) on
both sides a holed sphere body. On the other hand $\sigma\circ f(S_i)$
must bound (in $M$) an irreducible summand.  Therefore in $\hat M_0$
it must bound an once-holed ball, whose hole corresponds to some
$S_j$. We can then make $\sigma\circ f(\mathcal{S})=\mathcal{S}$.

If $f(S_i)=S_j$, $i\neq j$ we use a slide to interchange $\hat M_i$
and $\hat M_j$. A composition $\iota$ of these interchanging slides
yields $\iota\circ\sigma\circ f(S_i)=S_i$ for each $i$.

Now consider the collection of 1-handles $J=(\cup\hat M_i)\cap K$,
where we recall that $K$ is the union of the support of $\sigma$
with $\hat M_0$ ($K$ is preserved by $\sigma$). We note that there
is some $n$ such that, for any $i\geq 1$, $\iota^n|_{\hat
M_i}=\text{Id}$. In particular, $\iota^n(J)=J$. By perturbing the
paths along which the simple slides of $\sigma$ are performed we
can assume that $\iota^k(J)\cap J=\emptyset$ for any $k$, $1\le
k<n$.  Therefore $\cup_k\iota^k(J)$ consists of finitely many
1-handles.  Let $R=K\cup\left(\cup_k\iota^k(J)\right)$. Clearly
$R$ is preserved by $\iota$. The restriction of $\iota$ to $\cup_i
\hat M_i$ is periodic, therefore its restriction to
$\overline{M-R}\subseteq(\cup_i \hat M_i)$ also is. Finally, $R$
is a mixed body. Either $\iota$ or $\sigma$ is not the identity by
hypothesis, therefore $R\neq B$, the 3-ball. It is not hard to see
that it is essential. Indeed, if $S\subseteq R$ is an essential
sphere then an argument similar to the one above yields a sequence
of slides $\tau$ which preserves $R$ and takes $S$ to
$\tau(S)\subseteq \hat M_0$. Essentiality of $S$ implies
essentiality of $\tau(S)$ in $R$. Therefore $\tau(S)$ is essential
in $\hat M_0$ (which is itself essential in $M$), and hence in
$M$. But essentiality of $\tau(S)$ in $M$ implies essentiality of
$S$ in $M$. It remains to consider $S$ a sphere boundary component
of $R$, but this is parallel to an $S_i\subseteq \mathcal{S}$,
which is essential since $M$ has no sphere boundary components.

Let $\hslash=\iota\circ\sigma$, $h=\hslash^{-1}$ and $g=\hslash\circ
f$. The following are immediate from the construction:

1) $R$ is preserved by $\hslash$ (and hence by $h$).

2) $\hslash|_{(\overline{M-R})}=\iota|_{(\overline{M-R})}$, which is
periodic.

3) $g$ preserves each $\hat M_i$.

The other conclusions of the theorem are clear.  \end{proof}

\begin{proof} (Corollary \ref{AdjustingCorollary}.)  We follow the proof of Theorem
\ref{AdjustingTheorem}.  In that proof, a composition of simple slides
and interchanging slides has the effect of moving $f(S_i)$ to $S_i$.
A 4-dimensional handlebody $P_0$ is obtained from a 4-ball $K_0$ by
attaching $\ell$ 1-handles.  The boundary is a connected sum of $\ell$
$S^2\times S^1$'s.  The 4-dimensional compression body $Q$ is obtained
from a disjoint union of the products $M_i\times I$ and from $P_0$ by
attaching a 1-handle to $P_0$ at one end and to $M_i\times 1$ at the
other end, one for each $i\ge 1$.  Then $\bdry_eQ$ is the connected
sum of the $M_i$'s and the sphere body $\bdry P_0$.    The cocore of
each 1-handle has boundary equal to an $S_i$, $i=1,\ldots, k$.  Each
$S_i$ bounds a 3-ball $E_i$ in $Q$.  A simple slide in $\bdry_eQ=M$
using some $S_i$ can be realized as the (restriction to the) exterior
boundary of an automorphism of $Q$ as follows:  Cut a 1-handle of $Q$
on the ball $E_i$.  This gives a lower genus compression body $\hat Q$
with two 3-ball {\it spots} on its exterior boundary.  We slide one of
the spots along the curve $\alpha$ associated to the slide, extending
the isotopy to $\hat Q$ and returning the spot to its original
position, then we reglue the pair of spots.  This defines an
automorphism of the compression body fixing the interior boundary and
inducing the slide on the exterior boundary.  An interchanging slide
can be realized similarly.

Suppose now that $f$ is given as in the proof of Theorem
\ref{AdjustingTheorem}, and we perform the slides in the theorem
so that $\iota\circ\sigma\circ f(S_i)=S_i$ for each $i$.  Capping
each $S_i$ in $\bdry \hat M_i$ with a disc $E_i$ to obtain $M_i$,
and capping the other copy of $S_i$ in the holed sphere body $\hat
M_0$ with a disc $E_i'$, we obtain a disjoint union of the $M_i$'s
and the sphere body $\bdry P_0=M_0$.  The automorphism $f$ induces
an automorphism of this disjoint union, regarded as spotted
manifolds. Thus we have an induced automorphism of a spotted
manifold $ f_s:((\sqcup_iM_i\sqcup \bdry P_0),\cup_i(E_i\cup
E_i'))\to((\sqcup_iM_i\sqcup \bdry P_0),\cup_i(E_i\cup E_i'))$.
($f_s$ for $f$ spotted.) If we ignore the spots, we have an
automorphism $f_a$ of $(\sqcup_iM_i\sqcup \bdry P_0)$.  ($f_a$ for
$f$ adjusted.)  These automorphisms yield an automorphism $f_p$ of
the disjoint union of the spotted products $M_i\times I$ with
spots $E_i$ in $M_i\times 1$ where $f_p$ restricts to $(M_i\times
1,E_i\times 1)$ as $f_s$ and to $M_i\times 0$ as $f_a$. ($f_p$ for
$f$ product.) We can also extend $f_s$ to $f_p$ on all of
$(P_0,\cup_i E_i')$ so that the restriction of $f_p$ to $(\bdry
P_0, \cup_iE_i')$ is $f_s$.  Now $f_p$ is an automorphism of a
disjoint union of products $M_i\times I$ with spots on one end of
the product, disjoint union $P_0$ with spots on $\bdry P_0$.
Gluing $E_i$ to $E_i'$, we obtain $Q$ and $f_p$ yields an
automorphism $f_c$ of the compression body.

Now it only remains to apply the inverses of the automorphisms of $Q$
realizing $\iota\circ\sigma$ on $\bdry_eQ$ as a composition of slides
of both kinds.  This composition is $\bar\hslash\inverse$, and we have designed $\bar\hslash$ such that
$\bar\hslash\inverse\circ f_c$ restricts to $f$ on $\bdry_eQ=M$.  Then take $\bar f=\bar\hslash\inverse\circ f_c$.

Note that the support of the slides is actually in $R$ of the previous
proof.
\end{proof}

\section{One-dimensional invariant laminations}\label{Laminations}

If $f\from H\to H$ is a generic automorphism of a handlebody, there is
an invariant 1-dimensional lamination for $f$ described in
\cite{UO:Autos}, associated to an incompressible 2-dimensional
invariant lamination of a type which we shall describe below, see
Proposition \ref{GoodLaminationProp}.  The invariant ``lamination" is
highly pathological, not even being embedded. One of our goals in this
section is to describe various types of 1-dimensional laminations in
3-manifolds, ranging from ``tame" to ``wild embedded" and to the more
pathological ``wild non-embedded" laminations which arise as invariant
laminations for automorphisms of handlebodies.  All invariant
1-dimensional laminations we construct are dual to (and associated
with) an invariant 2-dimensional lamination.  In general, it seems
most useful to regard the 1-dimensional and the 2-dimensional
invariant laminations of a generic automorphism as a dual pair of
laminations.

In the case of an automorphism of a compression body, only a
2-dimensional invariant lamination was described in
\cite{UO:Autos}. In this section we shall also describe invariant
1-dimension laminations for (generic) automorphisms of compression
bodies.

\begin{defns}
A 1-dimensional {\it tame lamination}  $\Omega\embed M$ in a
3-manifold $M$ is a subspace $\Omega$ with the property that $\Omega$
is the union of sets of the form $T\times I$ in {\it flat charts} of
the form $D^2\times I$ in $M$ for a finite number of these charts
which cover $\Omega$.  ($T$ is the {\it transversal} in the flat
chart.)  A tame lamination is {\it smooth} if the charts can be chosen
to be smooth relative to a smooth structure on $M$.

A 1-dimensional {\it abstract lamination} is a topological space which
can be covered by finitely many flat charts of the form $T\times I$,
where $T$ is a topological space, usually, in our setting, a subspace
of the 2-dimensional disc.

A 1-dimensional {\it wild embedded lamination} in a 3-manifold $M$ is
an embedding of an abstract lamination.  Usually wild will also mean
``not tame."

A 1-dimensional {\it wild non-embedded lamination} in a 3-manifold is
a mapping of an abstract lamination into $M$.

\end{defns}

We observe that in general the 1-dimensional invariant laminations
constructed in \cite{UO:Autos} for $f\from H \to H$, $H$ a handlebody,
are wild non-embedded, though they fail to be embedded only at
finitely many {\it singular points} in $H$.

We will now construct invariant laminations for generic automorphisms
of arbitrary compression bodies.

We will deal with a connected, spotless compression body $(H,V)$.
``Spotless" means in particular that $V$ contains no disc components.
The surface $\bdry H-\intr (V)$ will always be denoted $W$.  Also,
throughout the section, we assume $\text{genus}(H,V)>0$, where
$\text{genus}(H,V)$ denotes the number of 1-handles which must be
attached to $V\times I$ to build $Q$.

Let $(H_0, V_0) \subseteq  \intr (H)$  be a ``concentric compression
body" with a  product structure on its complement.  If the compression
body has the form $V\times [0,1]$ with handles attached to $V\times
1$, $H_0=V\times [0,1/2]$ with handles attached to $V\times 1/2$ in
such a way that $H-H_0\cup \text{fr}(H_0)$ has the structure of a
product $W\times [0,1]$.  Here $V=V\times 0$ and $W= W\times 1$.

We shall sketch the proof of the following proposition, which is
proven in \cite{UO:Autos}.

\begin {proposition} \label{InvtLamProp1}
Suppose $f:(H,V) \to (H, V)$ is generic automorphism of an
$3$-dimen\-sional compression body.  Then there is a 2-dimensional
measured lamination $\Lambda\subseteq \intr (H)$ with transverse
measure $\mu$, such that $f(\Lambda,\mu) = (\Lambda,\lambda \mu)$, up
to isotopy, for some stretch factor $\lambda > 1$.  The leaves of
$\Lambda$ are planes.  Further, $\Lambda$ ``fills" $H_0$, in the sense
that each component of $H_0-\Lambda$ is either contractible or
deformation retracts to $V$.
\end{proposition}

The statement of Proposition \ref{InvtLamProp1} will be slightly
expanded later, to yield Proposition \ref{InvtLamProp2}.

An automorphism $f: (H,V)\to (H,V)$  is called {\it outward expanding}
with respect to $(H_0,V)$ if $f|_{W \times I} = f|_{W\times 0} \times
h$, where  $h \from I \to I$ is a homeomorphism moving every point
towards the fixed point 1, so that $h(1)=1$ and $h(t)>t$ for all
$t<1$.  We define $H_t=f^t(H_0)$ for all integers $t\ge 0$ and
reparametrize the interval $[0,1)$ in the product $W\times [0,1]$ such
that $\bdry H_t=W\times t=W_t$, and the parameter $t$ now takes values
in $[0, \infty)$.  Finally, for any $t\ge 0$ we define $(H_t,V)$ to be
the compression body cut from $(H,V)$ by $W\times t=W_t$.

Since any automorphism $f\from (H,V)\to (H,V)$ of a compression body,
after  a suitable isotopy, agrees within a collar $W \times I$ of
$W\subseteq \bdry H$  with a product homeomorphism, a further
``vertical" isotopy of  $f$  within this collar gives:

\begin{lemma}
Every automorphism of a compression body is isotopic to an outward
expanding automorphism.
\end{lemma}

Henceforth, we shall always assume that automorphisms $f$ of
compression bodies have been isotoped such that they are outward
expanding.

Let ${\cal E}=\{E_i,i=1\ldots q\}$. be a collection of discs essential
in $(H_0,V)$, with $\bdry E_i\subseteq W_0$, cutting $H_0$ into a
product of the form $V\times I$, possibly together with one or more
balls. Such a collection of discs is called a {\it complete}
collection of discs.  (When $V=\emptyset$, $\cal E$ cuts $H$ into one
or more balls.)    We abuse notation by also using $\cal E$ to denote
the union of the discs in $\cal E$.  We further abuse notation by
often regarding $\cal E$ as a collection of discs properly embedded in
$H$ rather than in $H_0$, using the obvious identification of $H_0$
with $H$.  Thus, for example, we shall speak of $W$-parallel discs in
$\cal E$, meaning discs isotopic to discs in $W_0$. Corresponding to
the choice of a complete $\cal E$, there is a dual object $\Gamma$
consisting of the surface $V$ with a graph attached to $\intr (V)$ at
finitely many points, which are vertices of the graph.  The edges
correspond to discs of $\cal E$; the vertices, except those on $V$,
correspond to the complementary balls; and the surface $V$ corresponds
to complementary product.  If $V=\emptyset$, and $H$ is a handlebody,
then $\Gamma$ is a graph.

Now $H_0$ can be regarded as a  regular neighborhood of $\Gamma$, when
$\Gamma$ is embedded in $H_1$ naturally, with $V=V\times0$.  By
isotopy of $f$ we can arrange that $f({\cal E})$ is transverse to the
edges of $\Gamma$ and meets $H_0=N(\Gamma)$ in discs, each isotopic to
a disc of $\cal E$.   Any collection $\cal E$ of discs properly
embedded in $H_0$ (or $H$) with $\bdry{\cal E}\subseteq W_0$ (or
$\bdry{\cal E}\subseteq W$), not necessarily complete and not
necessarily containing only essential discs, is called {\it
admissible} if every component of $f({\cal E})\cap H_0$ is a disc
isotopic to a disc of $\cal E$.  We have shown:

\begin{lemma}\label{AdmissibleExistsLemma}
After a suitable isotopy every outward expanding automorphism $f\from
(H, V) \to (H, V)$  admits a complete admissible collection  ${\cal E}
\subseteq H_0$ as above, where every  $E_i\in \cal E$  is a
compressing disc of $W$.
\end{lemma}

We shall refer to admissible collections $\cal E$ of discs $(E_i,
\bdry E_i)\embed (H,W)$ as {\it systems} of discs.  Sometimes, we
shall retain the adjective ``admissible" for emphasis, speaking of
``admissible systems."  A system may contain discs which are not
compressing discs of $W$. Also, even if a system contains these
$W$-parallel discs, the definition of completeness of the system
remains the same.  We will say a system is {\it $W$-\-parallel} if
every disc in the system is $W$-parallel.

Let $P_j$ denote the holed disc $f(E_j)-H_0$.  Let $m_{ij}$ denote the
number of parallel copies of $E_i$ in $f(E_j)$, and let $M = M({\cal
E})$ denote the matrix $(m_{ij})$, which will be called the {\it
incidence matrix} for $\cal E$ with respect to $f$.  A system $\cal E$
is {\it irreducible} if the incidence matrix is irreducible.  In terms
of the discs $E_i$, the system $\cal E$ is irreducible if for each
$i,j$ there exists a $k\ge 1$ with $f^k(E_j)\cap H_0$ containing at
least one disc isotopic to $E_i$. It is a standard fact that a matrix
$M$ with non-negative integer entries has an eigenvector $x$ with
non-negative entries and that the corresponding eigenvalue
$\lambda=\lambda({\cal E})$ satisfies $\lambda\ge 1$. If the matrix is
irreducible, the eigenvector is unique and its entries are positive.

It turns out that the lack of reducing surfaces for $f$ in $(H,V)$ is
related to the existence of irreducible complete systems. The
following are proven in \cite{UO:Autos}.

\begin{lemma}\label{SystemsReducingLemma}
Suppose $f\from (H,V)\to (H,V)$ is an automorphism of a compression
body, and suppose that there is an (admissible) system which is not
complete and not $W$-parallel.  Then there is a reducing surface for
$f$.
\end{lemma}

\begin{proposition}\label{IrreducibleReducingProp}
If the automorphism $f\from (H,V)\to (H,V)$ is generic, then there is
a complete irreducible system $\cal E$ for $f$.  Also, any
non-$W$-parallel complete system $\cal E$ has a complete irreducible
subsystem ${\cal E}'$ with no $W$-parallel discs.  Further,
$\lambda({\cal E}')\le \lambda(\cal E)$, and $\lambda({\cal E}')<
\lambda({\cal E})$ if $\cal E$ contains $W$-parallel discs.
\end{proposition}

We always assume, henceforth, that $f$ is generic, so there are no
reducing surfaces. Assuming $\cal E$ is {\it any} (complete)
irreducible system of discs for $f$ in the compression body $(H,V)$,
we shall construct a branched surface $B=B({\cal E})$ in $\intr
(H)$. First we construct $B_1=B\cap H_1$:  it is obtained from
$f({\cal E})$ by identifying all isotopic discs of $f({\cal E})\cap
H_0$.  To complete the construction of $B$ we note that $f(B_1-\intr
(H_0))$ can be attached to $\bdry B_1$ to obtain $B_2$ and
inductively, $f^i(B_1-\intr (H_0))$ can be attached to $B_i$ to
construct a branched surface $B_{i+1}$. Alternatively, $B_i$ is
obtained by identifying all isotopic discs of $f^i({\cal E})\cap H_j$
successively for $j=i-1,\ldots,0$.  Up to isotopy, $B_i\cap H_r=B_r,\
r<i$, so we define $B=\cup_iB_i$.  The branched surface $B$ is a
non-compact branched surface with infinitely many sectors.  (Sectors
are completions of components of the complement of the branch locus.)
Note that the branched surface $B$ does not have boundary on $W_i$;
this follows from the irreducibility of $\cal E$.  If $\cal E$ is
merely admissible, the same construction works, but $B$ may have
boundary on $W_i,\  i\ge 0$.

If $x$ is an eigenvector corresponding to the irreducible system $\cal
E$, each component $x_i$ of the eigenvector $x$ can be regarded as a
weight on the disc $E_i\in \cal E$. The eigenvector $x$ now yields an
infinite weight function $w$  which assigns a positive weight to each
sector:

$$w(E_i)=x_i\text{ and }w(f^t(P_i))=x_i/\lambda^{t+1}\text{ for } t\ge
1.$$

Recall that $P_i$ is a planar surface, $P_i=f(E_i)-\intr (H_0)$. As is
well known, a weight vector defines a measured lamination $(\Lambda,
\mu)$ carried by $B$ if the entries of the weight vector  satisfy the
{\it switch conditions}.  This means that the weight on the sector in
$H_t-H_{t-1}$ adjacent to a branch circle in $W_t$ equals the sum of
weights on sectors in $H_{t+1}-H_t$ adjacent to the same branch
circle, with appropriate multiplicity if a sector abuts the branch
circle more than once. It is not difficult to check that our weight
function satisfies the switch conditions, using the fact that it is
obtained from an eigenvector for $M({\cal E})$.

At this point, we have constructed a measured lamination
$(\Lambda,\mu)$ which is $f$-invariant up to isotopy and fully carried
by the branched surface $B$.  Applying $f$, by construction we have
$f(\Lambda,\mu)=(\Lambda,\lambda\mu)$.  We summarize the results of
our construction in the following statement, which emphasizes the fact
that the lamination depends on the choice of non-$W$-parallel system:

\begin{proposition}\label{InvtLamProp2}
Suppose $f:(H,V)\to (H,V)$ is a generic automorphism of a
$3$\-dimensional compression body.  Given an irreducible system $\cal
E$ which is not $W$-parallel, there exists a lamination
$(\Lambda,\mu)$, carried by $B({\cal E})$, which is uniquely
determined, up to isotopy of $\Lambda$ and up to scalar multiplication
of $\mu$.

The lamination $\Lambda$ ``fills" $H_0$, in the sense that each
component of the complement is either contractible or deformation
retracts to $V$.  Also, $\Lambda\cup W$ is closed.
\end{proposition}

It will be important to choose the best possible systems $\cal E$ to
construct our laminations.  With further work, it is possible to
construct an invariant lamination with good properties as described in
the following proposition. This is certainly not the last word on
constructing ``good" invariant laminations, see \cite{LNC:Generic}.

\begin{proposition}\label{GoodLaminationProp}
Suppose $f\from (H,V) \to (H,V) $ is a generic automorphism of a
compression body. Then there is a system $\cal E$ such that there is
an associated 2-dimensional measured lamination $\Lambda\subseteq
\intr (H)$ as follows:   It has a transverse measure $\mu$ such that,
up to isotopy, $f((\Lambda,\mu))=(\Lambda,\lambda \mu)$ for some
$\lambda> 1$. Further, the lamination has the following properties:

1) Each leaf $\ell$ of $\Lambda$ is an open $2$-dimensional disc.

2) The lamination $\Lambda$ fills $H_0$, in the sense that each
component of $H_0-\Lambda$ is contractible or deformation retracts to
$V$

3) For each leaf $\ell$ of $\Lambda$,  $\ell-\intr (H_0)$ is
incompressible in $H-(\intr (H_0)\cup \intr(V))$.

4) $\Lambda\cup (\bdry H-\intr (V))$ is closed in $H$.
\end{proposition}

Given a ``good" system $\cal E$, e.g. as in Proposition
\ref{GoodLaminationProp}, we now construct a dual 1-dimensional
invariant lamination.  This was done in \cite{UO:Autos} for
automorphisms $f:H\to H$ where $H$ is a handlebody.  We now construct
an invariant lamination in the case of an automorphism $f\from
(H,V)\to (H,V)$ where $(H,V)$ is a compression body.

\begin{thm}\label{CompressionBodyClassificationTheorem}
Suppose $f\from (H,V)\to (H,V)$ is a generic automorphism of a
compression body with invariant 2-dimensional measured lamination
$(\Lambda,\mu)$ satisfying the incompressibility property. There is a
1-dimensional abstract measured non-compact 1-dimensional lamination
$(\Omega,\nu)$, and a map $\omega\from \Omega\to H_0-V$ such that
$f(\omega(\Omega,\nu))=\omega(\Omega,\nu/\lambda)$.  The map $\omega$
is an embedding on $\omega\inverse\big(N(\Lambda)\big)$ for some
neighborhood $N(\Lambda)$. The statement that
$f(\omega(\Omega,\nu))=\omega(\Omega,\nu/\lambda)$ should be
interpreted to mean that there is an isomorphism $h\from
(\Omega,\nu)\to (\Omega,\nu/\lambda)$ such that
$f\circ\omega=\omega\circ h$.
\end{thm}

\begin{proof}.
As in \cite{UO:Autos} we define $H_n$, $\mathcal{E}_n$ (for
$n\in\mathbb{Z}$) to be $f^n(H_0)$ and $f^n({\cal E})$
respectively. We choose a two-dimensional invariant lamination
$(\Lambda, \mu)$ with the incompressibility property. It can be shown,
as in \cite{UO:Autos}, that $\Lambda$ can be put in Morse position
without centers with respect to a height function $t$ whose levels
$t\in \integers$ correspond to $\bdry_e H_t$. It can also be shown as
in \cite{UO:Autos} that for $n>t$ $f^n({\cal E})\cap H_t$ consists of
discs belonging to a system $\mathcal{E}_t$ in $H_t$.  This implies
that for all real $t$, $(H_t,V)$ is a compression body homeomorphic to
$(H,V)$, and we can take a quotient of $(H_t,V)$ which gives
$\Gamma_t$, a graph attached to $V$, the quotient map being $q_t\from
H_t\to \Gamma_t$.  More formally, let $H_t^\Lambda$ denote the
completion of the complement of $\Lambda$ in $H_t$.  Points of the
graph in $\Gamma_t$ correspond to leaves of the 2-dimensional
lamination $\Lambda_t$ in $(H_t,V)$, or they correspond to balls in
$H_t^\Lambda$.  Points of $V\subseteq \Gamma_t$ correspond to interval
fibers in the product structure for the product component of
$H_t^\Lambda$.  We also clearly have maps $\omega_{st}\from
\Gamma_s\to \Gamma_t$ such that if $i_{st}:H_s\to H_t$ is the
inclusion, then $q_t\circ i_{st}=\omega_{st}\circ q_s$.  Then the
inverse limit of the maps $\omega_{st}$, for integer values of $t=s+1$
say, yields a 1-dimensional lamination limiting on $V$.  The
1-dimensional lamination is the desired lamination.  More details of
the ideas in this proof can be found in the proof in \cite{UO:Autos}
of the special case when $V=\emptyset$.
\end{proof}

We next investigate the 1-dimensional invariant lamination $\Omega$
further.  In \cite{UO:Autos} it is pointed out that (in the case of an
automorphism of a handlebody) $\Omega$ has, in general, much
self-linking.  An issue which was not addressed is whether the leaves
of $\Omega$ might be ``knotted."  This sounds unlikely, and in fact
knotting does not occur, but it is not even clear what it means for a
leaf to be knotted.  To clarify this notion, we first make some
further definitions.

\begin{defn}
A {\it spotted ball} is a pair $(K,\mathcal{D})$ where $K$ is a ball
and $\mathcal{D}$ is a collection of discs (or {\em spots}) in
$\partial K$.  A {\it spotted product} is a triple $(K,V,\mathcal{D})$
of the form $K=V\times I$ with $V=V\times 0$ and with $\mathcal{D}$ a
collection of discs in the interior of $V\times 1$. The spotted
product can also be regarded as a pair $(K,R)$, where $R=\bdry
K-\intr(V)-\intr(\mathcal{D})$. If $(K,\mathcal{D})$ is a spotted ball
we let the corresponding $V=\emptyset$.  If $(K,V,\mathcal{D})$ is a
spotted ball or connected product, (where $V=\emptyset$ if $K$ is a
ball) we say that it is a {\em spotted component}.
\end{defn}

Spotted balls occur naturally in our setting as follows:  Suppose $H$
is a handlebody.  Define $\widehat{H}_1$ to be $H_1$ cut open along
$\mathcal{E}_1$, or $H_1|\mathcal{E}_1$. We let
$\Gamma_t\cap\widehat{H}_1$ (or $H_t\cap\widehat{H}_1$) denote
$\Gamma_t|\mathcal{E}_1$ (or $H_t|\mathcal{E}_1$). The disc system
$\mathcal{E}_1$ determines a collection of spots
$\widehat{\mathcal{E}}_1\subseteq\partial\widehat{H}_1$ ($2$ spots for
each disc of $\mathcal{E}_1$). If $K$ is a component of
$\widehat{H}_1$ and $\mathcal{D}=K\cap\widehat{\mathcal{E}}_1$ then it
is clear that $(K,\mathcal{D})$ is a spotted ball.  If we make similar
definitions for a compression body $(H,V)$, then we also obtain
spotted product components $K$ in $\widehat{H}_1$.

In general, for $t\leq 1$, the triple $(K_t,V_t,\mathcal{D}_t)$ will
represent a spotted component where $K_t$ is a component of
$H_t\cap\widehat{H}_1$ and
$\mathcal{D}_t=K_t\cap\widehat{\mathcal{E}}_1$. Usually a spotted ball
$(K_t,\mathcal{D}_t)$ is a ball-with-two-spots, corresponding to arc
components $L_t$ of $\Gamma_t\cap\widehat{H}_1$. The other spotted
balls have more spots, $K_t$ corresponding to star components $L_t$ of
$\Gamma_t\cap\widehat{H}_1$. We denote $L_t\cap\partial K_t$ by
$\partial L_t$. It is clear that the spots of $\mathcal{D}_t$
correspond to points of $\partial L_t$ or to the edges of $L_t$.  We
also note that $L_t\subseteq K_t$ is uniquely defined, up to isotopy
(rel $\mathcal{D}_t$) by the property that it is a star with a vertex
in each spot and that it is ``unknotted''. We say $L_t$ is {\em
unknotted} in $K_t$ if it is isotopic in $K_t$ (rel $\partial K_t$) to
$\bdry K_t$.

When we are dealing with a compression body ($V_t\ne \emptyset$), we
also have component(s) of $L_t$ consisting of a component $X$ of $V$
with some edges attached at a single point of $X$.  In fact, it is
often useful to consider the union of star graphs $\hat L_t$, which is
the closure of $L_t-V$ in $L_t$.    These components of $L_t$
correspond to connected spotted products $K_t$ and the edges
correspond to spots.  The natural equivalence classes for the $ L_t$'s
are isotopy classes (rel $\mathcal{D}_t$), and for the $\hat L_t$'s
isotopy classes (rel $\mathcal{D}_t\cup X$). Given $K_t$ corresponding
to a spotted product, there are many such equivalence classes, even if
restricted to ``unknotted'' $\hat L_t$'s.  Here our $\hat L_t$ is
unknotted if it can be isotoped rel $\mathcal{D}_t$ to $V\times 1$.

We will often use ``handle terminology."    $H_t$ is constructed from
the 0-handles $K_t$ (or spotted products), with 1-handles
corresponding to edges of $\Gamma_t$ attached to the spots.

\begin{defn}
If $\Omega$ is the 1-dimensional dual lamination associated to a
2-dimensional invariant lamination $\Lambda$ for a generic
automorphism $f\from (H,V)\to (H,V)$, which in turn is associated to a
complete system $\cal E$, then a leaf $\ell$ of $\Omega$ is {\it
unknotted} if every component of $\ell\cap\widehat H_1$, contained in
a component $K$, say of $\widehat H_1$, is $\bdry$-parallel in $K$, or
isotopic to an arc embedded in $\bdry K$, whether $K$ is a spotted
ball or a spotted product.  In a spotted product  $(K,V,\mathcal{D})$,
we require that the arc be parallel to $R=\bdry
K-\intr(V)-\intr(\mathcal{D})$.
\end{defn}

\begin{thm}\label{UnknottedTheorem}
Let $f\from (H,V)\to (H,V)$ be a generic automorphism of a compression
body.  The 1-dimensional dual lamination $\Omega$ associated to a
system $\mathcal{E}$ is unkotted in $(H,V)$.
\end{thm}

To prove the theorem, we prove some lemmas from which the theorem
follows almost immediately.  We first make a few more definitions.

\begin{figure}[ht]
\centering \scalebox{1.0}{\includegraphics{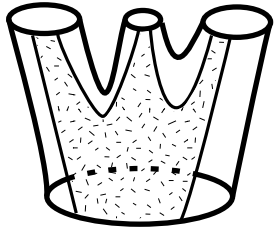}}
\caption{\small Polygonal disc.  The complement is $P_t$.}
\label{Aut3HalfSurface}
\end{figure}

If a spotted ball $K_t$ has $k$ spots then we can embed a polygonal
 disc $T_t$ with $2k$ sides in $R_t=\partial K_t -
 \intr({\mathcal{D}}_t)$, see Figure \ref{Aut3HalfSurface}.
 Alternate sides of $T_t$ are identified with embedded subarcs of each
 component of $\mathcal{D}_t$ in $\partial K_t$. We let  $P_t=\bdry
 K_t -\intr(T_t)$.  Then the pair $(P_t,\mathcal{D}_t)$ is a {\em
 \half for $(K_t,\mathcal{D}_t)$}.  From another point of view, $P_t$
 is a disc embedded in $\bdry K_t$ containing $\mathcal{D}_t$ and with
 $\bdry P_t$ meeting the boundary of each disc in $\mathcal{D}_t$ in
 an arc.  The reason we call $(P_t,\mathcal{D}_t)$ a spinal pair is
 that there is a homotopy equivalence from $(P_t,\mathcal{D}_t)$ to
 $(L_t,\text{closed edges of }L_t)$, which first collapses $K_t$ to
 $P_t$, then replaces the discs of $\mathcal{D}_t$ by edges, and
 $P_t-\mathcal{D}_t$ by a central vertex for the star graph $L_t$.
 Thus the \spinal is something like a spine.  There is of course a
 similar homotopy equivalence between $(K_t,\mathcal{D}_t)$ and $L_t$.

If $(K_t,V_t,\mathcal{D}_t)$ is a spotted product, we let  in
$P_t=\bdry K_t-\intr(V_t)$ (which is just $V_t\times 1$ when $K_t$ is
given the product structure $K_t=V_t\times I$) and the {\it \half} for
$(K_t,V_t,\mathcal{D}_t)$ is $(P_t, \mathcal{D}_t)$.  In the case of
spotted products, there is no choice for $P_t$.  Again, there is a an
obvious homotopy equivalence between the \half and $\Gamma_t$, taking
discs of $\mathcal{D}_t$ to closed edges.

Finally, it will be convenient later in the paper to have another
notion, used mostly in the case of a spotted ball
$(K_t,\mathcal{D}_t)$.  A {\it subspinal pair} for
$(K_t,V_t,\mathcal{D}_t)$ is a pair $(\breve P_t,\breve{\mathcal{D}}_t)$
embedded in a spinal pair $(P_t,\mathcal{D}_t)$, such
that $\breve{\mathcal{D}}_t$ is (the union of) a subcollection of the
discs in $\mathcal{D}_t$, and such that   $\breve P_t\embed P_t$ is an
embedded subdisc whose boundary meets each boundary of a disc in
$\breve{\mathcal{D}}_t$ in a single arc.  Exceptionally, we also allow
a subspinal pair consisting of a single disc component $D$ of
$\mathcal{D}_t$.  In this case the subspinal pair is $(D,D)$.

\begin{definition}  Fix $s\leq t$ and consider a spotted component $(K_{s},\mathcal{D}_s))$
contained in a spotted component $(K_t,\mathcal{D}_t))$. We say that
$(K_{s},\mathcal{D}_s))$ is {\em $\bdry$-parallel} to
$(K_t,\mathcal{D}_t))$ if there exist \halfs $P_s$ and $P_t$ such that
there is an isotopy of $(K_{s},\mathcal{D}_{s})$ in
$(K_t,\mathcal{D}_t)$ such that, after the isotopy,
$(P_{s},\mathcal{D}_{s})$ is contained in $(P_t,\mathcal{D}_{t})$.
Then we say the \spinal $(P_{s},\mathcal{D}_{s})$ is {\em parallel} to
the \spinal $(P_t,\mathcal{D}_{t})$.  Similarly, we say a {\it
subspinal pair}  $(\breve P_s,\breve{\mathcal{D}}_s)$ is {\em parallel}
to $(P_t,\mathcal{D}_{t})$ if there is an isotopy from the first pair
into the second.  Finally, we also say that the graph $(L_s, \bdry
L_s)$ is {\em parallel} to $(P_t,\mathcal{D}_{t})$ if it can be
isotoped in $(K_t,\mathcal{D}_{t})$ to  $(P_t,\mathcal{D}_{t})$.

We will need one more variation of the definition.  In Section
\ref{MixedBodies} we will have a spotted product $K_t$ with closed
interior boundary and contained in a handlebody $K_t\cup H'$, with
$H'$ a handlebody; the interior boundary of $K_t$ is  ``capped" with a
handlebody $H'$.  Then we speak of $(L_s, \bdry L_s)$,
$(P_{s},\mathcal{D}_{s})$, or $(\breve P_s,\breve{\mathcal{D}}_s)$
being {\it parallel in $K_t\cup H'$} to $(P_t,\mathcal{D}_{t})$.
\end{definition}

\begin{lemma} \label{GraphToSpinalLemma}
For $s<t\le 1$, suppose $K_t$ is a spotted component of $H_t\cap \widehat
H_1$; suppose $(P_t,\mathcal{D}_t)$ is a \spinal for $K_t$; and
suppose $K_s$ is a spotted component of $H_s\cap \widehat H_1$ with
$K_s\subseteq K_t$.   (If $K_t$ is a spotted product, it might be
capped by $H'$ on its interior boundary.)

Then:

i) If $K_s$ is a spotted ball, and the star graph $L_s$ is parallel to
$(P_t,\mathcal{D}_{t})$ (in $K_t\cup H'$) then there exists a \spinal
$(P_{s},\mathcal{D}_{s})$ parallel to $(P_t,\mathcal{D}_t)$ (in
$K_t\cup H'$).

ii) If $K_s$ is a spotted ball, and the embedded star subgraph $\breve
L_s$ of $L_s$  with ends in discs of $\breve{\mathcal{D}}_s\subseteq
\mathcal{D}_s)$ is parallel to $(P_t,\mathcal{D}_{t})$ in $K_t$ (in
$K_t\cup H'$) then there exists a subspinal pair $(\breve
P_s,\breve{\mathcal{D}}_s)$ parallel to $(P_t,\mathcal{D}_t)$ in $K_t$
(in $K_t\cup H'$).

iii) If $K_s$ is a spotted product, and the embedded star subgraph
$\breve L_s$ of $\hat L_s$  with ends in discs of
$\breve{\mathcal{D}}_s\subseteq \mathcal{D}_s)$ is parallel to
$(P_t,\mathcal{D}_{t})$ in $K_t$ (in $K_t\cup H'$) then there exists a
subspinal pair $(\breve P_s,\breve{\mathcal{D}}_s)$ parallel to
$(P_t,\mathcal{D}_t)$ in $K_t$ (in $K_t\cup H'$).  In particular, if
$\breve L_s=\hat L_s$, then we obtain a subspinal pair $(\hat
P_s,\mathcal{D}_s)$ including all spots of $K_s$.

iv) If $K_t$ is a spotted ball, a choice of isotopy of pairs of
$(L_t,\bdry L_t)$ in $(K_t,\mathcal{D}_t)$ to $\bdry K_t$ corresponds
to a choice of a \spinal $P_t$.

v)  If $K_t$ is a connected spotted product, a choice of isotopy of
pairs of $(\hat L_t,\bdry \hat L_t)$ in $(K_t,\mathcal{D}_t)$ to
$(P_t,\mathcal{D}_t)$ corresponds to a choice of a subspinal pair
$(\hat P_t,\mathcal{D}_t)$, including all spots, for $K_t$.
\end{lemma}

\begin{proof}
i) The argument is the same regardless $K_t$ is capped with $H'$ or
not.  Isotope $L_s$ to $L\subseteq P_t$ (possibly passing through
$H'$).  Now note that a regular neighborhood $N(L)$ of $L$ in $K_t$ is
a spotted ball isotopic to $K_{s}$, the component of $H_s\cap
\widehat{H}_t$ corresponding to (and containing) $L_{s}$.  We let
$P=N(L)\cap \bdry K_t$, then $(P,P\cap \mathcal{D}_t)$ is a spinal
pair for the  spotted ball $(N(L),N(L)\cap \mathcal{D}_t)$.  We can
use the inverse of the isotopy taking $L_{s}$ to $L$ to obtain an
isotopy taking $N(L)$ to $K_{s}$. This isotopy takes $(P,P\cap
\mathcal{D}_t)$ to a \spinal $(P_s,\mathcal{D}_s)$ in $K_{s}$. By
construction, $(K_{s},\mathcal{D}_{s})$ is parallel to
$(K_t,\mathcal{D}_{t})$ as desired and we have proved the isotopy of
$(L_s,\bdry L_s)$ to $(P_t,\mathcal{D}_t)$ determines a $P_s$ such
that $(P_s,\mathcal{D}_s)$ is isotopic to  $(P_t, \mathcal{D}_t)$.

ii), iii), The proofs are very similar to the proof of i).

iv) If we consider $L_t$ as a core graph of $K_t$, an isotopy of
$(L_t,\bdry L_t)$ in $(K_t,\mathcal{D}_t)$ to $\bdry K_t$  gives a
graph $(L,\bdry L)\subseteq (\bdry K_t,\mathcal{D}_t)$.  Again we let
$P=N(L)\cap \bdry K_t$ containing $\mathcal{D}_t$ to form a \spinal
$(P,\mathcal{D}_t)$ for the spotted ball $(N(L),\mathcal{D}_t)$ .
$N(L)$ is the same as $K_t$ up to isotopy, and $N(L)\cap \bdry K_t$
gives a \spinal for $N(L)$, so we have a \spinal for $K_t$.  It is
easy to recover the $\hat L_t$ and $L_t$ from the spinal pair.

v) The proof is similar to the proof of iv).
\end{proof}

The following lemma analyzes the ``events" which change the $L_t$'s
corresponding to spotted balls or products $K_t$.

\begin{figure}[ht]
\centering
 \psfrag{L}{\fontsize{\figurefontsize}{12}$L$} \psfrag{H}{\fontsize{\figurefontsize}{12}$H$} \psfrag{t}{\fontsize{\figuresmallfontsize}{12}$t$} \psfrag{s}{\fontsize{\figuresmallfontsize}{12}$s$}

\scalebox{1.0}{\includegraphics{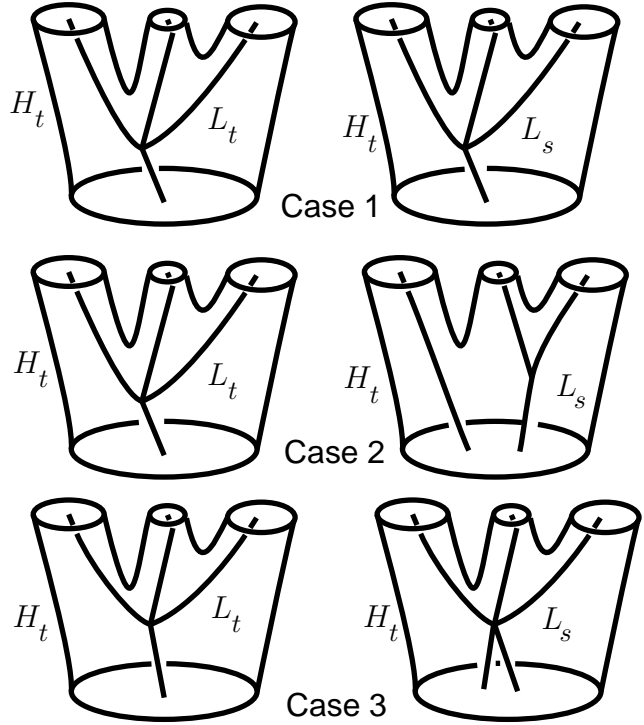}} \caption{\small
Event in $H_t$.} \label{Aut3Splitting}
\end{figure}

\begin{lemma} \label{EventLemma}
Recall $K_t$ is a component of $H_t\cap \widehat H_1$ and $L_t$ is the corresponding
component of $\Gamma_t\cap \widehat H_1$.  Suppose $s$ and $t$ are
consecutive regular values, in the sense that there is exactly one
critical value between them. Then one of the following holds (see
Figure \ref{Aut3Splitting}):

\begin{enumerate}

\item If the corresponding critical point is not contained in $K_t$ then
$L_t=\Gamma_t\cap K_t$ is the same as $L_s=\Gamma_s\cap  K_t$, otherwise

\item $L_t$ is obtained from two components $L_s$ and
$L_s'$ by isotoping an edge of  $L_s$ to an edge of  $L_s'$  (in the
complement of the remainder of  $\Gamma_s\cap  K_t$) and identifying
the two edges, or

\item $L_t$ is obtained from a component $L_s$ of
$\Gamma_s\cap \widehat H_1$ by isotoping an edge of  $L_s$ to an other
edge of  $L_s$  (in the complement of the remainder of  $\Gamma_s\cap
K_1$) and identifying the two edges.
\end{enumerate}

Figure \ref{Aut3Splitting} shows the events in case $K_t$ is a spotted
ball.
\end{lemma}

\begin{proof}   Recall that in the context of Theorem \ref{UnknottedTheorem} we
are assuming that the discs $\mathcal{E}_1$ are in Morse position
without centers with respect to the height function in
$H_1-\intr({H}_0)$.   The lemma follows from the fact that there is
exactly one critical value between $s$ and $t$  for the Morse height
function on the system of discs in $H_1$.
\end{proof}

We classify critical values according to the three cases of Lemma
\ref{EventLemma}.

In the following statement ``consecutive" regular values are again
regular values separated by exactly one critical value.

\begin{lemma} \label{SpinalPairLemma}   For $s<t\le 1$ consecutive regular values, suppose $K_t$ is a spotted component
of $H_t\cap \widehat H_1$ and suppose $K_s$ is a spotted component of
$H_s\cap \widehat H_1$ with $K_s\subseteq K_t$.  Then for any choice of
\spinal $(P_{t},\mathcal{D}_{t})$ for $(K_t,\mathcal{D}_t)$ there is a
choice of  \spinal $(P_{s},\mathcal{D}_{s})$ for $(K_s,\mathcal{D}_s)$
such that  the  \spinal $(P_{s},\mathcal{D}_{s})$ is parallel to the
\spinal $(P_{t},\mathcal{D}_{t})$, and hence $(K_s,\mathcal{D}_s)$ is
$\bdry$-parallel in $(K_t,\mathcal{D}_t)$.
\end{lemma}

\begin{proof}

We first consider the case that  $K_t$ is a spotted ball, and we
choose a \spinal $(P_t, \mathcal{D}_t)$.

We are assuming that the discs $\mathcal{E}_1$ are in Morse position
without centers with respect to the height function in
$H_1-\intr({H}_0)$.  We assume $s$ and $t$ are ``consecutive''
regular values, in the sense that there is exactly one critical value
between them, and that the corresponding critical point lies in $K_t$.
Lemma \ref{EventLemma} classifies the events corresponding to the
critical value: A component $L_{s}$ of $\Gamma_s\cap K_t$ is obtained
from a component $L_t$ of $\Gamma_t\cap K_t$ as described in the
lemma.  In Case 1, the component $L_{s}$ coincides with the
corresponding $L_t$. In particular $L_{s}\subseteq L_t$ in this case.
In Case 2, where $L_{s}$ is isotopic to a subspace of $L_t$, which is
split, we can isotope it so that $L_{s}\subseteq L_t$.  In Case 3,
$L_{s}$ is obtained from $L_t$ by splitting the edge $e$, and we can
isotope $L_{s}$ so that it is contained in $L_t$ away from a
neighborhood $N(e)\subseteq K_t$ of $e$.

But $L_t$ is clearly {\em parallel to $P_t$}, in the sense that there
exists an isotopy of $(L_t,\partial L_t)$ in $(K_t,\mathcal{D}_t)$
taking $(L_t,\bdry L_t)$ to $(P_t,\mathcal{D}_t)$.   Therefore, in
Cases 1 and 2, when $L_{s}\subseteq L_t$, the parallelism of $L_t$
realizes the parallelism of $L_{s}$ and by Lemma
\ref{GraphToSpinalLemma}, we obtain a \spinal
$(P_{s},\mathcal{D}_{s})$ which is parallel to the \spinal
$(P_{t},\mathcal{D}_{t})$.

We deal with Case 3 as follows.  We can regard the isotopy of $L_t$
into $P_t$ as a composition of two isotopies: first isotope a
neighborhood $L_t\subseteq K_t$ (containing $e)$) into a neighborhood
$N(P_t)\simeq P_t\times I$. We may assume that $L_t$ is taken to a
horizontal slice $P_t\times\{x\}$. The isotopy of $L_t$ is completed
by vertical projection. Now we can assume $L_{s}\subseteq K_t$ is
contained in $L_t$ away from $N(e)$. After the first isotopy we then
have that $\big(L_{s}-N(e)\big)\subseteq P_t\times\{x\}$. But we can
isotope $L_{s}\cap N(e)$ (rel $\partial N(e)$) so that $L_{s}\subseteq
P_t\times\{x\}$. Vertical projection gives $L_{s}\subseteq P_t$.
Again, by Lemma \ref{GraphToSpinalLemma} we obtain a \spinal
$(P_{s},\mathcal{D}_{s})$ which is parallel to the \spinal
$(P_{t},\mathcal{D}_{t})$.

Now we consider the case that $K_t$ is a spotted product and $K_s$ is
a spotted ball.  There is no choice for the \spinal $(P_t,
\mathcal{D}_t)$. Again, analysis of the cases of Lemma
\ref{EventLemma} shows that in each case $L_s$ can be isotoped to
$P_t$, so by Lemma \ref{GraphToSpinalLemma} we again get a \spinal
$(P_{s},\mathcal{D}_{s})$ which is parallel to $(P_t, \mathcal{D}_t)$.

Finally, if both $K_t$ and $K_s$ are spotted products, then there is
no choice for either spinal pair $(P_t, \mathcal{D}_t)$ or $(P_s,
\mathcal{D}_s)$.  The surface $P_s$ is isotopic to $P_t$ and
$\mathcal{D}_s$ represents a subcollection of the discs of
$\mathcal{D}_t$.  Hence $(P_s, \mathcal{D}_s)$ is automatically
parallel to $(P_t, \mathcal{D}_t)$
\end{proof}

\begin{proposition}\label{ParallelLemma}  For all $s<t\le 1$, for any component $K_t$ of $H_t\cap \widehat H_1$ and for any
component $K_s$ of $H_s\cap \widehat H_1$ with $K_s\subseteq K_t$, $K_s$
is $\bdry$-parallel in $K_t$.
\end{proposition}

\begin{proof}  There is a sequence $t_i$, $1\le i\le n$ where
$t_1=s$, $t_n=t$ and with consecutive values $t_i$ separated by a
single critical value.  By Lemma \ref{SpinalPairLemma}, given a
\spinal $(P_{t_n}, \mathcal{D}_{t_n})$ for $K_{t_n}$, there is a
parallel \spinal $(P_{t_{n-1}}, \mathcal{D}_{t_{n-1}})$ for
$K_{t_{n-1}}$.  We use induction with decreasing $i$:  If we have
found a spinal pair $(P_{t_{i}}, \mathcal{D}_{t_{i}})$ for $K_{t_{i}}$
for $K_{t_i}$ which is parallel to $(P_{t_{i+1}},
\mathcal{D}_{t_{i+1}})$ for $K_{t_{i+1}}$ for $K_{t_{i+1}}$, then we
can find a \spinal $(P_{t_{i-1}}, \mathcal{D}_{t_{i-1}})$ for
$K_{t_{i-1}}$ which is parallel to the pair $(P_{t_{i}},
\mathcal{D}_{t_{i}})$ for $K_{t_{i}}$ for $K_{t_i}$.  This then proves
we have a sequence of parallel \spinals.  It is easy to see that
parallelism for \spinals is transitive, whence we obtain the statement
of the lemma.
\end{proof}

\begin{proof}[Proof of Theorem \ref{UnknottedTheorem}]  For a leaf $\ell$ of $\Omega$, every component of $\ell
\intersect \widehat H_1$ is the core of a spotted ball $K_s$ with two
spots, for some $s$.  $K_s$ is contained in some component $K$ of
$\widehat H_1$.  Letting $t=1$ in the previous lemma, we see that $K_s$ is
$\bdry$-parallel in $K_t$, which implies the arc of $\ell$ in $\widehat
H_1$ is $\bdry$-parallel, hence unknotted.
\end{proof}

\section{Automorphisms of sphere bodies, the Train Track Theorem}\label{TrainTrackSection}

Our goal in this section is to prove Theorem \ref{TrainTrackTheorem}.

We recall from the Introduction that $g$ is {\em train track generic}
if it is generic on $H_1$ and the corresponding invariant
2-dimensional lamination determines, through a quotient map $H\to G$
on a graph $G$, a homotopy equivalence which is a train track map.

Our strategy for proving Theorem  \ref{TrainTrackTheorem} is similar
to that of \cite{BH:Tracks} for improving homotopy equivalences of
graphs. Given an automorphism $f\from M\to M$ we first isotope it so
that it preserves the canonical splitting $M=H\cup H'$ (see Section
\ref{Splitting} and Theorem \ref{CanonicalHeegaardTheorem}). Let
$g=f|_{H}\colon H\to H$. By Theorem \ref{HandlebodyClassificationThm}
$g$ is either periodic, reducible or generic. If it is periodic or
reducible the proof is done, so assume that it is generic.  We
consider the corresponding 2-dimensional lamination and homotopy
equivalence of a graph. If it is a train-track generic map, the proof
is also done, so assume otherwise. This means that a power of the
induced homotopy equivalence on the graph admits ``back tracking'', in
the sense of \cite{BH:Tracks}. In that context, there is a sequence of
moves that simplify the homotopy equivalence.  Similarly in our
setting there are analogous moves that realize this
simplification. The result of this simplification is an automorphism
$g'$ of $H$ which is either periodic, reducible or generic. In the
last case, the growth rate of $g'$ is smaller than that of
$g$. Proceeding inductively, this process yields the proof of the
theorem. In general $g$ and $g'$, which are automorphisms of $H$, will
not be isotopic. The key here is that they are both restrictions of
automorphisms $f$, $f'$ of $M$ to the handlebody $H$. These $f$ and
$f'$ will be isotopic (through an isotopy that does not leave
invariant the splitting, see e.g.  Example \ref{NonUniqueExample}).

In the light of the sketch above, it is clear that Theorem
\ref{TrainTrackTheorem} will follow from the result below.

\begin{thm}\label{BackTrackTheorem}
Let $f\colon M\to M$ be an automorphism preserving the canonical
splitting $M=H\cup H'$. Assume that $g=f|_{H}$ is generic and not
train track. Then $f$ is isotopic to an $f'$ preserving the splitting
such that $g'=f'|_{H}$ is either periodic, reducible or generic with
growth rate smaller than $g$.
\end{thm}

There are two technical ingredients in proving Theorem
\ref{BackTrackTheorem}. The first is to prove the existence of
``half-discs'', which are geometric objects that represent some of the
homotopy operations performed in the algorithm  of \cite{BH:Tracks}
(these operations are ``folding'' and ``pulling-tight''). In special
situations, when the half-discs are ``unholed,'' we can use them to
realize these homotopies by isotopies. These isotopies are described
in \cite{UO:Autos} (a ``down-moves'' corresponds to a ``fold'' and a
``diversion'' to a ``pulling-tight''). We can then proceed to the next
step of the algorithm remaining in the same isotopy class, simplifying
the automorphism by either reducing the growth rate or turning it into
a periodic or reducible automorphism.

The second ingredient is an operation called ``unlinking''; the aim is
to isotope $f$ so that we can obtain these special unholed half-discs
from general ones. As noted above, this makes it possible to proceed
with the algorithm. An unlinking preserves the canonical splitting
but, unlike the other moves, it is an isotopy which does not leave
invariant the canonical splitting. In the unlinking process we will
decompose the given half-disc into ``rectangles'' and a ``smaller''
half-disc.

We have sketched the proof of Theorem \ref{BackTrackTheorem} (which
implies Theorem \ref{TrainTrackTheorem}). We now turn to details. We
fix $f\colon M\to M$ preserving the splitting $M=H\cup H'$ and
consider its restriction $g\colon H\to H$. We assume that $g$ is
generic. Now consider the construction of the invariant laminations,
determining $H_t$, $\mathcal{E}_t$ and $\Gamma_t$ for $t\in
\mathbb{R}$, as in Section \ref{Laminations}. We also use the notions
of spotted balls as in Section \ref{Laminations}, and we work in the
space $\widehat{H}_1$, i.e., $H_1|\mathcal{E}_1$ as before.  In
general, for $t\leq 1$, the pair $(K_t,\mathcal{D}_t)$ will represent
a spotted ball where $K_t$ is a component of $H_t\cap\widehat{H}_1$
and $\mathcal{D}_t=K_t\cap\widehat{\mathcal{E}}_1$.

\begin{defns}\label{HalfDisc}
A {\em half-disc} is a triple $(\Delta,\alpha,\beta)$ where $\Delta$
is a disc and $\alpha$, $\beta\subseteq\partial\Delta$ are
complementary closed arcs. We may abuse notation and say that $\Delta$
is a half-disc.

Let $s<t$. A {\em half-disc for $H_s$ in $H_t$} is a half-disc
$(\Delta,\alpha,\beta)$ satisfying:
\begin{enumerate}[\upshape (i)]
\item $\Delta\subseteq K_t$ is embedded for some $K_t$,
\item $\beta$ is contained in a spot $D_t$ of $K_t$,
\item $\alpha\subseteq\partial K_s$ for some $K_s\subseteq K_t$ and
$\Delta\cap K_s=\alpha$,
\item $\alpha$ {\em connects} (i.e., intersects the boundaries of) two
distinct spots of $(K_s,\mathcal{D}_s)$.
\end{enumerate}

We refer to $\alpha$, the arc in $\partial K_s$ as the {\em upper
boundary of the half-disc for $H_s$ in $H_t$, $s<t$}, denoted by
$\partial_u\Delta$, and $\beta$ as its {\em lower boundary}, denoted
by $\partial_l\Delta$.

Similarly, a {\em rectangle} is a triple $(R,\gamma,\gamma')$ where
$R$ is a disc and $\gamma$, $\gamma'\subseteq\partial R$ are disjoint
closed arcs. We may abuse notation and say that $R$ is a rectangle.

Let $s<t$. A {\em rectangle for $H_s$ in $H_t$} is a rectangle
$(R_s,\gamma_s,\gamma_t)$ satisfying:

\begin{enumerate} [\upshape (a)]
    \item $R_s\subseteq K_t$ is embedded for some $K_t$,
    \item $\gamma_s\subseteq\partial K_s$ for some $K_s\subseteq K_t$
    and $R_s\cap K_s=\gamma_s$,
    \item $\gamma_s$ connects distinct spots of $(K_s,\mathcal{D}_s)$,
    \item $\gamma_t\subseteq\partial K_t-\intr({\mathcal{D}}_t)$,
    \item $\partial R_s-(\gamma_s\cup\gamma_t)\subseteq\mathcal{D}_t$.
\end{enumerate}
\end{defns}

We note that a half-disc $\Delta$ (for $H_s$) may be {\em holed}, in
the sense that $H_s$ may intersect the interior of $\Delta$. If this
intersection is essential then $\Delta$ also intersects $H_0$. We
define {\it holed rectangles} similarly.

\begin{defn}
Let $0\leq s<t\leq 1$ and $\Delta$ a half-disc for $H_s$ in $H_t$.
Suppose that $s=t_0<t_1<\cdots<t_k=t$ is a sequence of regular
values. We say that $\Delta$ is {\em standard with respect to
$\{\,t_i\,\}$} if there is $0\leq m<k$ such that the following
holds. The half-disc $\Delta$ is the union of rectangles $R_{t_i}$,
$s\leq t_i\leq t_{m-1}$, and a half-disc $\Delta_{t_m}$ (see
Figure~\ref{Aut3StandardDisc}) with the properties that:
\begin{enumerate}
    \item each $R_{t_i}$, $s\leq t_i< t_{m}$, is a rectangle for
    $H_{t_i}$ in $H_{t_{i+1}}$,
    \item the holes $\intr({R}_{t_i})\cap H_{t_i}$ are essential
    discs in $H_{t_i}$,
    \item $\Delta_{t_m}$ is an unholed half-disc for $H_{t_m}$ in
    $H_{t_{m+1}}$.
\end{enumerate}
\begin{figure}[ht]
\centering
\psfrag{R}{\fontsize{\figurefontsize}{12}${\hskip-1.7pt}\raise0.8pt\hbox{$R$}$}
\psfrag{t}{\fontsize{\figuresmallfontsize}{12}$t$}\psfrag{m}{\fontsize{\figuresmallerfontsize}{12}$m$}
\psfrag{i}{\fontsize{\figuresmallerfontsize}{12}${\hskip-0.7pt}i$}
\psfrag{g}{\fontsize{\figurefontsize}{12}$\gamma$}\psfrag{D}{\fontsize{\figurefontsize}{12}$\Delta$}

\scalebox{1.0}{\includegraphics{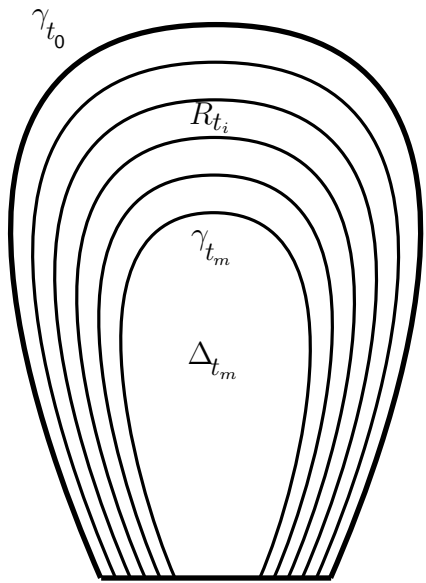}}
\caption{\small Standard half disc.} \label{Aut3StandardDisc}
\end{figure}
\end{defn}

We want to prove the following:

\begin{proposition}\label{P:half-disc}
Let $0=t_0<\dots<t_k=1$ be a complete sequence of regular values. If a
component of $H_0\cap\widehat{H}_1$ intersects a disc of
$\widehat{\mathcal{E}}_1$ in two distinct components $D_0$, $D_1$,
then there exists a half-disc $\Delta$ for $H_0$ connecting $D_0$ and
$D_1$ which is standard with respect to $\{\,t_i\,\}$.
\end{proposition}

To prove the proposition, we need:

\begin{lemma}\label{L:half-disc}
Suppose that there exist \halfs\ $P_{t_i}$ such that
$(K_{t_i},P_{t_i})$ is parallel to $(K_{t_{i+1}},P_{t_{i+1}})$. If a
component $K_0$ intersects a disc of $\mathcal{D}_1$ in two distinct
components $D_0$, $D_1$, then there exists a half-disc for $H_0$
connecting $D_0$ and $D_1$ which is standard with respect to
$\{\,t_i\,\}$.
\end{lemma}
\begin{proof}

Since $P_{t_0}=P_{0}$ is a \half\ it intersects both $D_0$, $D_1$.
Let $\gamma_{t_0}\subseteq P_{t_0}$ be an arc connecting $D_0$ to
$D_1$.

For any $i$, let $\gamma_{t_i}\subseteq P_{t_i}$ be an embedded arc
with endpoints in distinct spots of $K_{t_i}$. Consider an isotopy
realizing the parallelism between $(K_{t_i},P_{t_i})$ and
$(K_{t_{i+1}},P_{t_{i+1}})$, taking $P_{t_i}$ to a subsurface of
$P_{t_{i+1}}$. The arc $\gamma_{t_i}$ is taken to an arc
$\gamma_{t_{i+1}}\subseteq P_{t_{i+1}}$. It is clear that the isotopy
yields, via the Loop Theorem, a rectangle
$(R_{t_{i}},\gamma_{t_{i}},\gamma_{t_{i+1}})$ for $H_{t_{i}}$ in
$H_{t_{i+1}}$. Here we note that indeed $R_{t_i}\cap
K_{t_i}=\gamma_{t_i}$, as needed. It can happen, though, that
$H_{t_i}$ intersects $R_{t_i}$ in other components. In this case there
exists another $K_{t_i}'\neq K_{t_i}$ contained in
$K_{t_{i+1}}$. Consider then the graph $L_{t_i}'$ corresponding to
$K_{t_i}'$ and perturb $R_{t_i}$ so that $L_{t_i}'\cap R_{t_i}$ is
transverse. Now make $K_{t_i}'\cap R_{t_i}$ consist of essential discs
(in $H_{t_i}$).

From $\gamma_{t_0}$ defined in the beginning of the proof we use the
process described above to proceed inductively, obtaining rectangles
$(R_{t_{i}},\gamma_{t_{i}},\gamma_{t_{i+1}})$. We do that for as long
as $\gamma_{t_i}$ has endpoints in distinct spots of $K_{t_{i+1}}$
(and hence in distinct spots of $K_{t_i}$).

It is clear that eventually $\gamma_{t_{i}}$ will have both endpoints
in the same spot of the corresponding $K_{t_{i+1}}$. To see this
simply note that all $\gamma_{t_j}$'s are parallel (as pairs
$(\gamma_{t_j},\partial\gamma_{t_j})$ in $(K_1,\mathcal{D}_1)$) and
hence have endpoints in the same spot of $K_{1}$ by hypothesis.

So let $\gamma_{t_{m}}$ be the first such arc which has both endpoints
in the same spot of $K_{t_{m+1}}$, say $D$. In this case we still
build $\gamma_{t_{m+1}}\subseteq P_{t_{m+1}}$ through parallelism. Now
consider $\beta\subseteq (\partial D\cap\partial P_{t_{m+1}})$ with
$\partial\beta=\partial\gamma_{t_{m+1}}$. Then
$\beta\cup\gamma_{t_{m+1}}$ bounds a disc $D'\subseteq
P_{t_{m+1}}$. As before, the parallelism between $\gamma_{t_m}$ and
$\gamma_{t_{m+1}}$ defines a rectangle
$(R_{t_m},\gamma_{t_m},\gamma_{t_{m+1}})$. We consider
$\Delta_{t_m}=R_{t_m}\cup D'$ and push $(\Delta_{t_m},\beta)$ into
$(K_{t_{m+1}},D)$ as a pair. Now $(\Delta_{t_m},\gamma_{t_m},\beta)$
is a half-disc for $H_{t_m}$ in $H_{t_{m+1}}$.

From the construction it is clear that if $R_{t_m}$ is unholed then so
is $\Delta_{t_m}$. We claim that that is the case. Indeed, recall that
$\{\,t_i\,\}$ is complete and that $\gamma_{t_{m+1}}$ is the first
$\gamma_{t_i}$ to have both endpoints in the same spot of
$K_{t_{i+1}}$.  Therefore $K_{t_m}$ may be regarded as obtained from
$K_{t_{m+1}}$ by the event pictured Figure \ref{Aut3Splitting}, Case
3. In particular, there is no other $K'_{t_m}\neq K_{t_m} $ in
$K_{t_{m+1}}$. Therefore $R_{t_m}$ is unholed, so the same holds for
$\Delta_{t_m}$.

It is now easy to see that $\Delta=\big(\,\bigcup_{t_0\leq t_i\leq
t_{(m-1)}} R_{t_i}\,\big)\cup\Delta_{t_m}$ is 1) a half-disc for
$H_0$, 2) connects $D_0$ and $D_1$ and 3) is standard with respect to
$\{\,t_i\,\}$.
\end{proof}
\begin{remark}\label{UnholedHalf-Disc}
Note from the construction of $\Delta$ above that the rectangles
$R_{t_i}$ are holed only if there are two distinct $K_{t_i}$,
$K_{t_i}'$ contained in $K_{t_{i+1}}$, as pictured in Figure
\ref{Aut3Splitting}, Case 2.
\end{remark}

\begin{proof}[Proof of Proposition \ref{P:half-disc}]
Follows directly from the result above, Theorem \ref{UnknottedTheorem}
and Lemma \ref{ParallelLemma}.
\end{proof}

We wish to use a standard half disc to obtain a simplified
representative of the automorphism.

Let $\Delta$ be a half-disc for $H_0$, with upper boundary on some
$K_0$. If its interior did not intersect $H_0$, the half-disc $\Delta$
could be used to drag $H_0$ so it intersects ${\cal E}_1$ in fewer
discs, but recall that $\Delta$ may be holed. A solution would be to
isotope $1$-handles of $H_0$, to avoid $\Delta$. One might worry that
the handles that make holes in $\Delta$ may be ``linked'' with
$K_0$. More precisely, it may be impossible to isotope the handles of
$H_0$ in $H$ (leaving $K_0$ fixed) in order to avoid holes in
$\Delta$. We emphasize that the word ``linked'' is abused of here. The
property it describes is not intrinsic to the $K_0$'s in $K_1$, and
rather relative to the half-disc $\Delta$. In any case, these
``essential links'' do appear in the setting of handlebodies (see
\cite{UO:Autos,LNC:Generic, LNC:Tightness}). But if we allow the
isotopies to run through the whole manifold $M$, not just $H$, the
unlinking is always possible. The idea is to use isotopies of the kind
described in Example \ref{NonUniqueExample}.

\hop

More specifically we shall describe an isotopy of the type in Example
\ref {NonUniqueExample} and how it can be used to unlink.  We call an
isotopy of this type {\em an unlinking move} or just an {\em
unlinking}.

At this point we do not worry about the automorphism $f\colon M\to
M$. Although $f$ is the object of interest, we will intentionally
disregard it initially, thinking instead of the various $H_t\subseteq
H$, disc systems and sphere systems, as geometric objects embedded in
$M$.  We will bring $f$ back into the picture later in the argument.

As usual, let $K$ be a component of $\widehat{H}_1$. Let $0\leq s\leq
1$ and consider a component $K_s$ of $H_s\cap K$. Recall that a spot
$D_s$ of $(K_s,\mathcal{D}_s)$ corresponds to a disc $E_s\subseteq
H_s$. Such a disc extends to an essential sphere $S_s\embed M$ with
the following properties: 1) $S_s\cap H_s=E_s$ and 2) for $s<t\leq 1$,
$S_s\cap H_t$ consists of a single disc. We say that such a sphere is
{\em dual to $H_s$}, or just a {\em dual sphere} in  case the value of
$s$ is clear.

The unlinking move will be performed along a rectangle (coming from a
standard half-disc) and a sphere.  Let $K$ be a component of
$\widehat{H}_1$. Fix $0\leq s<t\leq 1$, with exactly one critical
point between $s$ and $t$ as usual. Consider a component $K_s$ of
$H_s\cap K$ and the component $K_t$ of $H_t\cap K$ containing
$K_s$. In this situation, most spots $K_t$ contain exactly one spot of
$K_s$, see Figure \ref{Aut3Splitting}. We shall take advantage of this
fact. Consider a rectangle $(R_s,\gamma_s,\gamma_t)$ for $H_s$ in
$H_t$ with $R_s\subseteq K_t$. Let $K_s$ be the component that
contains $\gamma_s\subseteq\partial R_s$. Assume that $\gamma_s$
connects distinct spots of $K_s$ (which is the case if the rectangle
comes from a standard half-disc). Finally, since most spots of $K_t$
contain exactly one spot of $K_s$, it is reasonable also to suppose
that at least one of the two sides of $R_s$ in a spot of $K_t$ lies in
a spot containing only one spot of $K_s$.  If this is the case, we can
assume that the dual sphere $S_s$ containing $D_s$ has the property
that $S_s\cap K_t$ is a spot of $K_t$. In other words, we are assuming
that $S_s$ is also dual to $H_t$, which is a strong restriction. Thus,
$S_s\cap R_s\subseteq\partial R_s$, a key property we need for the
unlinking move.

\begin{figure}[ht]
\centering
\psfrag{R}{\fontsize{\figurefontsize}{12}$R$}\psfrag{C}{\fontsize{\figurefontsize}{12}$C$} 
\psfrag{H  }{\fontsize{\figurefontsize}{12}$H$} \psfrag{S  }{\fontsize{\figurefontsize}{12}${\hskip2pt}S$}
\psfrag{t}{\fontsize{\figuresmallfontsize}{12}${\hskip2pt}t$}\psfrag{s}{\fontsize{\figuresmallfontsize}{12}$s$}
\psfrag{m}{\fontsize{\figuresmallerfontsize}{12}$m$}\psfrag{i}{\fontsize{\figuresmallerfontsize}{12}$i$}
\psfrag{g}{\fontsize{\figurefontsize}{12}$\gamma$}\psfrag{D}{\fontsize{\figurefontsize}{12}$D$}
\psfrag{N}{\fontsize{\figurefontsize}{12}${\hskip-3pt}N$}\psfrag{N'}{\fontsize{\figurefontsize}{12}$N'$}

\scalebox{1.0}{\includegraphics{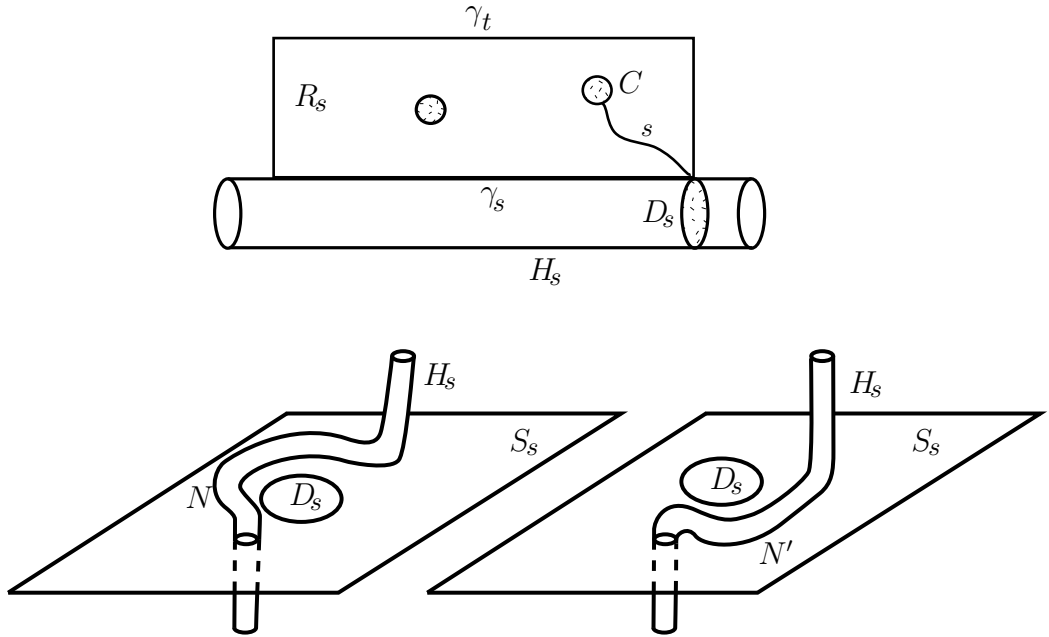}} \caption{\small
Unlinking move.} \label{Aut3Move}
\end{figure}

Consider a component $C$ of $H_s\cap R_s$ distinct from
$\gamma_s$. Assume that $C$ is an essential disc in $H_s$, which is
the case if the rectangle comes from a half-disc in standard
position. Let $K'_s$ be the component containing $C$. The move that we
will describe isotopes $K'_s$ along $R_s$ and $S_s$, see Figure
\ref{Aut3Move}. More precisely, consider an embedded path
$\sigma\subseteq R_s-\intr({H}_s)$ connecting $C$ to $D_s$. We first
isotope $C$ along $\sigma$ so it is close to $D_s\subseteq S_s$. This
move changes $K'_s$ in a neighborhood $N\subseteq K'_s$ of $C$. We
note that the condition imposed on $S_s$ that $S_s\cap K_t$ is a spot
of $K_t$ implies that $\sigma\cap S_s=\sigma\cap
D_s\subseteq\partial\sigma$. Hence the move described does not
introduce intersections of $K'_s$ with $S_s$.  So $S_s$ is preserved
as a dual sphere and the disc $\check D_s=S_s-\intr({D}_s)$ does not
intersect $H_s$ (except at $\partial\check D_s=\partial D_s$).

We now regard $N$ as the neighborhood of an arc $\alpha$ transverse to
$R_s$ and close to an arc $\alpha'\subseteq \partial D_s$ (this
$\alpha'$ also intersecting $R_s$). Let $\beta'\subseteq\partial D_S$
be the closed complement of $\alpha'$.

Now recall that the disc $\check{D}_s$ is disjoint from $H_s$. We can
push the arc $\alpha$ close to $\alpha'$ and isotope it along
$\check{D}_s$ to an arc $\beta\subseteq K_t-\intr({H}_s)$ close to
$\beta'$. More precisely, consider an arc $\beta\subseteq
\intr({K}_t)-{H}_s$ such that 1)
$\alpha\cap\beta=\partial\alpha=\partial\beta$ and that 2) the circle
$\alpha\cup\beta$ is parallel in $K_t-\intr({H}_s)$ to $\partial
D_s=\alpha'\cup\beta'$. Then $\alpha$ is isotopic to $\beta$ (along
the disc $\check{D}_s$). If $N'$ is a neighborhood of $\beta$ we can
isotope $N$ to $N'$. But $\beta\cap R_s=\emptyset$ then also $N'\cap
R_s=\emptyset$. Clearly, this operation does not change $H_s$ away
from $N$, therefore such an operation reduces $H_s\cap R_s$. We call
this move an {\em unlinking along $R_s$ and $S_s$}. We also say that
this unlinking is {\em performed in $H_t$}. We emphasize that $H_t$
remains unchanged and that the resulting $H_s$ still is contained in
$H_t$.

The following summarizes properties of an unlinking that are either
clear or were already mentioned above.

\begin{lemma}\label{L:preserve}
Consider $H_s'$ obtained from $H_s$ by an unlinking along $R_s$ and
$S_s$ realized in $H_t$. Then:
\begin{enumerate}[\upshape (i)]
    \item $|H_s'\cap R_s|<|H_s\cap R_s|$,
    \item $H_t\subseteq H_1$ remains unchanged and $H_s'\subseteq H_t$
    is isotopic to $H_s$ in $M$, and
    \item the spots ${D}_s$ of $K_s$ also remain unchanged, therefore
    $H_s'\cap\mathcal{E}_1=H_s\cap\mathcal{E}_1$.
\end{enumerate}
\end{lemma}

\begin{proposition}\label{P:unlinking}
Let $0=t_0<t_1<\dots<t_{k-1}<t_k=1$ be a complete increasing sequence
of regular values. Suppose that $\Delta$ is a half-disc for
$H_0\subseteq H_1$, standard with respect to $\{\,t_i\,\}$.  Then
there exists a sequence of unlinkings taking $H_0$ to $H_0'\subseteq
H_1$ for which $\Delta$ is an unholed half-disc. Moreover
$H_0'\cap\mathcal{E}_1=H_0\cap\mathcal{E}_1$.
\end{proposition}

\begin{proof}
Let $(R_{t_i},\gamma_{t_i}, \gamma_{t_{i+1}})$, $t_0\leq t_i\leq
t_{m-1}$ and $\Delta_{t_m}$ be the rectangles and half-disc contained
in $\Delta$ determined by $\{\,t_j\,\}$. For each of the rectangles
$R_{t_i}$ we will perform a sequence of unlinkings along $R_{t_i}$ and
a carefully chosen sphere $S_{t_i}$ attached to one side of
$R_{t_i}$. The effect will be to unlink $K_{t_i}$ in the $K_{t_{i+1}}$
containing $R_{t_i}$.  The unlinkings will be performed successively
using induction on $i$, this time with $i$ increasing. At the end of
the inductive process, all rectangles will be unholed, leading to an
unholed $\Delta$ (recall that $\Delta_{t_m}$ is unholed from the start
and this property will be preserved).

As we have mentioned before, an important part of the argument will be
aimed at the following. Each $R_{t_i}$ is contained in a $K_{t_{i+1}}$
and intersects a $K_{t_i}$. It is necessary to find a spot $D_{t_i}$
of $K_{t_i}$ which 1) intersects $R_{t_i}$ and 2) extends to a dual
sphere $S_{t_i}$ with the property that $S_{t_i}\cap K_{t_{i+1}}=
D_{t_{i+1}}$ for some spot $D_{t_{i+1}}$ of $K_{t_{i+1}}$ (see Lemma
\ref{L:spots} below). This means that one side of $R_{t_i}$ intersects
a spot of $K_{t_{i+1}}$ containing just one spot of $K_{t_i}$.  Recall
that this property is necessary to perform an unlinking.

Let $K$ be the component of $\widehat{H}_1$ that contains $\Delta$ and
let $K_0$ be the component of $H_0\cap K$ that contains the upper
boundary $\alpha$ of $\Delta$.  We assume that $\alpha$ intersects two
distinct $1$-handles of $H_0$ and thus a $0$-handle. If it intersects
just one $1$-handle (hence no $0$-handle) then we replace a disc of
$\mathcal{E}_0$, the one corresponding to this $1$-handle, by two
parallel copies, which has the effect of introducing a $0$-handle in
$K_0$.  This amounts to choosing a new system ${\cal E}$ a, and there
are choices in how $f({\cal E}_0)={\cal E}_1$ intersects $H_0$ in
discs parallel to discs of ${\cal E}_0$, i.e. choices for the
incidence matrix, see Section \ref{Laminations}.

Our modification ensures that each $K_{t_i}\supseteq K_0$, $t_i>0$,
corresponds to a 0-handle, with 1-handles giving distinct spots on
$\bdry K_t$. Therefore each rectangle $R_{t_i}$ has the following
property: the spots of $K_{t_i}$ that $\gamma_{t_i}\subseteq\partial
R_{t_i}$ connects correspond to distinct discs of
$\mathcal{E}_{t_i}$. That property will be important in finding the
sphere $S_{t_i}$ intersecting $K_{t_{i+1}}$ at a spot.

Consider two consecutive terms $s<t$ of the sequence (so there is
precisely one singular value between them). We need the following.

\begin{lemma}\label{L:spots}
A holed rectangle $R_s$ intersects a spot of $K_s$ with the property
that the corresponding dual sphere intersects $K_t$ in a single spot.
\end{lemma}
\begin{proof}
It clearly suffices to show that $R_s$ intersects a spot $D_s$ with
the following property. If $D_t$ is the spot of $K_t$ containing $D_s$
then $K_s\cap D_t=D_s$.

Recall that we are assuming that no disc of $\mathcal{E}_1\subseteq
H_1$ is represented twice as a spot of a single component $K_1$ of
$\widehat{H}_1$.

There are three cases.  The case analysis is exactly the same as in
the proof of \ref{ParallelLemma}, pictured in Figure
\ref{Aut3Splitting}. The simplest and most common is Case 1, where a
single component $K_s$ of $H_s\cap K$ is contained in $K_t$ and with
the same isotopy type (in $K$) as $K_t$. Another possibility is Case
2; that there are two distinct $K_s$'s in $K_t$. In Case 3, there is a
single $K_s$ in $K_t$ but $K_s$ is not isotopic to $K_t$.  In Cases 1
and 3, there is only one $K_s$ in $K_t$, so in fact there is no real
need for unlinking; there may be a rectangle, but there is no danger
that it is holed by a $K'_s$ in $K_t$, see Remark
\ref{UnholedHalf-Disc}.  Therefore we only need to deal with Case 2.

In Case 2 there are two distinct components of $H_s\cap K$ contained
in $K_t$: $K_s$ and $K'_s$. Here one spot $D_t$ of $K_t$ contains
precisely one spot of each of $K_s$, $K'_s$.  Denote these by $D_s$
and $D'_s$ respectively. The other spots of $K_t$ contain a single
spot, of either $K_s$ or $K'_s$. Now all spots of $K_s$, $K'_s$ other
than $D_s$, $D'_s$ extend to distinct spots of $K_t$.  But each of
$K_s$, $K_s'$ contains just one of these exceptional spots. Recalling
that $\gamma_s$ connects two distinct spots, at least one will have
the desired property.  The illustration of Case 2 in Figure
\ref{Aut3Splitting} is somewhat misleading, since it is not easy to
picture a holed rectangle there. Figure \ref{Aut3Linked} is more
realistic.  It shows a rectangle $R_s$ essentially intersected by
$K_s'$.

\begin{figure}[ht]
\centering
\psfrag{K}{\fontsize{\figurefontsize}{12}$K$}\psfrag{K'}{\fontsize{\figurefontsize}{12}$K'$}
\psfrag{R}{\fontsize{\figurefontsize}{12}${\hskip1pt}R$}
\psfrag{s}{\fontsize{\figuresmallfontsize}{12}${\hskip1pt}\raise1.8pt\hbox{$s$}$}
\psfrag{t}{\fontsize{\figuresmallfontsize}{12}${\hskip-1pt}\raise2.5pt\hbox{$t$}$}

\scalebox{1.0}{\includegraphics{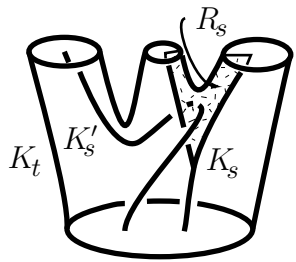}} \caption{\small
Linking and the unlinking move.} \label{Aut3Linked}
\end{figure}

\end{proof}

We continue the proof of Proposition \ref{P:unlinking}. We are
considering two consecutive terms $s<t$ of the sequence. Suppose that
for any $t_i<s$ the rectangle $R_{t_i}$ is unholed. Suppose also that
this process of removing holes did not alter $H_s$.  Therefore the
spheres dual to $H_s$ are preserved and so the conclusions of Lemma
\ref{L:spots} above still hold for $R_s$.

The goal now is to prove the induction step and remove components
(holes) of $H_s\cap R_s$ using unlinkings. We shall do this while also
preserving $H_t$. Consider $K_s$ containing $\gamma_s$. The arc
$\gamma_s\subseteq\partial R_s$ connects two spots of $K_s$. As
mentioned in the previous paragraph, at least one of these has the
property that the corresponding sphere $S_s$ intersects $K_t$ at a
spot, therefore an unlinking along $R_s$ and $S_s$ is possible.

We consider components of $H_s\cap R_s$ disjoint from $\gamma_s$.
Note that these are essential discs in $H_s$. Indeed, $\Delta$ was
standard in the beginning.  This was true before any unlinking moves
were done because $\Delta$ was standard.  But our inductive procedure
ensures that $H_s$ remains unchanged as we do unlinkings in lower
levels, therefore the property is preserved. We perform the unlinking
move to remove such a disc component of $H_s\cap R_s$.  As stated in
Lemma \ref{L:preserve}, this reduces $|H_s\cap R_s|$. The idea is to
continue performing unlinking moves inductively until $H_s\cap
R_s=\gamma_s$. Here again we note that unlinkings preserve spots of
$K_s$ and $K_t$ (Lemma \ref{L:preserve}), therefore the property in
the conclusion of Lemma \ref{L:spots} is preserved.

Lemma \ref{L:preserve} also ensures that each unlinking does not
introduce intersections of $H_0$ with the other rectangles $R_{t_i}$,
$t_i<s$, and that $\Delta$ is preserved as a half-disc.  It also shows
that the resulting $H_0\subseteq H_s$ is contained in $H_1$ and that
$H_0\cap\mathcal{E}_1$ is left unchanged. We now have all $R_{t_i}$
for $t_i\leq s$ unholed, preserving the other desired properties of
the conclusion of the proposition.  This proves the induction step.

We recall that $\Delta_{t_m}$ is unholed originally. Since unlinkings
along $R_{t_i}$, $t_i<t_m$ preserve $H_{t_m}$, $\Delta_{t_m}$ remains
unholed, completing the proof.
\end{proof}

\hop

So far, we have introduced the unlinking move and used it to obtain an
unholed half-disc out of a holed one. As mentioned before, we omitted
any reference to the automorphism $f\colon M\to M$, which is the
object of interest. We now use Proposition \ref{P:unlinking} to
obtained the desired unholed half-disc by isotopies of $f$.

The reader may have anticipated the idea for dealing with $f$:  A
sequence of unlinkings is isotopic to the identity, hence the
composition is realizable by an isotopy of $f$. The catch here is that
it is not obvious that the resulting automorphism preserves some
important features of the original one. Among these is the property of
preserving the canonical Heegaard splitting of $M$, restricting to an
automorphism of $H$. It is also needed that this handlebody
automorphism is outward expanding with respect to $H'_0$ and that the
growth rate does not increase (in fact, it will remain unchanged).

Recall that handlebodies $H\subseteq H'$ are ``concentric'' if
$H'-\intr(H)\simeq(\partial H\times I)$.  The main technical lemma
in this discussion is the following:

\begin{lemma}\label{L:concentric}
In the situation described in Proposition \ref{P:unlinking} we have
that $H_0'\subseteq H_1$ are concentric.
\end{lemma}
\begin{proof}
With the same complete increasing sequence
$0=t_0<t_1<\dots<t_{k-1}<t_k=1$ used in the proof of Proposition
\ref{P:unlinking}, suppose $s,t$ are consecutive terms of the complete
sequence.  Suppose we can show that an unlinking move replacing $H_s$
by $H_s'$ in $H_t$ leaves $H_s'$ concentric in $H_t$.  Then the
concentricity of each $H_{t_i}$ in $H_{t_{i+1}}$ is preserved by the
unlinking move, and inductively, we then show that all the sequence of
unlinkings preserves concentricity.

To prove that an unlinking move replacing  $H_s$ by $H_s'$ in $H_t$
leaves $H_s'$ concentric in $H_t$, we first observe that the
concentricity of $H_s$ in $H_t$ follows from the fact that for all but
one $K_t$, there is just one component $K_s\subseteq K_t$ of $H_s\cap
K_t$.  For all but two components $K_t$, the unique $K_s\subseteq K_t$
is concentric, i.e. there is a product structure for $K_t-K_s$.  The
two exceptional $K_t$'s are illustrated in Figure \ref{Aut3Splitting},
on the right in Cases 2 and 3.  The transition to a picture before the
event, replacing $H_s$ by $H_s'$, which intersects every $K_t$
concentrically, is made by adding a neighborhood of a pinching
half-disc, as illustrated in Figure \ref{Aut3PinchingHalfDisc}.

\begin{figure}[ht]
\centering
\psfrag{H}{\fontsize{\figurefontsize}{12}$H$}\psfrag{H'}{\fontsize{\figurefontsize}{12}${\hskip1.3pt}H'$}
\psfrag{s}{\fontsize{\figuresmallfontsize}{12}${\hskip1.5pt}\raise1pt\hbox{$s$}$}
\psfrag{t}{\fontsize{\figuresmallfontsize}{12}${\hskip-0.7pt}\raise2pt\hbox{$t$}$}
\psfrag{pinching half-disc}{\fontsize{\figurefontsize}{12}$\hbox{pinching half-disc}$}

\scalebox{1.0}{\includegraphics{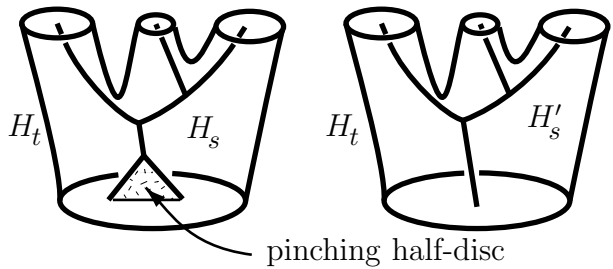}}
\caption{\small The pinching half-disc.}
\label{Aut3PinchingHalfDisc}
\end{figure}

Now we need only observe that the unlinking move does not affect the
pinching half-disc.  Namely, the rectangle $R_s$ can easily be chose
to be disjoint from the pinching half-disc.  Thus after the unlinking
move, we can still add a neighborhood of the pinching half disc to
obtain an $H_s'$ perfectly concentric in each $K_t$.
\end{proof}

\begin{proof}[Proof of Theorem \ref{BackTrackTheorem}]
Let $G$ be the quotient graph determined by the 2-dimensional
lamination invariant under $g$. Then $g$ determines a homotopy
equivalence $h_g\colon G\to G$. By hypothesis, $h_g$ is not
train-track. Consider the algorithm of Bestvina-Handel to simplify the
homotopy equivalence. It clearly suffices to show that we can carry
out the algorithm using isotopies of $f$ that preserve the canonical
splitting.

The following moves of \cite{BH:Tracks} are realizable by isotopies of
$g$ (and hence isotopies of $f$): {\em collapses}, {\em valence-one}
and {\em valence-two homotopies} and {\em subdivisions}. We refer the
reader to \cite{UO:Autos} for details on these isotopies.

We shall consider the other moves, {\em folds} and {\em
pulling-tights}, more carefully. For a fold, two edges $e$, $e'$ of
$G$ incident to the same vertex are mapped to the same edge. In the
context of $g\colon H\to H$ this corresponds to the following
situation.  There is a component $K$ of $\widehat{H}_1$ and a
$K_0\subseteq H_0\cap K$ with the property that two distinct spots
$D_0$, $D_0'$ of $K_0$ are contained in the same spot $D_1$ of $K$,
see Figure \ref{Aut3Fold}.

\begin{figure}[ht]
\centering
\psfrag{G}{\fontsize{\figurefontsize}{12}$G$}\psfrag{K}{\fontsize{\figurefontsize}{12}$K$}
\psfrag{D}{\fontsize{\figurefontsize}{12}$D$} 
\psfrag{D '}{\fontsize{\figurefontsize}{12}$D'$}\psfrag{e}{\fontsize{\figurefontsize}{12}$e$}
\psfrag{e'}{\fontsize{\figurefontsize}{12}$e'$}

\scalebox{1}{\includegraphics{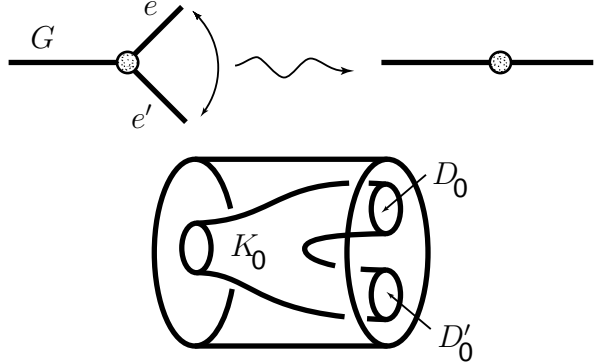}} \caption{\small A fold in
the quotient graph and the corresponding $H_0$ in $H_1$.}
\label{Aut3Fold}
\end{figure}

The spots $D_0$, $D_0'$ correspond to the edges $e$, $e'$ and $D_1$ to
the edge to which they are mapped. By Proposition \ref{P:half-disc},
there exists a standard half-disc $\Delta$ for $H_0$ connecting these
two spots.

If $\Delta$ is unholed, the fold is realizable by a {\em down-move}
\cite{UO:Autos}, which is an isotopy of $g$. Roughly speaking, such a
move changes $g$ so that the disc corresponding to the spot $D_1$,
which is the image of a disc $E\subseteq\mathcal{E}_0$, is isotoped
along $\Delta$. The two discs of $\mathcal{E}_0$ corresponding to the
spots $D_0$, $D_0'$ are replaced by their band-sum along the upper
boundary of $\Delta$. Such a move realizes the desired fold.

If $\Delta$ is holed then we use Proposition \ref{P:unlinking} and
obtain $H_0'$. Lemma \ref{L:preserve} assures that $H_1$ (and the
canonical splitting) are preserved. Lemma \ref{L:concentric} shows
that $H_0'\subseteq H$ and $H$ are concentric, so the resulting automorphism of
$M$ can be adjusted to preserve the canonical splitting $M=H\cup H'$.
Moreover the resulting restriction $g'\colon H\to H$ can be made
outward expanding with respect to $H_0'$ (see Section
\ref{Laminations}). The assertion of Proposition \ref{P:unlinking}
that $H_0'\cap\mathcal{E}_1=H_0\cap\mathcal{E}_1$ implies that the
original disc system yields invariant laminations with the same growth
rate. In fact more holds. Let $H\to G'$ be the new quotient map and
$h_{g'}\colon G'\to G'$ the new homotopy equivalence. Then $h_{g'}$
and $h_g$ are conjugate by a homeomorphism $G'\to G$.  Therefore the
half-disc $\Delta$ still corresponds to the fold being considered. But
now $\Delta$ is not holed so we fall in the previous case and realize
the fold.

For a pulling-tight the argument is similar, so we will skip
details. In this case there exists an edge of $G$ in whose interior
$h_g$ fails to be locally injective. As before, in the setting of the
handlebody $H$ that corresponds to a $K_0\subseteq K$ with two
distinct spots in the same spot of $K$ (in this case $K_0$ will have
just two spots). So, here again, that corresponds to a standard
half-disc $\Delta$. If it is not holed we perform a {\em diversion}
\cite{UO:Autos}, which realizes the pulling-tight move. If it is holed
we remove the holes with unlinkings and return to the previous case.
\end{proof}

\begin{remark}
We can say more about the resulting $f'$ and $g'=f'|_{H}$ in the
conclusion of Theorem \ref{TrainTrackTheorem}. Bestvina-Handel's
algorithm stops when either 1) the homotopy equivalence admits no
back-tracking and the corresponding incidence matrix is irreducible or
2) when the incidence matrix, after collapsing of invariant forests,
is reducible. In this last case $g'$ is either reducible or
periodic. In the first case, if the growth $\lambda=1$ then $g'$ may
be isotoped to permute (transitively) the discs of
$\mathcal{E}_0$. This means that a power $(g')^k$ fixes
$\mathcal{E}_0$. Therefore $(f')^k$ fixes the sphere system
$\mathcal{S}_0\subseteq M$ corresponding to $\mathcal{E}_0$. By
reducing $(f')^k$ along $\mathcal{S}_0$ we obtain automorphisms of
holed balls preserving each boundary component, therefore isotopic to
the identity. This proves that $(f')^k$ is isotopic to a composition
of twists on spheres of $\mathcal{S}_0$. In particular $(f')^{2k}$ is
isotopic to $Id_{M}$ (although in general $g'$ will not be
periodic). If $\lambda>1$ then $g'$ may be either reducible or
generic.  If it is generic, as stated in the theorem, it is
train-track.  If it is reducible, we believe we can prove that it is a
lift of a pseudo-Anosov automorphism of a surface $F$ to the total
space of an $I$-bundle over $F$ (in this case $f$ restricted to the
splitting surface ${\partial H}$ is reducible).  The proof of this
last statement is not written formally yet.
\end{remark}

\section{Automorphisms of mixed bodies} \label{MixedBodies}

In this section we suppose $M$ is a mixed body.  Automorphisms of
holed mixed bodies can be understood by the same method.  In fact,
recall from the introduction that the ``difference" between an
automorphism of a holed or spotted manifold and an unspotted manifold
is an automorphism of the spotted manifold which can be realized by an
isotopy of the unspotted manifold $M$.

Recall from Theorem \ref{CanonicalHeegaardTheorem2} that there exists
a canonical closed surface $F\subseteq M$ separating $M$ into a
handlebody $H$ and a compression body $Q$. The splitting $M=H\cup_F Q$
has the property that every curve in $F$ bounding a disc in $Q$ also
bounds a disc in $H$. A curve in $F$ bounding a disc in $H$ may either
bound a disc in $Q$ or be parallel to $\partial_i Q=\partial M$. As in
the special case treated in Section \ref{TrainTrackSection} we can use
the canonicity of $F$ to isotope an automorphism $f\colon M\to M$ to
preserve the splitting, so it induces restricted automorphisms on $Q$
and $H$.  Here, again, $F$ is not ``rigid'' (e.g., a double twist on a
sphere $S$ intersecting $F$ on an essential curve is isotopic to the
identity, see Example \ref{NonUniqueExample}).

But, in this case, there will be another canonical surface $V\subseteq
H$ which can also be made invariant. From another point of view, we
shall prove that the restriction of $f$ to $H$ is rigidly
reducible. More specifically, we shall prove that if $f$ preserves
$F$, then $f|_{H}$ always preserve a compression body $Q'\subseteq H$
with $\partial_e Q'=\partial H$ and $\partial_i Q'\neq\emptyset$ (see
Section \ref{Intro}). Roughly speaking, $Q'$ is the ``mirror image''
of $Q$ through $F$ and $V$ is the mirror image of $\partial
M=\partial_i Q$.

\begin{proposition}\label{FurtherReducingSurface}
Let $f$ be an automorphism of $M$. Then $f$ can be isotoped to
preserve the canonical splitting $M=H\cup_F Q$. Moreover, $f$ can be
further isotoped rel $F$ to preserve a surface $V\subseteq H$ with the
following properties. It separates $H$ into a (possibly disconnected)
handlebody $H'$ and a compression body $Q'$ such that a curve in $F$
bounds a disc in $Q$ if and only if it bounds a disc in $Q'$. Then $f$
preserves $Q$, $Q'$ and $H'$. Also, $f|_{Q}$ and $f|_{Q'}$ are
conjugate, and $f$ restricted to $V$ (and hence $H'$) is unique, in
the sense that it depends only on the isotopy class of $f$, and not on
$f|_H$.
\end{proposition}

To prove the proposition above we consider a compression body
$Q'\subseteq M$ with the following properties: 1) $Q'\subseteq H$, 2)
$\partial_e Q'=\partial H=F$ and 3) a curve in $F$ bounds a disc in
$Q'$ if and only if it bounds a disc in $Q$. We can regard $Q\cup Q'$
as the double of $Q$ along $\partial_e Q$.

\begin{lemma}\label{MirrorCompressionBody}
The compression body $Q'$ is well-defined up to isotopy. In
particular, so is $V=\partial_i Q'$.
\end{lemma}
\begin{proof}
Consider a complete system of essential discs $\mathcal{E}$ in $Q$
i.e., $\mathcal{C}=\partial\mathcal{E}\subseteq F$ and $\mathcal{E}$
cuts $Q$ open into a union of a product $(\partial_i Q)\times I$ and
balls. Then $\mathcal{C}$ also bounds discs $\mathcal{E}'\subseteq
H$. Let $Q'\subseteq H$ be the compression body determined by $F$ and
$\mathcal{E}'$. More precisely, consider a neighborhood $N\subseteq H$
of $F\cup\mathcal{E}'$. Let $B\subseteq H$ be the union of ball
components of $\overline{H-N}$. Then $Q'=N\cup B$.

We need to prove ``symmetry'', i.e. that curves of $F$ bound discs in
$Q'$ if and only if they also bound in $Q$. Here we require that the
system $\mathcal{E}$ is minimal, in the sense that the removal of any
discs yields a system which is not complete.  It is clear that any
essential disc is contained in some such system. It is a standard fact
that any pair of such minimal disc systems in a compression body are
connected by a sequence of {\em disc slides} (see
e.g. \cite{FB:CompressionBody}, appendix B). But a disc slide in any
of the compression bodies $Q$, $Q'$ corresponds to a disc slide in the
other whose results coincide in $F$. We need not consider
$\partial$-parallel discs in the compression bodies, so the proof that
$Q$ and $Q'$ are symmetric is complete.

The above argument with disc slides, together with irreducibility of
$H$, also proves that the construction of $Q'$ does not depend on
$\mathcal{E}$, implying its uniqueness.
\end{proof}

Let $f\colon M\to M$ be an automorphism of a mixed body. Denote by
$\mathcal{A}_f$ the class of automorphisms $f'\colon M\to M$ which are
isotopic to $f$ and preserve $Q$, $Q'$, $H'\subseteq M$.

\begin{lemma}\label{RigidHandlebody}
If $f$, $f'\in\mathcal{A}_f$ then $f|_{V}$ and $f'|_{V}$ are isotopic
(hence also $f|_{H'}$ and $f'|_{H'}$ are isotopic).
\end{lemma}
\begin{proof}
Let $h=(f^{-1}\circ f')\in\mathcal{A}_{\rm Id}$. Therefore
$h|_{\partial M}={\rm Id}_{\partial M}$.

It is a simple and known fact that an automorphism of a compression
body is uniquely determined by its restriction to the exterior
boundary. Recall that $F$ is the exterior boundary of both $Q$ and
$Q'$. Symmetry of $Q$ and $Q'$ implies that $h|_{Q\cup Q'}$ may be
regarded as the double of $h|_{Q}$. Now $\partial_i Q=\partial M$,
then $h|_{\partial_i Q}={\rm Id}$.  Therefore $h|_{\partial_i Q'}={\rm
Id}$, completing the proof.
\end{proof}

\begin{proof}[Proof of Proposition \ref{FurtherReducingSurface}]
Use Lemma\ref{MirrorCompressionBody} and note the following. By
cutting $Q'$ along the cocores of its $1$-handles it is easy to see
that $V=\partial_i Q'$ separates $H$ into $Q'$ and a handlebody
$H'=\overline{H-Q'}$. This proves the conclusions on the existence of
the canonical decomposition. The conclusions on $f|_{V}$ and $f|_{H'}$
follow from Lemma \ref{RigidHandlebody}.  The conclusion on conjugacy
of $g=f|_{Q}$ and $g'=f|_{Q'}$ is clear from the symmetry of $Q$ and
$Q'$.
\end{proof}

\begin{remark}
We emphasize that while $V$ is a rigid reducing surface for $f|_{H}$
it is not a rigid reducing surface for $f$, even with unique $f|_{H'}$
and $f|_{V}$. The reason is that $f$ need not restrict uniquely to the
complement $Q\cup Q'$.

Recall from the introduction that some of the automorphisms $f\colon
M\to M$ being considered arise from adjusting automorphisms of general
reducible $3$-manifolds. In this case $f|_{\partial M}$ is periodic,
permuting components. The argument above in Lemma
\ref{RigidHandlebody} shows, then, that in this case $f$ will always
be periodic on the (possibly disconnected) handlebody $H'\subseteq
H\subseteq M$. But we also consider cases without the periodic
hypothesis on $f|_{\partial M}$. Here there is a constraint on
$f|_{\partial M}$. To see this recall from the argument that
$f|_{\partial M}$ determines $f|_{V}$ uniquely. But $V=\partial H'$,
therefore $f|_{V}$ is a surface automorphism that extends over the
handlebody $H'$, which is a constraint.
\end{remark}

The goal is to find a nice representative for the class of an
automorphism $f\colon M\to M$. Use Proposition
\ref{FurtherReducingSurface} to decompose $f\colon M\to M$ into
automorphisms $f|_{Q}$, $f|_{Q'}$ and $f|_{H'}$. The first two are
conjugate and they will be dealt with using unlinkings analogous to
those of Section \ref{TrainTrackSection}. The last is dealt with in
\cite{UO:Autos} (see also Sections \ref{Intro} and \ref{Laminations}).

We now proceed to find desirable restrictions of $f$ to the canonical
compression bodies $Q$, $Q'$. The two restrictions are conjugate so it
suffices to deal with one of them. We choose $Q'$ because $Q'\cup H'$
is a handlebody, what shall be of a technical relevance. So consider $g=f|_{Q'}$.  We need to define a
convenient notion of ``back-tracking'' for an automorphism of
$Q'$. Here suppose that $g\colon Q'\to Q'$ is generic. Similarly to
what is done in Section \ref{TrainTrackSection}, consider concentric
$Q_i\subseteq Q'$, disc systems $\mathcal{E}_i\subseteq Q_i$,
laminations $\Lambda$ and $\Omega$, and growth rate $\lambda$ as in
Section \ref{Laminations}. Also assume that $\Lambda$ has the
incompressibility property. Consider $\widehat{Q}_1=Q_1|\mathcal{E}_1$
and let $K_1$ be a component.  Recall that $(K_1,V_1,\mathcal{D}_1)$
denotes either a spotted ball or spotted product, with spots
determined as duplicates of $\mathcal{E}_1$.

The new phenomenon for back-tracking in $Q'$ is roughly the
following. If one follows a leaf of the mapped 1-dimensional
lamination it may enter the product part of $Q'$ through a spot $D$
and also leave through $D$. Such a leaf may ``back-track'' in $Q'$, in
the sense that it can be pulled tight (homotopically) in $Q'$ to
reduce intersections with that spot. Other (portions of) leaves may
not admit pulling tight in this way, being somehow linked with
$\partial_i Q'$. But among these some can still be pulled tight if
they are allowed to pass through the handlebody $H'$ (recall that
$H=Q'\cup H'$).

\begin{defn}\label{BackTrackingCB}
Suppose that $\sigma\subseteq\Omega$ is an arc such that
$\sigma'=\omega(\sigma)\subseteq Q_0\subseteq Q'$ has the following
properties: 1) $\partial\sigma'$ is contained in a single disc $E_0$
of the system $\mathcal{E}_0\subseteq Q_0$, and 2) $\sigma'$ is
isotopic in $H$ (rel $\partial\sigma'$) to an arc in $E_0$. If there
exists such an arc we say that $\Omega$ {\em back-tracks}. In this
case we say that $g=f|_{Q'}$ {\em admits back-tracking in $H$} or just
{\em admits back-tracking}.
\end{defn}
\begin{remark}
Property 1) above is intrinsic to $Q'$ while property 2) is not, since
it refers to the entire handlebody $H$.
\end{remark}

\begin{lemma}\label{L:back-track}
If $\Omega$ back-tracks then there exists an integer $i\leq 0$ such
that $Q_i$ has the following properties. There is a component $K_i$ of
$Q_i\cap K_1$ which is a twice spotted ball (in other words, $K_i$ is
contained in a 1-handle of $Q_i$), both contained in the same spot
$D_1$ of $\mathcal{D}_1$. Moreover if $(P_i,\mathcal{D}_i)$ is the
unique \spinal for $Q_i$, then it is parallel in $K_1\cup H'$ to the
subspinal pair $(D_1,D_1)$ for $(K_1,V_1,\mathcal{D}_1)$, where $D_1$
is a disc of $\mathcal{D}_1$.
\end{lemma}

We skip the proof, which follows directly from unknottedness and
parallelism (see Section \ref{Laminations}).

\begin{remark}
The converse is obviously true.
\end{remark}

In the situation of the lemma above we say that $Q_i$ {\em back-tracks
in $Q_1$}. More generally, we say that {\em $Q_i$ back-tracks in
$Q_{j+1}$}, $i<j$, if $g^{-j}(Q_i)$ back-tracks in
$g^{-j}(Q_{j+1})=Q_1$.

\begin{defn}
If $f\in\mathcal{A}_f$ and $g=f|_{Q'}$ is generic and admits no
back-tracking in $H$ then we say that $g$ is {\em train-track generic
in $H$}, or just {\em train-track generic}.
\end{defn}

The importance of considering train-track generic automorphisms comes
from the result below. Its proof is an adaptation of an argument of
\cite{LNC:Tightness}.

\begin{thm}\label{MinimalGrowth}
Let $f$ be an automorphism of a mixed body $M$. If $g=f|_{Q'}$ is
train-track generic in $H$ then its growth rate is minimal among the
restrictions to $Q'$ of all automorphisms in $\mathcal{A}_f$.
\end{thm}

\begin{proof}
We prove the contrapositive, so assume that there exists
$f'\in\mathcal{A}_f$ such that $g'=f'|_{Q'}$ has growth $\lambda'$
less than the growth $\lambda$ of $g$.

Consider the original disc system $\mathcal{E}_0\subseteq Q_0$,
determining a dual object $\Gamma_0$. For simplicity denote it by
$\Gamma$.  Recall that $\Gamma=\hat\Gamma\cup V$, where $\hat \Gamma$
is a graph with one vertex on each component of $V$.  Consider the
mapped-in 1-dimensional lamination $\omega\colon\Omega\to Q_0$ (see
Section \ref{Laminations}, Theorem
\ref{CompressionBodyClassificationTheorem}). Recall that $\nu$ denotes
its transverse measure. Use $(\Omega,\nu)$ to give weights $\hat\nu$
to the edges of $\Gamma$. More precisely, if $E_i\subseteq
\mathcal{E}_0$ is the disc corresponding to an edge $e_i$ of $\Gamma$
then the weight on $e_i$ is $\int_{(\omega^{-1}(E_i))}\nu$. For
simplicity, call the pair $(\Gamma,\hat\nu)$ a {\em weighted graph},
although it is really $\hat \Gamma$ which is a weighted graph.

Extend $\mathcal{E}_0$ to a disc system $\mathcal{E}\subseteq Q'$
through the product in $Q'-\intr({Q}_0)$. Apply powers of $g^{-1}$ to
$\Gamma$. One can ``push forward'' the weights on $\Gamma$ to
$\Gamma_n=g^{-n}(\Gamma)$ through $g^{-n}$. More precisely, the weight
on the image $g^{-n}(e)$ of an edge $e$ is the same weight as that on
$e$. Denote the resulting weighted graph as $(\Gamma_n,\hat\nu)$. It
follows that
\begin{equation}\label{Eq:eigenvalue}
\frac{\weight{\mathcal{E}}{\Gamma_n}{\hat\nu}}{\lambda^n}=c>0
\end{equation}
is constant, where $\weight{\mathcal{E}}{\Gamma_n}{\hat\nu}$ denotes
the sum of points of intersection $\mathcal{E}\cap\Gamma_n$ weighted
by $\hat\nu$.

Now consider $f'$ and $g'=f'|_{Q'}$, with corresponding $Q_0'$, disc
system $\mathcal{E}'_0\subseteq Q_0'$ with dual graph $\Gamma'$ and
growth $\lambda'$. As before, the corresponding 1-dimensional
lamination determines a weighted graph $(\Gamma',\hat\nu')$.

But $Q_0$, $Q_0'$ are concentric in $Q'$, therefore one can apply an
ambient isotopy rel $V$ taking $\Gamma$ to some $\bar\Gamma\subseteq
Q_0'$. We can also assume that $\bar\Gamma$ does not ``back-track'' in
the 1-handles of $Q_0'$, in the sense that no vertex of $\bar\Gamma$
is contained in such a 1-handle $K$ and that the edges of $\bar\Gamma$
are transverse to the foliation of $K$ by cocores.  The isotopy from
$\Gamma$ to $\bar\Gamma$ pushes forward $\hat\nu$ to yield
$(\bar\Gamma,\hat\nu)$.

We are interested in applying powers of $(g')^{-1}$ to
$\bar\Gamma$. Denote $\bar\Gamma_n=(g')^{-n}(\bar\Gamma)$. By the
construction of $\bar\Gamma\subseteq Q_0'$,
\[
\frac{\weight{\mathcal{E}}{\bar\Gamma_n}{\hat\nu}}{(\lambda')^n}<\bar{C}
\]
is bounded.  This is true because the growth rate of intersections
with any essential disc with the images under $(g')^{-n}$ of {\em any}
weighted graph is $\lambda'$.  Since $\lambda'<\lambda$,
\begin{equation}\label{Eq:convergence}
\frac{\weight{\mathcal{E}}{\bar\Gamma_n}{\hat\nu}}{(\lambda)^n}\to 0.
\end{equation}
By combining (\ref{Eq:eigenvalue}) and (\ref{Eq:convergence}) it
follows that there exists $N$ such that
\begin{equation}\label{Eq:inequality}
\frac{\weight{\mathcal{E}}{\bar\Gamma_N}{\hat\nu}}{\lambda^N}<
\frac{\weight{\mathcal{E}}{\Gamma_N}{\hat\nu}}{\lambda^N}
\end{equation}

But recall that there is an isotopy rel $V$ taking $\Gamma$ to
$\bar\Gamma$. There is also an isotopy (not necessarily rel $V$)
taking $g^{-N}(\bar\Gamma)$ to $\bar\Gamma_N$ ($f$ and $f'$, which
restrict to $Q'$ as $g$ and $g'$, are isotopic). Therefore there is an
isotopy $\iota$ taking $\Gamma_N$ to $\bar\Gamma_N$.  Note that this
is an isotopy in $M$, not just $Q'$.

Regard $Q_{-N}$ as a neighborhood of $\Gamma_N$. Extend the isotopy
$\iota$ to $Q_{-N}$, taking it to a neighborhood $\bar Q_N$ of
$\bar\Gamma_N$. Define $\bar\omega=\iota\circ\omega \colon\Omega \to
Q_{-N}$. The inequality (\ref{Eq:inequality}) is restated in terms of
laminations as
\[
\int_{(\bar\omega)^{-1}(\mathcal{E})}\nu<\int_{\omega^{-1}(\mathcal{E})}\nu.
\]

Now let $\mathcal{S}\subseteq M$ be a ``dual sphere system''
corresponding to $\mathcal{E}$, i.e. $\mathcal{S}\cap
Q'=\mathcal{E}$. The inequality above clearly implies
\begin{equation}\label{Eq:inequality2}
\int_{(\bar\omega)^{-1}(\mathcal{S})}\nu<\int_{\omega^{-1}(\mathcal{S})}\nu.
\end{equation}

The isotopy $\iota$ determines a map $p\colon \Omega\times I\to M$
such that $p|_{\Omega\times\{0\}}=\omega$ and
$p|_{\Omega\times\{1\}}=\bar\omega$. Perturb it to put it in general
position. Consider $p^{-1}(\mathcal{S})$. Let $L\subseteq\Omega$ be a
leaf that embeds through $\omega$. Then $(L\times
I)\subseteq(\Omega\times I)$ is a product $\mathbb{R}\times I$ and
$(L\times I)\cap p^{-1}(\mathcal{S})$ consists of embedded closed
curves and arcs. Ignore closed curves for the moment.  Let $A_L$ be
the set of such arcs $\alpha$ intersecting the lower boundary
$L\times\{0\}$. Recall that the set of leaves that do not embed has
measure zero. If for all leaves $L$ that embed and all $\alpha\in A_L$
the other endpoint of $\alpha$ is in $L\times\{1\}$ then
\[
\int_{(\bar\omega)^{-1}(\mathcal{S})}\nu\geq\int_{\omega^{-1}(\mathcal{S})}\nu,
\]
contradicting inequality (\ref{Eq:inequality2}). Therefore there
exists $L$ that embeds and $\alpha\in A_L$ with both endpoints at
$L\times\{0\}$.

Let $\alpha$ be an edgemost such arc. Then there is an arc
$\beta\subseteq (L\times\{0\})$ such that
$\partial\beta=\partial\alpha$. Let $D$ be the half-disc bounded by
$\alpha\cup \beta$ in $L\times I$.  Then $p|D$ is a map which sends
$\alpha$ to a sphere $S$, $S \subseteq\mathcal{S}$, and embeds $\beta$
in the leaf $L$.  We can homotope the map $p|D$ rel $\beta$ to obtain
a map $d:D\to M$ such that $d(\alpha)\subseteq E_0\subseteq S$, where
$E_0\subseteq \mathcal{E}_0$ is the disc in $S$.  Thus one has
$d(\alpha\cup\beta)\subseteq Q_0$.

Recall the Heegaard surface $F$ separates $M$ into $Q$ and $H$.  Let
$d$ be transverse to $F$ and let $\Delta=d^{-1}(Q)\subseteq $D$ $. On
$\Delta$ redefine $d$ to be its mirror image $r\circ d$, where $r$ is
the reflection map $r:Q\to Q'$.  Then homotope $d$ by pushing slightly
into $Q'$.  By post-composing $d$ with an isotopy along the $I$-fibers
of the product $Q'-\intr({Q}_{-1})$ one obtains $d$ mapping into
$Q_0\cup H'$.

Using the irreducibility of $Q_0\cup H'$ (which is a handlebody)
we remove closed curves from $d^{-1}(\mathcal{E}_0)$, which will
then consist of $\alpha$. By unknottedness (Theorem
\ref{UnknottedTheorem}), the embedded image of
$d|\beta\subseteq\Omega$, which lies in a leaf, is parallel to the
subspinal pair $(D_0,D_0)$ for some spotted component
$(K_0,V_0,\mathcal{D}_0)$, where $D_0$ is a duplicate of $E_0$,
proving that $\Omega$ back-tracks.
\end{proof}

\begin{thm}[Compression Body Train Track Theorem]\label{CompressionTrainTrackTheorem}
Let $f\colon M\to M$ be an automorphism of a mixed body. Then there
exists $f'\in\mathcal{A}_f$ such that $g=f'|_{Q'}$ is

1) periodic,

2) rigidly reducible or

3) train-track generic in $H$.

On the other pieces, $Q$ and $H'$, we have:

4) $f'|_{Q}$ is isotopic to a conjugate of $g$, and

5) the class of $f'|_{H'}$ does not depend on
$f'\in\mathcal{A}_f$.
\end{thm}

The conclusions concerning $f'|_{Q}$, $f'|_{H'}$ are proved in
Proposition \ref{FurtherReducingSurface}. We proceed now to prove that
one can always find $f'\in\mathcal{A}_f$ such that $f'|_{Q'}$ is
either periodic, reducible or train-track generic, i.e. admits no
back-tracking. To do so we need to perform unlinkings analogous to
those of Section \ref{TrainTrackSection}. As before, assume that
$g=f|_{Q'}$ is generic and that the leaves of the 2-dimensional
lamination are in Morse position with respect to the height function
determined by the product structure in $Q'-\intr(Q_0)$. We refer
the reader to Section \ref{Laminations} and \cite{UO:Autos} for
details.

For $t\leq 1$ consider spotted balls or products $(K_t,\mathcal{D}_t)$
determined by $Q_t\cap \widehat{Q}_1$.

As before, to define linking we need half-discs. Let
$(\Delta,\alpha,\beta)$ be a half-disc and $(R,\gamma,\gamma')$ a
rectangle (Definitions \ref{HalfDisc}).

\begin{defns}
Let $s<t$. A {\em half-disc for $Q_s$ in $Q_t$} is a half-disc
$(\Delta,\alpha,\beta)$ satisfying:
\begin{enumerate}[\upshape (i)]
\item $\Delta\subseteq K_t\cup H'$ is embedded for some $K_t$,
\item $\beta$ is contained in a spot $D_t$ of $K_t$,
\item $\alpha\subseteq\partial K_s$ for some $K_s\subseteq K_t$ and
$\Delta\cap K_s=\alpha$,
\item $\alpha$ connects two distinct spots of $(K_s,\mathcal{D}_s)$.
\end{enumerate}

As before we refer to $\alpha$ as the {\em upper boundary of
$\Delta$}, denoted by $\partial_u\Delta$, and $\beta$ as its {\em
lower boundary}, denoted by $\partial_l\Delta$.

A {\em rectangle for $Q_s$ in $Q_t$} is a rectangle
$(R_s,\gamma_s,\gamma_t)$ satisfying:

\begin{enumerate}[\upshape (a)]
    \item $R_s\subseteq K_t\cup H'$ is embedded for some $K_t$,
    \item $\gamma_s\subseteq\partial K_s$ for some $K_s\subseteq K_t$
    and $R_s\cap K_s=\gamma_s$,
    \item $\gamma_s$ connects distinct spots of $(K_s,\mathcal{D}_s)$,
    \item $\gamma_t\subseteq \partial K_t -
    \intr(\mathcal{D}_{t_{i+1}})- \intr(V_t)$ (recall that when $K_t$
    is a spotted ball then $V_t=\emptyset$),
    \item $\partial R_s-(\gamma_s\cup\gamma_t)\subseteq\mathcal{D}_t$.
\end{enumerate}
\end{defns}

Note that a half-disc $\Delta$ in $Q_t$ has the property that
$\partial\Delta\cap\partial_i Q'=\emptyset$. Similarly $\partial
R\cap\partial_i Q'=\emptyset$, where $R$ is a rectangle in $Q_t$.

Also, holes now may appear due not just to intersections with $Q_s$
but also with $H'$.

\begin{defn}\label{D:standard2}
Let $0\leq s<t\leq 1$ and $\Delta$ a half-disc for $Q_s$ in $Q_t$.
Suppose that $s=t_0<t_1<\cdots<t_k=t$ is a sequence of regular
values. We say that $\Delta$ is {\em standard with respect to
$\{\,t_i\,\}$} if there is $0\leq m<k$ such that the following
holds. The half-disc $\Delta$ is the union of rectangles $R_{t_i}$,
$s\leq t_i\leq t_{m-1}$, and a half-disc $\Delta_{t_m}$ (see
Figure~\ref{Aut3StandardDisc}) with the properties that:

\begin{enumerate}
    \item each $R_{t_i}$, $s\leq t_i< t_{m}$, is a rectangle for
    $Q_{t_i}$ in $Q_{t_{i+1}}$,
    \item the holes $\intr({R}_{t_i})\cap Q_{t_i}$ are essential
    discs in $Q_{t_m}\cup H'$, each consisting of either an essential
    disc in $Q_{t_m}$ or of a vertical annulus in the product part of
    $Q_{t_m}$ union an essential disc in $H'$,
    \item $\Delta_{t_m}$ is an unholed half-disc for $K_{t_m}$ in
    $H_{t_{m+1}}$.
\end{enumerate}
\end{defn}

The following is the analogue of Proposition \ref{P:half-disc} needed
here.

\begin{proposition}\label{P:half-disc2}
Let $0=t_0<\dots<t_k=1$ be a complete sequence of regular values.
Suppose that a component $K_0$ of $Q_0\cap \widehat{Q}_1$ intersects a
disc $D$ of $\widehat{\mathcal{E}}_1$ in two distinct components
$D^0$, $D^1$. If an arc $\sigma\subseteq K_0$ having an endpoint at
each $D^i$ is parallel (rel $\partial\sigma$) in $\widehat{Q}_1\cup
H'$ to an arc $\sigma'\subseteq D$ then there exists a half-disc
$\Delta$ for $Q_0$ connecting $D^0$ and $D^1$ which is standard with
respect to $\{t_i\}$.
\end{proposition}

\begin{proof}
The proof is similar to that of Proposition \ref{P:half-disc}, using
Theorem \ref{UnknottedTheorem}, Lemma \ref{SpinalPairLemma} and
techniques of Lemma \ref{L:half-disc}.

If $K_1$ is a spotted ball then all  $K_{t_i}$'s containing $\sigma$
also are spotted balls and the result follows directly from the proof
of Lemma \ref{L:half-disc}. In fact, suppose that $t_m$ is the biggest
$t_i$ with the property that $D^0$ and $D^1$ are contained in distinct
spots $D^0_{t_i}$, $D^1_{t_i}$ of $K_{t_i}$. If this $K_{t_m}$ is a
spotted ball then $K_{t_{m+1}}$ also is a spotted ball, as one can
verify by considering the possible events described in Lemma
\ref{EventLemma} changing $K_{t_m}$ to $K_{t_{m+1}}$.  The argument of
Lemma \ref{L:half-disc} still works.

We can therefore assume that $K_{t_m}$ is a spotted product. In this
case $K_{t_{m+1}}$ is also a spotted product.  Let $D'$ be the spot of
$K_{t_{m+1}}$ that contains both $D^0_{t_m}$, $D^1_{t_m}$. There is a
single arc $\gamma_{t_{m+1}}\subseteq D'-\intr(D^0_{t_m}\cup
D^1_{t_m})$ connecting $D^0_{t_m}$ to $D^1_{t_m}$ (up to isotopy in
$D'-\intr(D^0_{t_m}\cup D^1_{t_m})$, rel $\partial(D^0_{t_m}\cup
D^1_{t_m})$).  Recall that there are anambiguous spinal pairs
$(P_{t_m},\mathcal{D}_{t_m})$, $(P_{t_{m+1}},\mathcal{D}_{t_{m+1}})$
for $K_{t_m}$, $K_{t_{m+1}}$ respectively,  and that they are parallel.
We use this parallelism to build $\gamma_{t_m}\subseteq P_{t_m}$
isotopic to $\gamma_{t_{m+1}}$. Note that the parallelism determines
an unholed half-disc $(\Delta_{t_m},\gamma_{t_m},\gamma_{t_{m+1}})$
using the Loop Theorem.

Now proceed inductively ``downwards'', building $\gamma_{t_{i-1}}$
from $\gamma_{t_i}$ in the analogous manner for as long as
$K_{t_{i-1}}$ is a spotted product. More precisely, to each spotted
product  $K_{t_i}$ corresponds a spinal pair
$(P_{t_{i}},\mathcal{D}_{t_{i}})$.  The parallelism of spinal pairs
$(P_{t_{i-1}},\mathcal{D}_{t_{i-1}})$ and
$(P_{t_i},\mathcal{D}_{t_{i}})$ yields $\gamma_{t_{i-1}}$ from
$\gamma_{t_i}$. Note that this parallelism defines a rectangle
$(R_{t_{i-1}},\gamma_{t_{i-1}},\gamma_{t_i})$, which is also holed
only if there is a spotted ball $K_{t_{i-1}}'\subseteq K_{t_i}$
(i.e. in a situation analogous to Case 2 in Figure
\ref{Aut3Splitting}). In any case, $\gamma_{t_{i-1}}$ has the property
of connecting the spots $D^0_{t_{i-1}}$ and $D^1_{t_{i-1}}$ that
contain $D^0$, $D^1$.

If no $K_{t_j}$ containing $\sigma$ is a spotted ball then
$K_0$ is a spotted product and we obtain an arc $\gamma_0\subseteq
P_0$. But $\gamma_0$ clearly connects the spots $D^0$ and $D^1$ of
$K_0$. Hence the union of all rectangles $R_{t_i}$'s and the half-disc
$\Delta_{t_m}$ yields the desired standard half-disc. This takes care
of the case when there is no spotted ball $K_{t_j}$ containing
$\sigma$.

Therefore assume that $\sigma$ is contained in a spotted ball
$K_{t_k}$. Assume that the spotted ball $K_{t_k}$ is maximal, in the
sense that it is contained in a spotted product $K_{t_{k+1}}$. We also
assume from the construction above that there is an arc
$\gamma_{t_{k+1}}\subseteq P_{t_{k+1}}$ isotopic (rel $D'$) to
$\gamma_{t_{m+1}}$.  The original arc $\sigma$ determines an embedded
arc $\sigma_{t_k}\subseteq L_{t_k}$.  Now let $H''$ be the component
of $H'$ which intersects $K_{t_{k+1}}$. Then
$H_{t_{k+1}}=K_{t_{k+1}}\cup H''$ is a (spotted) handlebody. We first
note that $\gamma_{t_{k+1}}$ is isotopic (rel $D'$) in $H_{t_{k+1}}$
to $\sigma$. Therefore $\sigma_{t_k}$ is isotopic in $H_{t_{k+1}}$ to
$\gamma_{t_{k+1}}$.  But since $\sigma_{t_k}$ is a subgraph of
$\Gamma_{t_k}$, we can apply Lemma \ref{GraphToSpinalLemma} to obtain
a subspinal pair $(\breve P_{t_k}, \breve{\mathcal{D}}_{t_k})$ (a disk
with two spots) which is parallel in  $H_{t_{k+1}}$ to
$(P_{t_{k+1}},\mathcal{D}_{t_{k+1}})$.

The parallelism described above determines $\gamma_{t_k}\subseteq
P_{t_k}$. It also determines a rectangle
$(R_{t_k},\gamma_{t_k},\gamma_{t_{k+1}})$ which may be holed. In this
case the holes come from intersections with $H'_{t_k}$, which is the
spotted handlebody corresponding to the spotted product component
$K'_{t_k}\subseteq K_{t_{k+1}}$. As usual, these holes may be adjusted
to be essential discs either in $Q_{t_k}$ or in $H'_{t_k}$, each of
these last consisting of a disc in $H'$ and a vertical annulus in
$K'_{t_k}$.

Now enlarge $(\breve P_{t_k}, \breve{\mathcal{D}}_{t_k})$ to a \spinal
for $K_{t_k}$. We use Lemma \ref{SpinalPairLemma} to obtain \spinals
$(P_{t_i},\mathcal{D}_{t_i})$ for each $t_i<t_k$. We now proceed as in
Lemma \ref{L:half-disc}. Recalling, consider a path
$\gamma_{t_0}=\gamma_0\subseteq P_0$ connecting $D^0$ and $D^1$.  Now
use parallelism to go ``upwards'', obtaining $\gamma_{t_i}\subseteq
P_{t_i}$ for any $t_i<t_k$. In the same way $\gamma_{t_{k-1}}$
determines through parallelism an arc $\gamma_{t_k}'$ in
$P_{t_k}$. Now $\gamma_{t_k}$, $\gamma_{t_k}'\subseteq P_{t_k}$ both
connect the same spots of $K_{t_k}$ therefore they are isotopic (rel
$\partial\mathcal{D}_{t_k}$) in $P_{t_k}$. We adjust them to coincide,
obtaining $\gamma_{t_k}'=\gamma_{t_k}$.

To finalize, let $\Delta$ be the union of the rectangles $R_{t_i}$'s
and the half-disc $\Delta_{t_m}$, which is the desired standard
half-disc in this last case.
\end{proof}

\begin{remark}\label{UnholedHalf-Disc2}
Note from the argument that, as in the sphere body case, a rectangle
$R_{t_i}$ is holed only when $K_{t_i}\subseteq K_{t_{i+1}}$ is in a
situation like Case 2 of Figure \ref{Aut3Splitting}.
\end{remark}

The next step is to describe the unlinking move in this setting of
mixed bodies. Suppose that there is a half-disc $\Delta$ for $Q_0$ in
$Q_1$, which we assume to be in standard position with respect to a
complete sequence of regular values $\{t_i\}$. The goal is to remove
holes of the rectangles $R_{t_i}$ and the half-disc
$\Delta_{t_m}$. Let $s<t$ be consecutive terms of the sequence and
consider $R_s\subseteq K_t$. If $K_t$ is a spotted ball then the
operation works precisely as in Section \ref{TrainTrackSection}.

The new situation appears when $K_t$ is a spotted product, and even
there the argument is essentially the same, so we will be brief.
Recall from Section \ref{TrainTrackSection} that the main condition
for the unlinking to be possible was that one side of the holed
rectangle $R_{t_i}$ needed to intersect a spot of $K_{t_{i+1}}$ in a
spot containing a single spot $D$ of $K_{t_i}$. If that is the case,
as before, we consider the sphere $S$ in the mixed body $M$ which is
dual to $D$. The holes of $R_{t_i}$ are removed by ``isotopying them''
along $R_{t_i}$ close to $S$ and then along $S-\intr(D)$. Here the
only new feature is that a hole may correspond to a disc $E\subseteq
Q_{t_i}\cup H'$.

Let $N$ be a neighborhood (in $Q_{t_i}\cup H'$) of a component $E$ of
$R_{t_i}\cap (Q_{t_i}\cup H')$. One can consider it as a neighborhood
of an arc $\alpha$ either in $H'$ (in case $E$ is a disc intersecting
$H'$) or in $Q_{t_i}$ (otherwise). Isotope $\alpha$ (rel
$\partial\alpha$) along $R_{t_i}$ bring it close to $D$. Now perform
the isotopy along $S-\intr(D)$ as before. Note that this isotopy does
not introduce intersections of $\alpha$ with $R_{t_i}$, for $S$
intersects $R_{t_i}$ only at its boundary.  Call the composition of
these isotopies $\iota$ and let $\alpha'=\iota(\alpha)\subseteq
K_{t_i}$. Extend $\iota$ to $N$ and let $N'=\iota(N)\subseteq
K_{t_i}$, a neighborhood of $\alpha'$.

We separate two cases. If $\alpha\subseteq Q_{t_i}$ the unlinking
move is done, as in Section \ref{TrainTrackSection}. If
$\alpha\subseteq H'$ note that $\alpha$ and $\alpha'$ are isotopic
(rel $\partial\alpha$) as arcs in $K_{t_{i+1}}\cup H'$. Consider
an ambient isotopy $\kappa$ of $K_{t_{i+1}}\cup H'$ such that
$\kappa(\alpha')=\alpha$ and fixing $\partial K_{t_{i+1}}\cup H'$.
Adjust $\kappa$ so that $\kappa\circ\iota(N)=N$. Now replace
$K_{t_i}$ by $\kappa(K_{t_i})$, so $H'$ is preserved. This
completes the description of the operation, which we call an {\em
unlinking along $R_{t_i}$ and $S$}. Clearly such an unlinking
reduces the number of holes in $R_{t_i}$.

\begin{proposition}\label{P:unlinking2}
Let $0=t_0<t_1<\dots<t_{k-1}<t_k=1$ be a complete sequence of regular
values. Suppose that $\Delta$ is a half-disc for $Q_0\subseteq Q_1$,
standard with respect to $\{\,t_i\,\}$. Then there exists a sequence
of unlinkings taking $Q_0$ to $Q_0'\subseteq Q_1$ for which $\Delta$
is an unholed half-disc.  Moreover
$Q_0'\cap\mathcal{E}_1=Q_0\cap\mathcal{E}_1$.
\end{proposition}
\begin{proof}[Sketch of proof]
The proof is essentially the same as that of Proposition
\ref{P:unlinking}. We just note that also here, as Section
\ref{TrainTrackSection}, a rectangle $R_{t_i}$ will be holed only if
$K_{t_i}\subseteq K_{t_{i+1}}$ is in a situation like Case 2 of Figure
\ref{Aut3Splitting} (see Remark \ref{UnholedHalf-Disc2}
above). Therefore one can always find a spot of $K_{t_{i+1}}$
containing a single spot of $K_{t_i}$ and containing a side of
$R_{t_i}$. A finite sequence of unlinkings in $K_{t_{i+1}}$ yields an
unholed rectangle $R_{t_i}$. As before, the final unholed half-disc
$\Delta$ is obtained by induction on $i$. The other conclusion of the
theorem clearly also holds here.
\end{proof}

The last ingredient needed for the proof of Theorem
\ref{CompressionTrainTrackTheorem} is concentricity. The following is
the analogue of Lemma \ref{L:concentric} needed in this setting.

\begin{lemma}\label{L:concentric2}
In the situation described in Proposition \ref{P:unlinking2} we have
that $Q_0'\subseteq Q_1$ are concentric.
\end{lemma}
\begin{proof}[Sketch of proof]
As in Section \ref{TrainTrackSection} the rectangles $R_{t_i}$ in the
original half-disc $\Delta$ will be holed only in a situation like
Case 2 of Figure \ref{Aut3Splitting} (see Remark
\ref{UnholedHalf-Disc2}). Therefore, as in Lemma \ref{L:concentric},
an unlinking does not interfere with the pinching half-discs. The
argument is completed as before.
\end{proof}

\begin{proof}[Proof of Theorem \ref{CompressionTrainTrackTheorem}]

As usual, the proof is analogous to that of Theorem
\ref{TrainTrackTheorem} (which was done as the proof of
Theorem~\ref{BackTrackTheorem}).

Consider $f\in\mathcal{A}_f$. If $g=f|_{Q'}$ is reducible or periodic,
the proof is over, so assume that $g$ is generic, with growth rate
$\lambda$. If it is not train-track generic there is an integer $m<1$
such that $Q_m$ back-tracks in $Q_1$ (Lemma \ref{L:back-track}). Let
$K_m\subseteq K_1$ be the spotted ball in a $K_1$ (spotted ball or
product) as in the lemma and assume that $m$ is the greatest such
value. Let $D_m^0$, $D_m^1$ be the two spots of $K_m$ as in the lemma,
contained in a single spot $D$ of $K_1$. Consider the increasing
sequence of consecutive integers $m=i_0, i_1,\dots,i_{1-m}=1$. Then
there are two spots $D^0$, $D^1$ of $K_{i_{(2-m)}}=K_0$ containing
$D_m^0$, $D_m^1$ respectively (hence $D^0$, $D^1\subseteq D$). The
parallelism in the conclusion of Lemma \ref{L:back-track} gives an arc
$\sigma\subseteq\partial K_0$ as in the statement of Proposition
\ref{P:half-disc2}, yielding a half-disc $\Delta$ for $Q_0$ in $Q_1$
connecting $D^0$ and $D^1$.

Choose a complete sequence $0=t_0<\dots<t_k=1$ of regular values and
assume that $\Delta$ is in standard position with respect to
$\{t_i\}$. Apply Proposition \ref{P:unlinking2} to obtain $\Delta$
unholed. Lemma \ref{L:concentric2} assures that an unlinking preserves
$Q'$, $H'$. This operation does not change the growth rate $\lambda$
of $g$.

Skipping details, one proceeds as in \cite{UO:Autos} and Theorem
\ref{TrainTrackTheorem}, using $\Delta$ to realize down-moves or
diversions, respective analogues of folds and pulling tights, as in
\cite{BH:Tracks} (note that ``splittings'' may be needed as a
preparation for these moves). After performing a down-move the
resulting automorphism $g'$ may be either periodic, reducible or
generic. In the first two cases the proof is done, so assume it is
generic. The operation increases the greatest value $m$ in the
beginning of the proof, therefore repetition of this process has to
stop eventually (when $m=0$).

After performing a diversion, which is the case when $m=0$, the
resulting $g'$ may also be periodic, reducible or generic. Again
assuming it is generic generic its growth $\lambda'$ is strictly
smaller than $\lambda$, therefore repetition of this process also has
to stop.

The procedure above, which can be carried on for as long as $g'$ is
generic admitting back-tracking, yields $g'$ either periodic,
reducible or generic without back-tracking in $H$, completing the
proof for $g'$.

We recall that the conclusions on $f'|_{Q}$, $f'|_{H'}$ follow from
Proposition \ref{FurtherReducingSurface}.
\end{proof}

\hop

\hop

\bibliographystyle{hplain} \bibliography{ReferencesUO3}

\begin{thebibliography}{10}

\bibitem{BH:Tracks}
M.~Bestvina and M.~Handel.
\newblock Train tracks and automorphisms of free groups.
\newblock {\em Ann. Math.}, 135:1--51, 1992.

\bibitem{FB:CompressionBody}
F.~Bonahon.
\newblock Cobordism of automorphisms of surfaces.
\newblock {\em Ann. Sci. \'Ecole Norm. Sup. (4)}, 16(2):237--270, 1983.

\bibitem{LNC:Generic}
L.~Navarro de~Carvalho.
\newblock Generic automorphisms of handlebodies.
\newblock Dissertation, Rutgers University, 2003.

\bibitem{LNC:Tightness}
L.~Navarro de~Carvalho.
\newblock Tightness and efficiency of irreducible automorphisms of
  handlebodies.
\newblock preprint, available at arXiv:math.GT/0408353, 2004.

\bibitem{EC:Automorphisms}
E.~C{\'e}sar de~S{\'a}.
\newblock Automorphisms of 3-manifolds and representations of 4-manifolds.
\newblock Dissertation, Warwick University, 1977.

\bibitem{FL:ConnectedS2TimesS1}
Albert Fathi and Fran{\c{c}}ois Laudenbach.
\newblock Diff\'eomorphismes pseudo-{A}nosov et d\'ecomposition de {H}eegaard.
\newblock {\em C. R. Acad. Sci. Paris S\'er. A-B}, 291(6):A423--A425, 1980.

\bibitem{HG:64}
H.~B. Griffiths.
\newblock Automorphisms of a {$3$}-dimensional handlebody.
\newblock {\em Abh. Math. Sem. Univ. Hamburg}, 26:191--210, 1963/1964.

\bibitem{AH:HomologicalStability}
Allen Hatcher.
\newblock Homological stability for automorphism groups of free groups.
\newblock {\em Comment. Math. Helv.}, 70(1):39--62, 1995.

\bibitem{WJPS:Characteristic}
W.~Jaco and P.~Shalen.
\newblock Seifert fibered spaces in 3-manifolds.
\newblock {\em Mem. Amer. Math. Soc}, 220:1--192, 1979.

\bibitem{JWW:Homeomorphisms}
B.~Jiang, S.~Wang, and Y.~Wu.
\newblock Homeomorphisms of 3-manifolds and the realization of {N}ielsen
  number.
\newblock {\em Comm. Anal. Geom.}, 9(4):825--877, 2001.

\bibitem{KJ:Characteristic}
K.~Johannson.
\newblock {\em Homotopy equivalences of 3-manifolds with boundaries}.
\newblock Number 761 in LNM. Springer Verlag, 1979.

\bibitem{FL:Spheres}
F.~Laudenbach.
\newblock Sur les 2-sph{\`e}res d'une vari{\'e}t{\'e} de dimension 3.
\newblock {\em Ann. of Math.}, 97(2):57--81, 1973.

\bibitem{FL:HomotopyIsotopy}
Fran{\c{c}}ois Laudenbach.
\newblock {\em Topologie de la dimension trois: homotopie et isotopie}.
\newblock Soci\'et\'e Math\'ematique de France, Paris, 1974.
\newblock With an English summary and table of contents, Ast\'erisque, No. 12.

\bibitem{DM:Elliptic}
D.~McCullough.
\newblock Isometries of elliptic 3-manifolds.
\newblock available ArXiv math.GT/0010077.

\bibitem{DM:MappingSurvey}
D.~McCullough.
\newblock Mappings of reducible 3-manifolds.
\newblock {\em Geometric and Algebraic Topology, Banach Center Publications},
  18:61--76, 1986.

\bibitem{UO:Autos}
U.~Oertel.
\newblock Automorphisms of 3-dimensional handlebodies.
\newblock {\em Topology}, 41:363--410, 2002.

\bibitem{PS:Geometries}
P.~Scott.
\newblock The geometries of $3$-manifolds.
\newblock {\em Bull. London Math. Soc.}, 15:401--487, 1983.

\bibitem{WPT:Notes}
W.~Thurston.
\newblock The geometry and topology of 3-manifolds.
\newblock notes, 1978.

\bibitem{FW:Heegaard}
Friedhelm Waldhausen.
\newblock Heegaard-{Z}erlegungen der {$3$}-{S}ph\"are.
\newblock {\em Topology}, 7:195--203, 1968.

\end{thebibliography}

\end{document}